\documentclass[12pt,english]{article}
\usepackage[T1]{fontenc}
\usepackage[latin9]{inputenc}
\usepackage{geometry}
\usepackage{color}
\usepackage{babel}
\usepackage{float}
\usepackage{mathtools}
\usepackage{enumitem}
\usepackage{amsmath}
\usepackage{amssymb}
\usepackage{cancel}
\usepackage{stackrel}
\usepackage{graphicx}
\usepackage{wasysym}
\usepackage{esint}
\usepackage[unicode=true,pdfusetitle,
 bookmarks=true,bookmarksnumbered=false,bookmarksopen=false,
 breaklinks=false,pdfborder={0 0 0},pdfborderstyle={},backref=false,colorlinks=true]
 {hyperref}

\makeatletter

\providecommand{\tabularnewline}{\\}

\@ifundefined{date}{}{\date{}}
\usepackage{lmodern}
\usepackage[T1]{fontenc}
\usepackage{caption}
\numberwithin{equation}{section}
\usepackage{upgreek}
\usepackage{titlesec}
\titlespacing*{\section}
{0pt}{5.5ex plus 1ex minus .2ex}{4.3ex plus .2ex}
\titlespacing*{\subsection}
{0pt}{5.5ex plus 1ex minus .2ex}{4.3ex plus .2ex}
\titlespacing*{\subsubsection}
{0pt}{5.5ex plus 1ex minus .2ex}{4.3ex plus .2ex}
\titlespacing*{\paragraph}
{0pt}{5.5ex plus 1ex minus .2ex}{4.3ex plus .2ex}

\@ifundefined{showcaptionsetup}{}{%
 \PassOptionsToPackage{caption=false}{subfig}}
\usepackage{subfig}
\makeatother

\begin{document}
\title{A globally conservative finite element MHD code and its application
to the study of compact torus formation, levitation and magnetic compression }
\author{Carl Dunlea$^{1,2,(a)}$, and  Ivan Khalzov$^{2}$, General Fusion
Team, Akira Hirose}

\maketitle
$^{1}$University of Saskatchewan, 116 Science Pl, Saskatoon, SK S7N
5E2, Canada 

$^{2}$General Fusion, 106 - 3680 Bonneville Pl, Vancouver, BC V3N
4T5, Canada

$^{(a)}$Current address: Tokamak Energy Ltd., Milton Park, Oxfordshire,
OX14 4SD, England. Author to whom correspondence should be addressed:
cpd716@mail.usask.ca
\begin{abstract}
The DELiTE (Differential Equations on Linear Triangular Elements)
framework was developed for spatial discretisation of partial differential
equations on an unstructured triangular grid in axisymmetric geometry.
The framework is based on discrete differential operators in matrix
form, which are derived using linear finite elements and mimic some
of the properties of their continuous counterparts. A single-fluid
two-temperature MHD model is implemented in this framework. The inherent
properties of the operators are used in the code to ensure global
conservation of energy, particle count, toroidal flux, and (in some
scenarios) angular momentum. The code was applied to study a novel
experiment in which a compact torus (CT), produced with a magnetized
Marshall gun, is magnetically levitated off an insulating wall and
then magnetically compressed through the action of currents in the
levitation/compression coils located outside the wall. We present
numerical models for CT formation, levitation, and magnetic compression,
and comparisons between simulated and experimental diagnostics.\\
\\
\end{abstract}

\section{Introduction}

Ideally, a numerical discretisation method should, as well as accurately
representing the continuous form of the mathematical equations that
describe a particular physical system, reproduce the physical properties
of the system being modelled. In practice, such properties are often
expressed as conservation laws, and their maintenance can be just
as important as standard numerical method assessment gauges like convergence,
stability, accuracy, and range of applicability. Numerical solutions
that contradict basic physical principles by, for example, destroying
mass or energy are inherently unreliable when applied to novel physical
regimes. On the other hand, many numerical methods do not have strong
conservation properties, or conserve only some naturally conserved
quantities, but can still regularly produce informative results that
have a resemblance to real-world observations. However, as more complicated
physical systems are modelled, physically incorrect numerical solutions
can go unnoticed. This is especially true in simulating complex physical
models, such as fluid dynamics or magnetohydrodynamics (MHD), where
the maintenance of conservation laws puts more constraints on numerical
schemes and, thus, helps to avoid spurious solutions. It is best to
deal with the physical fidelity of the model at the numerical method
design level \cite{Perot2}, and try to replicate, in the discretised
form, as many of the conservation laws inherent to the original physical
system as possible. 

Numerical methods with discrete conservation properties are well known
in computational fluid dynamics (for example, \cite{Perot2,Caramana,Margolin,Zhang}),
and in computational MHD \cite{Derigs,Liu,Franck}. In this paper
we present a novel numerical scheme for a two-dimensional (axisymmetric)
compressible MHD system, which is based on a continuous Galerkin linear
finite element method applied to an unstructured triangular mesh.
By construction, the scheme has global (for the whole domain) conservation
of mass, energy, toroidal flux and angular momentum. A novelty of
the code is that all discrete spatial differential operators are represented
as matrices, and the discretized MHD equations are obtained by simply
replacing the original continuous differentiations with the corresponding
matrix operators. 

Note that by global conservation of a quantity in our numerical method,
we imply that there is a discretised analogue of the continuity equation
for that quantity, and, when integrated over the volume, its fluxes
are completely cancelled in the interior of the domain, even though
the explicit form of these fluxes are not always given. As shown in
\cite{Hughes} and \cite{Perot2}, global conservation for a method
with local support ($i.e.,$ local stencil) also implies local conservation.
To enable the development of a numerical formulation with the aforementioned
conservation properties, the discrete differential operators must
obey a property equivalent, for a scalar field $u(\mathbf{r})$ and
a vector field $\mathbf{p}(\mathbf{r})$, to 
\begin{equation}
\int u\nabla\cdot\mathbf{p}\,dV+\int\mathbf{p}\cdot\nabla u\,dV=\int\nabla\cdot\left(u\mathbf{p}\right)\,dV=\int u\mathbf{p}\cdot d\boldsymbol{\Gamma}\label{eq:500.0}
\end{equation}
so that the discrete forms of the differential product rule and divergence
theorem are satisfied \cite{Perot}. This is the essence of mimetic
schemes, also known as support operator methods, to which this work
is closely related \cite{Perot2,Shashkov}.

In section \ref{sec:Finite-element-discretisation}, we use the finite
element method to derive various discrete differential matrix operators
for a two-dimensional domain with axial symmetry in cylindrical coordinates,
and discuss their properties. These operators satisfy the discrete
form of equation \ref{eq:500.0}, and constitute the base of the DELiTE
(Differential Equations on Linear Triangular Elements) framework,
which is implemented in MATLAB. In section \ref{sec:Discrete-form-of},
the discrete single-fluid two-temperature MHD equations are constructed
within the framework. The formal process of equation development is
detailed in \cite{thesis}. Here, using inherent properties of the
differential operators, we demonstrate global conservation of mass,
toroidal flux, angular momentum, and energy. In section \ref{sec:Experiment-overview},
an overview of the magnetic compression experiment is presented -
more details can be found in \cite{thesis,exppaper}. The magnetic
compression experiment at General Fusion was a repetitive non-destructive
test to study plasma physics applicable to magnetic target fusion
compression. A compact torus (CT) is formed with a co-axial gun into
a containment region with an hour-glass shaped inner flux conserver,
and an insulating outer wall. External coil currents hold the CT off
the outer wall (radial \textquotedbl levitation\textquotedbl ) and
then rapidly compress it inwards. In section \ref{sec:Model-for-CT},
we present the details of our numerical models for simulating CT formation,
magnetic levitation, and compression within the DELiTE framework.
The model includes coupling between the poloidal external vacuum field
solution in the part of the domain representing the insulating wall,
and the full MHD solutions in the remainder of the domain. In section
\ref{sec:Simulation-results,-and}, we present principal code input
parameters and simulation results, and simulated diagnostics are compared
with the experimentally measured counterparts. With inclusion of the
insulating wall in the model, the effect of reduced plasma/wall interaction
with an improved levitation/compression coil configuration, as observed
in the experiment, is reproduced with MHD simulations. In the conclusion,
section \ref{sec:Conclusion}, we present a summary of principal code
features and suggest possible further extensions. A detailed derivation
of the first order node-to-element differential operator is presented
in Appendix \ref{sec:Appendix-A}. Further node-to element and node-to-node
operator identities are demonstrated in Appendices \ref{sec:Appendix-B}
and \ref{sec:Appendix-C} respectively. In Appendix \ref{sec:Appendix-D},
it is shown how energy and momentum conservation can be maintained
when artificial density diffusion is included in the model. 

\section{Finite element discretisation and differential operators\label{sec:Finite-element-discretisation}}

In this section, we start with a brief description of the discretisation
and finite element method, and introduce notations, then we derive
several types of discrete differential operators used in the code
and discuss their properties. 

\subsection{Discretisation and finite element method overview\label{subsec:Finite-element-method}}

The two-dimensional computational domain, with azimuthal symmetry
in cylindrical coordinates, is represented by an unstructured triangular
mesh. To develop the finite element discretisation, we drew inspiration
from material presented in \cite{PICwebsite}, which in turn, is based
on material in \cite{Strang}, in which a finite element method is
used to solve Laplace and Poisson equations in two dimensions. The
freely-available mesh generator DISTMESH \cite{MeshPersson1,MeshPersson2,MeshPersson3}
was adapted to provide the computational grid. Nodes are located at
triangle vertices. In the linear finite element method, any continuous
function $u(\mathbf{r})$ is approximated as a piecewise continuous
function $U(\mathbf{r})$ that is linear across each triangular element:
\begin{equation}
u(\mathbf{r})\approx U(\mathbf{r})=\overset{N_{e}}{\underset{e=1}{\Sigma}}\,U^{e}(\mathbf{r})\label{eq:500.01}
\end{equation}
where $N_{e}$ is the number of elements, and $U^{e}(\mathbf{r})$
is the contribution to $U(\mathbf{r})$ from within element $e$:
\begin{equation}
U^{e}(\textbf{\ensuremath{\mathbf{r}}})=A^{e}+B^{e}r+C^{e}z\label{eq:500.02}
\end{equation}
Here, $A^{e},\,B^{e}$ and $C^{e}$ are constants that are specific
to element $e$. The mechanics of the derivations of these coefficients
is presented briefly in in appendix \ref{sec:detailed-derivation-of},
and an overview of the matrix operator asssembly is shown. Equivalent
to equation \ref{eq:500.01}, $U(\mathbf{r})$ may be defined using
a combination of basis functions as 
\begin{equation}
u(\mathbf{r})\approx U(\mathbf{r})=\overset{N_{n}}{\underset{n=1}{\Sigma}}\,U_{n}\,\phi_{n}(\mathbf{r})\label{eq:500.03}
\end{equation}
where $N_{n}$ is the number of nodes, $U_{n}$ represents the values
of $U(\mathbf{r})$ at node $n$, and $\phi_{n}(\mathbf{r})\equiv\phi_{n}(r,\,z)$
is the basis function associated with node $n$. 
\begin{figure}[H]
\centering{}\includegraphics[width=8cm,height=6cm]{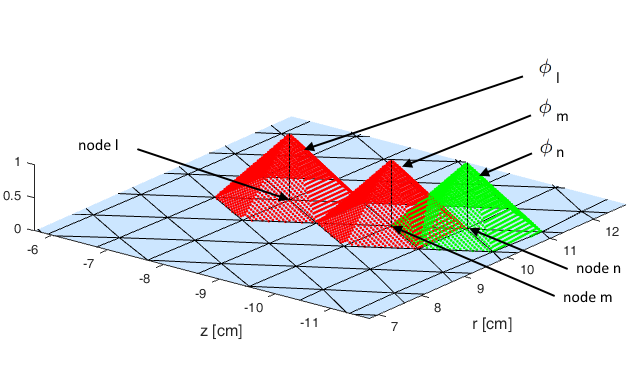}\caption{\label{fig:Linear-basis-function}Linear basis function depiction
for triangular elements}
\end{figure}
The linear basis functions have the forms of pyramids, as depicted
in figure \ref{fig:Linear-basis-function}, which indicates a portion
of a computation mesh, with depictions of the basis functions for
example nodes $l,\,m$ and $n$. Each node $n$ is associated with
a pyramid function $\phi_{n}$, which has an elevation of one above
the node, and falls linearly to zero at the immediately surrounding
nodes. Each pyramid function $\phi_{n}$ has $e_{n}$ sides, where
$e_{n}$ is the node-specific number of triangular elements surrounding
node $n$. Note that the pyramid shapes associated with basis functions
for nodes on the domain boundary have, in addition, one or more vertical
sides at the boundaries. The remaining individual pyramid-side functions
are defined by sections of planes that are tilted relative to the
$r-z$ plane. Thus, each node is associated with $e_{n}$ tilted planes.
In total, there are $K$ tilted planes defined in the solution domain,
where $K=\overset{N_{n}}{\underset{n=1}{\Sigma}}e_{n}=3N_{e}$. Each
triangular element is associated with three tilted planes. The tilted
planes are truncated to constitute pyramid sides by defining the functions
representing the pyramid sides as 
\begin{equation}
\psi_{n}^{e}(\textbf{\ensuremath{\mathbf{r}}})=a_{n}^{e}+b_{n}^{e}r+c_{n}^{e}z\label{eq:502.2}
\end{equation}
with the additional truncating property that $\psi_{n}^{e}(\textbf{\ensuremath{\mathbf{r}}})=0$
for all points located at $\mathbf{r}$ which lie outside the triangular
element associated with $\psi_{n}^{e}$. The notation $\psi_{n}^{e}$
indicates that each pyramid side function is associated with a particular
node $n$, \emph{and} with a particular triangular element $e$. The
coefficients $a_{n}^{e},\,b_{n}^{e}$, and $c_{n}^{e}$ are also associated
with a particular node and element, and are such that $\psi_{n}^{e}=1$
at its associated node $n$, and $\psi_{n}^{e}=0$ at the other two
nodes in the triangular element $e$. In summary, the basis functions
$\phi_{n}(\textbf{\ensuremath{\mathbf{r}}})$, which define pyramid
shapes with a peak elevation of one at node $n$, can be expressed
as 
\begin{equation}
\phi_{n}(\textbf{\ensuremath{\mathbf{r}}})=\overset{e_{n}}{\underset{}{\Sigma}}\,\psi_{n}^{e}(\textbf{\ensuremath{\mathbf{r}}})\label{eq:502.3}
\end{equation}
where the summation is over the pyramid side functions associated
with node $n$, each of which is non-zero only over its associated
triangular element. The basis functions have the property that $\phi_{n}(r_{j},\,z_{j})=\delta_{nj}$,
where $\delta_{nj}$ is the Kronecker delta, and $(r_{j},\,z_{j})$
are the coordinates of node $j$. Noting that the volume of a pyramid
is given by $V_{pyramid}=sH/3$, where $s$ is the pyramid base area
and $H$ is the pyramid height, this leads naturally, for any continuous
function $U(\mathbf{r})$ (including the piecewise linear approximation),
to the property 
\begin{equation}
\int\phi_{n}(\mathbf{r})\,U(\mathbf{r})\,dr\,dz=(U_{n}+\mathcal{O}(h_{e}))\,\int\phi_{n}(\mathbf{r})\,dr\,dz\approx\frac{U_{n}\,s_{n}}{3}\label{eq:502.4}
\end{equation}
where $U_{n}=U(r_{n},z_{n})$, $s_{n}$ is the support area of node
$n$ (area of the base of the pyramid function defined by $\phi_{n}(\mathbf{r})$),
pyramid height $H=1,$ and $h_{e}$ is the typical element size. This
identity is analogous to the integral property of the Dirac delta
function. In deriving the property, it is assumed that the function
$U(\mathbf{r})$ is sufficiently smooth that it is approximately constant
(to order $h_{e}$) in the support area of node $n$. We neglect the
term of order $h_{e}$ because our numerical scheme has overall the
first order accuracy, as defined by the use of linear basis functions. 

Another important property of the basis functions is related to partition
of unity - the sum of all the basis functions in the domain, at any
point in the domain ($i.e.,$ at non-nodal locations, as well as at
nodal locations), is equal to one. This property also hold for the
pyramid side functions, $i.e.,$ $\overset{N_{n}}{\underset{n=1}{\Sigma}}\phi_{n}(\textbf{\ensuremath{\mathbf{r}}})=\overset{N_{n}}{\underset{n=1}{\Sigma}}\psi_{n}(\textbf{\ensuremath{\mathbf{r}}})=1$. 

\subsection{Notations for matrices and operations }

The notation $\overline{X}$ (or $\overline{x})$ is used to denote
vectors of dimensions $[N_{n}\times1]$ that contain node-associated
quantities, while $\widehat{X}\,[N_{e}\times1]$ denotes vectors of
element-associated quantities. The notation $\overline{\widehat{X}}$will
be used to denote matrices of dimensions $[N_{e}\times N_{n}]$ that
operate on vectors of nodal quantities $\overline{X}$, to produce
vectors of elemental quantities $\widehat{Y}$, while $\widehat{\overline{X}}\,[N_{n}\times N_{e}]$
operates on vectors of elemental quantities $\widehat{X}$ to produce
vectors of nodal quantities $\overline{Y}$. Finally, $\overline{\overline{X}}$
and $\widehat{\widehat{X}}$ denote square matrices with dimensions
$[N_{n}\times N_{n}]$ and $[N_{e}\times N_{e}]$ respectively. Notations
defining the various matrix dimensions are collected in table \ref{tab:Notations-for-matrices}.
\begin{table}[H]
\begin{centering}
\includegraphics[scale=0.15]{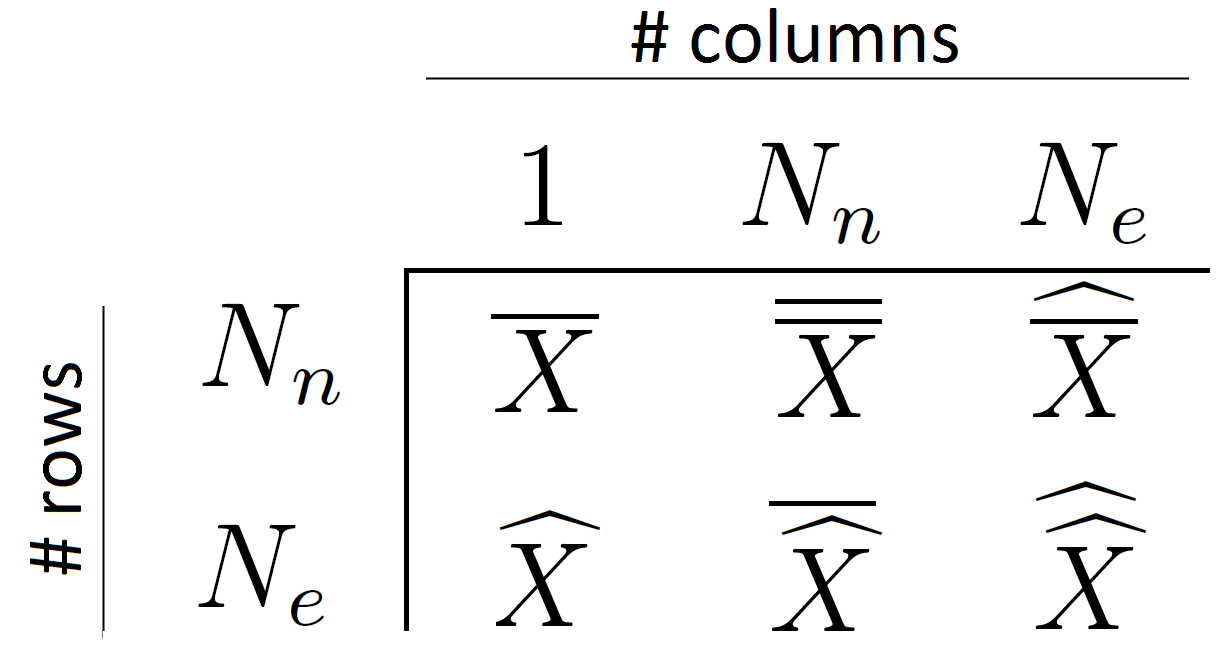}\caption{\label{tab:Notations-for-matrices}Notations for matrices of various
dimensions}
\par\end{centering}
\end{table}

Defining operation notations, a space between two arrays of equal
size represents the Hadamard product, piecewise element-by element
multiplication ($e.g.,$ $(\overline{a}\,\,\overline{b})_{i}=a_{i}b_{i}$),
and the symbol $\,/\,$ represents Hadamard division, piecewise element-by
element division ($e.g.,$ $(\overline{a}\,/\,\overline{b})_{i}=a_{i}/b_{i}$).
The symbol $*$ between two matrices ($e.g.,$ $\overline{\overline{C}}=\overline{\overline{A}}\,*\,\overline{\overline{B}}$)
implies regular matrix multiplication. The superscript $T$ implies
the transpose operation, $e.g.,$ $\overline{\overline{B}}=\overline{\overline{A}}^{T}$,
and the superscript $-1$ implies the inverse operation. 

\subsection{Element-node mapping operators }

The connectivity operator $\overline{\widehat{M}}$, used to map between
nodal and elemental quantities, has dimensions $[N_{e}\times N_{n}]$.
Each row of $\overline{\widehat{M}}$ corresponds to a particular
triangular element, and has only three non-zero entries. $\overline{\widehat{M}}(e,\,n)=1$
for the column indexes $n$ corresponding to the indexes $n$ of the
nodes located at the vertices of the triangle with index $e.$ We
introduce node-to-element averaging as
\begin{equation}
<\overline{X}>^{e}=\frac{1}{3}\overline{\widehat{M}}*\overline{X}\label{eq:502.31}
\end{equation}
For example, the radial coordinates of the element centroids are defined
with $\widehat{r}=(r_{1}^{e},\,r_{2}^{e},...r_{N_{e}}^{e})^{T}=<\overline{r}>^{e}$,
where $\overline{r}=(r_{1},r_{2},...r_{N_{n}})^{T}$ contains the
radial coordinates of the nodes. Similarly, the axial centroid coordinates
are defined as $\widehat{z}=<\overline{z}>^{e}$. The vector of nodal
support areas $\overline{s}$, containing the areas of the bases of
the pyramid functions associated with the nodes, is related to $\widehat{s}$,
the vector of elemental areas, as $\overline{s}=\overline{\widehat{M}}^{T}*\widehat{s}$. 

Volume integrals over the computational domain can be approximated
in two ways, corresponding to nodal or elemental representations of
the integrand function 
\begin{align*}
\int u(\mathbf{r})\,dV & \approx\overline{dV}^{T}*\overline{U}\,\,\,\mbox{ or }\,\,\,\int u(\mathbf{r})\,dV\approx\widehat{dV}^{T}*\widehat{U}
\end{align*}
where $\overline{dV}=\frac{2\pi}{3}\,\overline{s}\,\,\overline{r}$
contains the volumes associated with each node, which are found by
integrating the node-associated areas from $0$ to $2\pi$ in the
toroidal direction, and $\widehat{dV}=2\pi\,\widehat{s}\,\,\widehat{r}$
contains the elemental volumes. The factor of three in the former
expression arises because each elemental area is shared by three nodes.
Note that these two approximations do not give identical results for
vectors related by equation \ref{eq:502.31}, unless the original
integrand function $u(\mathbf{r})$ is constant. \\

Defining $\overline{\overline{R}}$, $\overline{\overline{S}}$, $\widehat{\widehat{R}}$
and $\widehat{\widehat{S}}$ as the diagonal arrays constructed from
$\overline{r}$, $\overline{s}$, $\widehat{r}$, and $\widehat{s}$,
we define a volume-averaging operator\\
\[
\widehat{\overline{W}}=\overline{\overline{R}}^{-1}*\overline{\overline{S}}^{-1}*\overline{\widehat{M}}^{T}*\widehat{\widehat{S}}*\widehat{\widehat{R}}
\]
 that is used to map element-based quantities to node-based quantities,
as
\begin{equation}
<\widehat{U}>=\widehat{\overline{W}}*\widehat{U}\label{eq:515.1}
\end{equation}
This operator satisfies the following identity: 
\begin{equation}
\overline{dV}^{T}*\left(\overline{Q}\,\,\left(\widehat{\overline{W}}*\widehat{U}\right)\right)=\widehat{dV}^{T}*\left(\widehat{Q}\,\,\widehat{U}\right)\label{eq:516}
\end{equation}
where $\widehat{U}$ and $\widehat{Q}$ are the discrete representations,
defined at the element centroids, of the approximations to the continuous
functions $u(\mathbf{r})$ and $q(\mathbf{r})$, and $\overline{Q}$
is the discrete representation, defined at the nodes, of $q(\mathbf{r})$,
and is related to $\widehat{Q}$ by equation \ref{eq:502.31}. A proof
of this identity follows:

\begin{align*}
\overline{dV}^{T}*\left(\overline{Q}\,\,\left(\widehat{\overline{W}}*\widehat{U}\right)\right) & =\frac{2\pi}{3}\left(\overline{s}\,\,\overline{r}\right)^{T}*\left(\overline{Q}\,\,\left(\overline{\overline{R}}^{-1}*\overline{\overline{S}}^{-1}*\overline{\widehat{M}}^{T}*\widehat{\widehat{S}}*\widehat{\widehat{R}}*\widehat{U}\right)\right)\\
 & =\frac{2\pi}{3}\overline{Q}^{T}*\left(\overline{\widehat{M}}^{T}*\widehat{\widehat{S}}*\widehat{\widehat{R}}*\widehat{U}\right)\\
 & =2\pi\left(\widehat{\widehat{S}}*\widehat{\widehat{R}}*\widehat{U}\right)^{T}*\left(\frac{1}{3}\overline{\widehat{M}}*\overline{Q}\right) & \mbox{(transpose the scalar)}\\
 & =\widehat{dV}^{T}*\left(\widehat{Q}\,\,\widehat{U}\right) & \mbox{(use equation \ref{eq:502.31})}
\end{align*}
Note that the matrix transpose relation

\begin{equation}
\left(\mathbb{A}*\mathbb{B}\right){}^{T}=\mathbb{B}{}^{T}*\mathbb{A}{}^{T}\label{eq:505.1}
\end{equation}
for matrices $\mathbb{A}$ and $\mathbb{B}$, is used to transpose
the scalar on the right side of the equation in the second last step
of the derivation above. In the particular case with $\overline{Q}=\overline{1}$,
the identity becomes 
\begin{equation}
\overline{dV}^{T}*\left(\widehat{\overline{W}}*\widehat{U}\right)=\widehat{dV}^{T}*\widehat{U}\label{eq:516.1}
\end{equation}

\subsection{First order node-to-element differential operators\label{subsec:Dre Dze}}

\begin{figure}[H]
\centering{}\includegraphics[width=8cm,height=2.5cm]{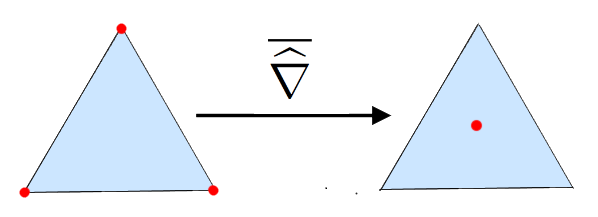}\caption{\label{fig:Drc} Node-to-element differential operator mechanism }
\end{figure}
The node-to-element derivative operator matrices are defined such
that 
\begin{alignat*}{1}
\widehat{U}_{r}' & =\overline{\widehat{Dr}}*\overline{U}\\
\widehat{U}_{z}' & =\overline{\widehat{Dz}}*\overline{U}
\end{alignat*}
Referring for example to the radial derivative operator, application
of $\overline{\widehat{Dr}}$, a matrix with dimensions $[N_{e}\times N_{n}]$,
to the vector $\overline{U}$, returns $\widehat{U}_{r}'\,[N_{e}\times1],$
containing the values for the radial derivatives of $U(\mathbf{r})$
inside the triangular elements. The node-to-element gradient and divergence
operations, for cylindrical coordinates with azimuthal symmetry are
\\
 
\begin{align}
\overline{\widehat{\nabla}}\,\,\overline{U} & =\left(\overline{\widehat{Dr}}*\overline{U}\right)\hat{\mathbf{r}}+\left(\overline{\widehat{Dz}}*\overline{U}\right)\hat{\mathbf{z}}\label{eq:505.01}
\end{align}
and 
\begin{align*}
\overline{\widehat{\nabla}}\cdot\overline{\mathbf{P}} & =\left(\left(\overline{\widehat{Dr}}*\left(\overline{r}\,\,\overline{P}_{r}\right)\right)+\left(\overline{\widehat{Dz}}*\left(\overline{r}\,\,\overline{P}_{z}\right)\right)\right)\,/\,\widehat{r}
\end{align*}
where where $\hat{\mathbf{r}}$ and $\hat{\mathbf{z}}$ are the unit
vectors in the radial and axial directions, and $\overline{P}_{r}$
and $\overline{P}_{z}$ are the nodal representations of the $r$
and $z$ components the continuous vector field $\mathbf{p}(\mathbf{r})$.
A schematic of the node-to-element gradient operation mechanism is
indicated in figure \ref{fig:Drc}. To derive the $\overline{\widehat{Dr}}$
and $\overline{\widehat{Dz}}$ operators we use the elemental representation
of $U(\mathbf{r})$, equations \ref{eq:500.01} and \ref{eq:500.02}.
Equation \ref{eq:500.02} defines the values of the contributions
to $U(\mathbf{r})$ at the nodes of a particular element as 
\begin{equation}
U_{\varepsilon}^{e}=A^{e}+B^{e}r_{\varepsilon}+C^{e}z_{\varepsilon}\label{eq:503-1}
\end{equation}
where $\varepsilon=i,\,j,\,k$ denotes the indexes of the nodes at
the vertices of triangular element $e$. As shown in appendix \ref{sec:detailed-derivation-of},
the expressions for the spatial derivatives of $U(\mathbf{r})$ at
the interior of each element are:
\begin{align}
\frac{\partial U^{e}}{\partial r} & =\underset{\varepsilon}{\Sigma}U_{\varepsilon}^{e}b_{\varepsilon}^{e}\nonumber \\
\frac{\partial U^{e}}{\partial z} & =\underset{\varepsilon}{\Sigma}U_{\varepsilon}^{e}c_{\varepsilon}^{e}\label{eq:505-1}
\end{align}
 where the coefficients $b_{\epsilon}^{e}$ and $c_{\epsilon}^{e}$
are defined in equation \ref{eq:502.2}. $\overline{\widehat{Dr}}$
and $\overline{\widehat{Dz}}$ are initially defined as sparse all-zero
arrays. The values $b_{i}^{e},\,b_{j}^{e}$ and $b_{k}^{e}$ are inserted
in row $e$ of $\overline{\widehat{Dr}}$ with placements at the column
indexes $i,\,j,\,k$. Similarly, the values $c_{i}^{e},\,c_{j}^{e}$
and $c_{k}^{e}$, for each element, are used to assemble $\overline{\widehat{Dz}}$.
The resulting derivative operators produce exact derivatives for nodal
functions with linear dependence on $r$ and $z$, and have first
order accuracy ($i.e.,\,\mathcal{O}(h_{e})$) when applied to nonlinear
functions. The operators introduced in the following subsections are
all based on these node-to-element derivative operators, and so they
all have the same accuracy. \\
\\

\subsection{First order element-to-node differential operators\label{subsec:Drn}}

\begin{figure}[H]
\centering{}\includegraphics[width=6cm,height=2.5cm]{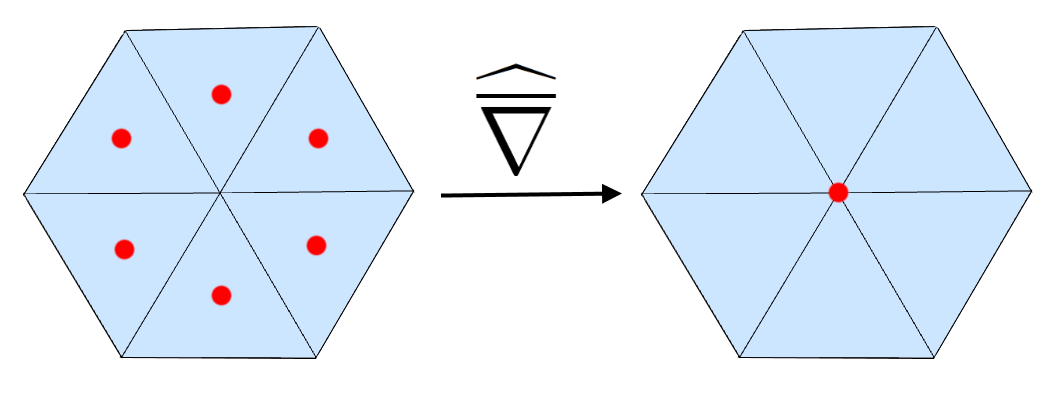}\caption{\label{fig:Drn}Element-to-node differential operator mechanism}
\end{figure}
The element-to-node derivative operator matrices are defined such
that 
\begin{align*}
\overline{U}_{r}' & =\widehat{\overline{Dr}}*\widehat{U}\\
\overline{U}_{z}' & =\widehat{\overline{Dz}}*\widehat{U}
\end{align*}
Referring for example to the radial derivative operator, application
of the $\widehat{\overline{Dr}}$ operator, which is a matrix with
dimensions $[N_{n}\times N_{e}]$, to the vector $\widehat{U}\,[N_{e}\times1]$,
which contains the values of $U(\mathbf{r})$ at the element centroids,
returns $\overline{U}_{r}'\,[N_{n}\times1],$ the values for the radial
derivatives of $U(\mathbf{r})$ at the nodes. The element-to-node
gradient and divergence operations, for cylindrical coordinates with
azimuthal symmetry are 
\begin{align}
\widehat{\overline{\nabla}}\,\,\widehat{U} & =\left(\widehat{\overline{Dr}}*\widehat{U}\right)\hat{\mathbf{r}}+\left(\widehat{\overline{Dz}}*\widehat{U}\right)\hat{\mathbf{z}}\label{eq:515.01}
\end{align}
and 
\begin{align}
\widehat{\overline{\nabla}}\cdot\widehat{\mathbf{P}} & =\left(\left(\widehat{\overline{Dr}}*\left(\widehat{r}\,\,\widehat{P_{r}}\right)\right)+\left(\widehat{\overline{Dz}}*\left(\widehat{r}\,\,\widehat{P_{z}}\right)\right)\right)\,/\,\overline{r}\label{eq:515.02}
\end{align}
A schematic of the element-to-node gradient operation mechanism is
indicated in figure \ref{fig:Drn}. The element-to-node differential
operators are derived in relation to the node-to-element operators
by requiring the discrete form of equation \ref{eq:500.0} to hold:
\begin{equation}
\overline{dV}{}^{T}*\left[\overline{\mathbf{P}}\cdot\left(\widehat{\overline{\nabla}}\,\,\widehat{U}\right)\right]+\widehat{dV}{}^{T}*\left[\widehat{U}\,\,\left(\overline{\widehat{\nabla}}\cdot\overline{\mathbf{P}}\right)\right]=0\label{eq:515.03}
\end{equation}
or 
\begin{equation}
\overline{dV}{}^{T}*\left[\overline{U}\,\,\left(\widehat{\overline{\nabla}}\cdot\widehat{\mathbf{P}}\right)\right]+\widehat{dV}{}^{T}*\left[\widehat{\mathbf{P}}\cdot\left(\overline{\widehat{\nabla}}\,\,\overline{U}\right)\right]=0\label{eq:515.031}
\end{equation}
In this derivation, it is assumed that the boundary term $\int u\mathbf{p}\cdot d\boldsymbol{\Gamma}=0$,
$i.e.,$ there is no flux of the continuous vector field $u\mathbf{p}$
through the boundary. For arbitrary discrete nodal representations
$\overline{\mathbf{P}}$ (equation \ref{eq:515.03}) and $\overline{U}$
(equation \ref{eq:515.031}), a consequence of this assumption is
that the element-to-node gradient operator produces valid results
at the boundary nodes only if the original continuous function $u$
is zero at the boundary $(i.e.,\,u|_{\Gamma}=0)$. Similarly, the
element-to-node divergence operator produces valid results at the
boundary nodes only if the original continuous function $\mathbf{p}$
has no component perpendicular to the boundary $(i.e.,\,\mathbf{p_{\perp}}|_{\Gamma}=0)$.
In the following, we will refer to these conditions as the natural
boundary conditions. For the terms involving radial derivatives, equation
\ref{eq:515.031} implies that

\begin{align*}
 & \overline{dV}^{T}*\left(\overline{U}\,\,\left(\left(\widehat{\overline{Dr}}*\left(\widehat{r}\,\,\widehat{P_{r}}\right)\right)\,/\,\overline{r}\right)\right) &  & =-\widehat{dV}^{T}*\left(\widehat{P_{r}}\,\,\left(\overline{\widehat{Dr}}*\overline{U}\right)\right)\\
 & \Rightarrow\frac{2\pi}{3}\,\overline{s}^{T}*\left(\overline{U}\,\,\left(\widehat{\overline{Dr}}*\left(\widehat{r}\,\,\widehat{P_{r}}\right)\right)\right) &  & =-2\pi\,\widehat{s}^{T}*\left(\widehat{r}\,\,\widehat{P_{r}}\,\,\left(\overline{\widehat{Dr}}*\overline{U}\right)\right)\\
 & \Rightarrow\frac{1}{3}\overline{U}^{T}*\left(\overline{\overline{S}}*\left(\widehat{\overline{Dr}}*\left(\widehat{r}\,\,\widehat{P_{r}}\right)\right)\right) &  & =-\left(\widehat{r}\,\,\widehat{P_{r}}\right)^{T}*\left(\widehat{\widehat{S}}*\overline{\widehat{Dr}}*\overline{U}\right)\\
 & \Rightarrow\frac{1}{3}\overline{U}^{T}*\left(\overline{\overline{S}}*\left(\widehat{\overline{Dr}}*\left(\widehat{r}\,\,\widehat{P_{r}}\right)\right)\right) &  & =-\overline{U}^{T}*\left(\overline{\widehat{Dr}}^{T}*\widehat{\widehat{S}}*\left(\widehat{r}\,\,\widehat{P_{r}}\right)\right) & \mbox{(transpose scalar on RHS)}
\end{align*}
Hence, the element-to-node derivative operators are 
\begin{align}
\widehat{\overline{Dr}} & =-3\overline{\overline{S}}^{-1}*\left(\overline{\widehat{Dr}}^{T}*\widehat{\widehat{S}}\right)\nonumber \\
\widehat{\overline{Dz}} & =-3\overline{\overline{S}}^{-1}*\left(\overline{\widehat{Dz}}^{T}*\widehat{\widehat{S}}\right)\label{eq:515.0}
\end{align}
These definitions can alternatively be obtained using equation \ref{eq:515.03}.
Noting the particular case when $\overline{U}=\overline{1}$, then
equation \ref{eq:515.031} leads to the identity 
\begin{align}
\overline{dV}{}^{T}*\left[\widehat{\overline{\nabla}}\cdot\widehat{\mathbf{P}}\right] & =0\label{eq:515.04}
\end{align}

\subsection{Second order node-to-node differential operators\label{subsec:Lapl_delstar}}

\begin{figure}[H]
\centering{}\includegraphics[width=6cm,height=2.5cm]{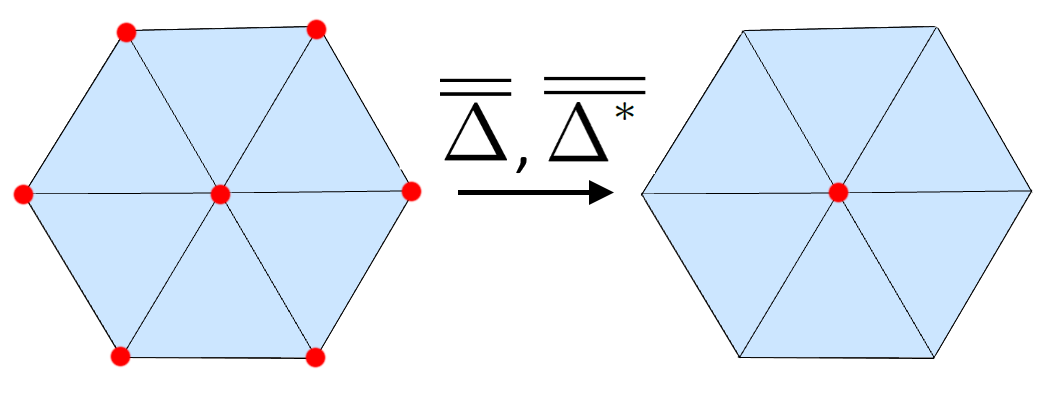}\caption{\label{fig:LaplDelstar} Node-to-node second order differential operator
mechanism}
\end{figure}
Discrete forms of node-to-node second order differential operators,
with dimensions $[N_{n}\times N_{n}]$, can be constructed using combinations
of the first order node-to-element and element-to-node operators.
For example, the discrete form of the Laplacian operator in cylindrical
coordinates, with azimuthal symmetry, is 
\begin{equation}
\overline{\overline{\Delta}}\,\,\overline{U}=\widehat{\overline{\nabla}}\cdot\overline{\widehat{\nabla}}\,\,\overline{U}=\left(\left(\widehat{\overline{Dr}}*\left(\widehat{r}\,\,\left(\overline{\widehat{Dr}}*\overline{U}\right)\right)\right)+\left(\widehat{\overline{Dz}}*\left(\widehat{r}\,\,\left(\overline{\widehat{Dz}}*\overline{U}\right)\right)\right)\right)\,/\,\overline{r}\label{eq:514.7}
\end{equation}
Note that $\overline{\overline{\Delta}}\,\,\overline{U}$ may be expressed
as $\widehat{\overline{\nabla}}\cdot\widehat{\mathbf{P}}$, where
$\widehat{\mathbf{P}}=\overline{\widehat{\nabla}}\,\,\overline{U}$.
Therefore, equation \ref{eq:515.04} implies that
\begin{align}
\overline{dV}{}^{T}*\left[\overline{\overline{\Delta}}\,\,\overline{U}\right] & =0\label{eq:515.71}
\end{align}
The condition $\mathbf{p_{\perp}}|_{\Gamma}=0$, required for correct
evaluation of $\widehat{\overline{\nabla}}\cdot\widehat{\mathbf{P}}$
at the boundary nodes, is equivalent to having the normal component
of the gradient of the continuous function $u$ equal to zero at the
boundary. Thus, the operation $\overline{\overline{\Delta}}*\overline{U}$
produces correct results at the boundary nodes only if the natural
boundary condition is satisfied, which in this case is $\left(\nabla_{\perp}u\right)|_{\Gamma}=0.$ 

The discrete form of the elliptic operator used in the Grad-Shafranov
equation, in cylindrical coordinates with azimuthal symmetry is

\begin{equation}
\overline{\overline{\Delta^{^{*}}}}\,\,\overline{U}=\overline{r}^{2}\,\,\left(\widehat{\overline{\nabla}}\cdot\left(\left(\overline{\widehat{\nabla}}\,\,\overline{U}\right)\,/\,\widehat{r}^{2}\right)\right)=\overline{r}\,\,\left(\widehat{\overline{Dr}}*\left(\left(\overline{\widehat{Dr}}*\overline{U}\right)\,/\,\widehat{r}\right)+\widehat{\overline{Dz}}*\left(\left(\overline{\widehat{Dz}}*\overline{U}\right)\,/\,\widehat{r}\right)\right)\label{eq:515.72}
\end{equation}
The $\overline{\overline{\Delta}}$ and $\overline{\overline{\Delta^{^{*}}}}$
operator mechanisms are indicated in figure \ref{fig:LaplDelstar}.
Note that equation \ref{eq:515.04} implies that 
\begin{align}
\overline{dV}{}^{T}*\left[\left(\overline{\overline{\Delta^{^{*}}}}\,\,\overline{U}\right)\,/\,\overline{r}^{2}\right] & =0\label{eq:515.71-1-1}
\end{align}
Again, the operation $\overline{\overline{\Delta^{^{*}}}}\,\,\overline{U}$
produces correct results at the boundary nodes only if the natural
boundary condition is satisfied, which in this case is $\left(\left(\nabla_{\perp}u\right)/r^{2}\right)|_{\Gamma}=0.$
As shown in \cite{thesis}, the standard Galerkin method may be used
to derive these second order node-to-node differential operators directly.

\subsection{First order node-to-node differential operators \label{subsec:Dr_Dz}}

\begin{figure}[H]
\centering{}\includegraphics[width=6cm,height=2.5cm]{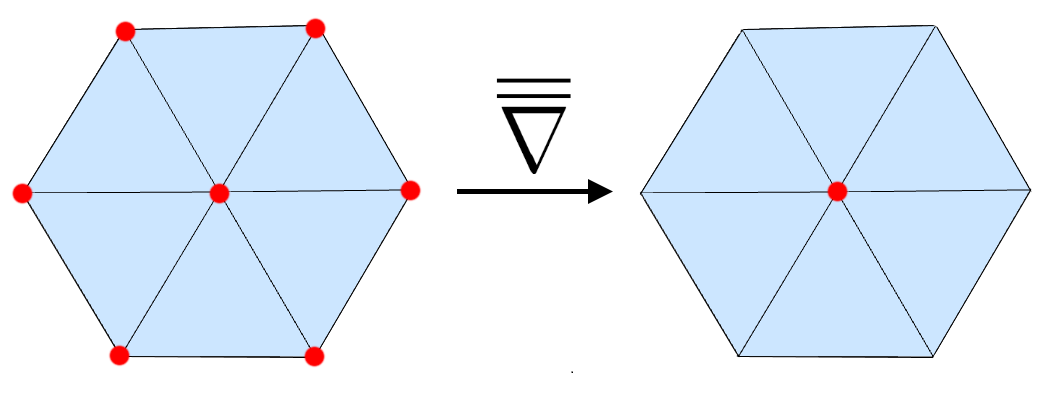}\caption{\label{fig:Dr} Node-to-node differential operator mechanism}
\end{figure}
The node-to-node derivative operator matrices $\overline{\overline{Dr}}$
and $\overline{\overline{Dz}}$, with dimensions $[N_{n}\times N_{n}]$,
are used to evaluate, at the nodes, the radial and axial derivatives
of any function $U(\mathbf{r})$ that is defined at the nodes, and
are defined as 
\begin{align*}
\overline{U}_{r}' & =\overline{\overline{Dr}}*\overline{U}\\
\overline{U}_{z}' & =\overline{\overline{Dz}}*\overline{U}
\end{align*}
The node-to-node gradient and divergence operations, for cylindrical
coordinates with azimuthal symmetry are 
\begin{align}
\overline{\overline{\nabla}}\,\,\overline{U} & =\left(\overline{\overline{Dr}}*\overline{U}\right)\hat{\mathbf{r}}+\left(\overline{\overline{Dz}}*\overline{U}\right)\hat{\mathbf{z}}\label{eq:510.21-1}
\end{align}
and 
\begin{align}
\overline{\overline{\nabla}}\cdot\overline{\mathbf{P}} & =\left(\left(\overline{\overline{Dr}}*\left(\overline{r}\,\,\overline{P}_{r}\right)\right)+\left(\overline{\overline{Dz}}*\left(\overline{r}\,\,\overline{P}_{z}\right)\right)\right)\,/\,\overline{r}\label{eq:510.22-1}
\end{align}
A schematic of the node-to-node gradient operation mechanism is indicated
in figure \ref{fig:Dr}. Here, we will look at the derivation of $\overline{\overline{Dr}}$,
the derivation of $\overline{\overline{Dz}}$ is analogous. The node-to-node
operators can be expressed in terms of the node-to-element operators
by considering the expansion of $U(\mathbf{r})$ in terms of the piecewise-linear
functions $U^{e}(\mathbf{r})$: 
\begin{align*}
U(\mathbf{r}) & =\overset{}{\underset{e}{\Sigma}}\,U^{e}(\textbf{\ensuremath{\mathbf{r}}})
\end{align*}
This implies that
\[
U_{r}'(\mathbf{r})=\overset{}{\underset{e}{\Sigma}}\,U_{r}^{e\,'}(\mathbf{r})
\]
Applying the Galerkin method, we obtain 
\begin{equation}
\int\phi{}_{j}(\mathbf{r})\,U_{r}'(\mathbf{r})\,dr\,dz=\int\phi{}_{j}(\textbf{\ensuremath{\mathbf{r}}})\,\overset{}{\underset{e}{\Sigma}}U_{r}^{e\,'}(\mathbf{r})\,dr\,dz\label{eq:507.01}
\end{equation}
where integral is taken over the support area of the basis function
$\phi_{j}(\mathbf{r})$. Note that $U_{r}^{e\,'}=B^{e}$ (equation
\ref{eq:500.02}). Expanding the right side of equation \ref{eq:507.01}
in terms of pyramid side functions according to equation \ref{eq:502.3},
and using equation \ref{eq:502.4} to expand the left side, this implies
that 
\begin{align*}
\frac{U_{rj}'s_{j}}{3} & =\overset{e_{j}}{\underset{}{\Sigma}}\int\psi_{j}^{e}(\textbf{\ensuremath{\mathbf{r}}})\,B^{e}\,dr\,dz\\
 & =\overset{e_{j}}{\underset{}{\Sigma}}\,B^{e}\int\psi_{j}^{e}(\textbf{\ensuremath{\mathbf{r}}})\,dr\,dz\\
 & =\overset{e_{j}}{\underset{}{\Sigma}}\,B^{e}\frac{s^{e}}{3}
\end{align*}
Therefore, in our method the derivative at a node is defined as the
area-weighted average of the elemental derivatives in the support
of that node. With use of connectivity operator $\overline{\widehat{M}}$,
this definition can be written as:

\begin{align}
\overline{U}_{r}' & =\left(\overline{\overline{S}}^{-1}*\left(\overline{\widehat{M}}^{T}*\widehat{\widehat{S}}*\overline{\widehat{Dr}}\right)\right)*\overline{U}\nonumber \\
\Rightarrow\overline{\overline{Dr}} & =\overline{\overline{S}}^{-1}*\left(\overline{\widehat{M}}^{T}*\widehat{\widehat{S}}*\overline{\widehat{Dr}}\right)\label{eq:510.1-1-1}
\end{align}
The derivation of $\overline{\overline{Dz}}$ follows from that of
$\overline{\overline{Dr}}$, with $\overline{\overline{Dz}}=\left(\overline{\overline{S}}^{-1}*\left(\overline{\widehat{M}}^{T}*\widehat{\widehat{S}}*\overline{\widehat{Dz}}\right)\right)$.
\\
\\
\\
The element-to-node differential operators are defined from the basic
node-to-element counterpart operators by enforcing the disrete form
of the combined expression for the divergence therorem and the differential
product rule (equation \ref{eq:500.0}). The node to node differential
operators obey the same laws: \\
\\

\begin{align}
\overline{dV}{}^{T}*\left[\left(\overline{U}\,\,\left(\overline{\overline{\nabla}}\cdot\overline{\mathbf{P}}\right)+\left(\overline{\overline{\nabla}}\,\,\overline{U}\right)\cdot\overline{\mathbf{P}}\right)\right] & =0 & \left(\mbox{if }\ensuremath{\ensuremath{\overline{U}\,|_{\Gamma}}=0}\mbox{ or }\ensuremath{\overline{\mathbf{P}}_{\perp}|_{\Gamma}=0}\right)\label{eq:511.04-1}
\end{align}
For the particular case of $\overline{U}=\overline{1}$ , this implies
that 
\begin{align}
\overline{dV}{}^{T}*\left[\overline{\overline{\nabla}}\cdot\overline{\mathbf{P}}\right] & =\overline{0} & \left(\mbox{if }\ensuremath{\overline{\mathbf{P}}_{\perp}|_{\Gamma}=0}\right)\label{eq:511.041-1}
\end{align}
More details and the proof can be found in appendix \ref{sec:ID_proff_n2n}.

\section{Discrete form of MHD equations and conservation laws\label{sec:Discrete-form-of}}

In this section we present the discretised equations for a two-temperature
MHD model, which are constructed so as to preserve the global conservation
laws inherent to the original continuous system of equations.

\subsection{Axisymmetric two-temperature MHD model\label{subsec:Axisymmetric-two-temperature-MHD}}

A two-temperature plasma MHD model, in which the electrons and ions
are allowed to have different temperatures, is a compromise between
the single-fluid and two-fluid MHD models. It has the relative computational
simplicity of the single-fluid model, but from the energy point of
view it treats electrons and ions as two separate fluids, and allows
implementation of different thermal diffusion coefficients for the
electrons and ions. In the following, SI units will be used, but temperatures
will be expressed in Joules unless otherwise indicated. The continuous
forms of the two-temperature MHD equations are:\\
\\
 
\begin{equation}
\dot{n}=-\nabla\cdot(n\mathbf{v})\label{eq: 478}
\end{equation}

\begin{equation}
\dot{\mathbf{v}}=-\mathbf{v}\cdot\nabla\mathbf{v}+\frac{1}{\rho}\left(-\nabla p-\nabla\cdot\underline{\boldsymbol{\pi}}+\mathbf{J\times}\mathbf{B}\right)\label{eq:479}
\end{equation}

\begin{equation}
\dot{p_{i}}=-\mathbf{v}\cdot\nabla p_{i}-\gamma p_{i}\,\nabla\cdot\mathbf{v}+(\gamma-1)\left(-\nabla\cdot\mathbf{q}_{i}+Q_{ie}-\underline{\boldsymbol{\pi}}:\nabla\mathbf{v}\right)\label{eq:479.1}
\end{equation}

\begin{equation}
\dot{p_{e}}=-\mathbf{v}\cdot\nabla p_{e}-\gamma p_{e}\,\nabla\cdot\mathbf{v}+(\gamma-1)\left(-\nabla\cdot\mathbf{q}_{e}-Q_{ie}+\eta'\mathbf{J}^{2}\right)\label{eq:479.2}
\end{equation}

\begin{equation}
\dot{\mathbf{B}}=\nabla\times\left(\mathbf{v}\times\mathbf{B}\right)-\nabla\times\left(\eta\nabla\times\mathbf{B}\right),\,\,\nabla\cdot\mathbf{B}=0\label{eq:479.3}
\end{equation}
Here, a dot above a symbol implies a partial time derivative, $\mu_{0}$
is the magnetic constant, $\eta'[\Omega-\mbox{m}]$ is the plasma
resistivity, $\eta[\mbox{m}^{2}/\mbox{s}]=\eta'/\mu_{0}$ is magnetic
diffusivity, and the adiabatic index $\gamma=\frac{5}{3}$ is the
same for ions and electrons. Note that $n$ is the ion number density,
so that the electron number density is $n_{e}=Zn$, where $Z$ is
the (volume averaged) ion charge, and $\rho=m_{i}n$ is the plasma
fluid mass density. We assume ideal gas laws for ions and electrons:
\begin{eqnarray}
p_{i} & = & nT_{i}\nonumber \\
p_{e} & = & ZnT_{e}\label{eq:479.31}
\end{eqnarray}
Total plasma fluid pressure is $p=p_{i}+p_{e}$. The remaining notations
are standard. The boundary conditions and closure for this model (namely,
definitions of thermal fluxes $\mathbf{q}_{i}$ and $\mathbf{q}_{e}$,
viscous stress tensor $\underline{\boldsymbol{\pi}}$ and ion-electron
heat exchange rate $Q_{ie}$) will be discussed in section \ref{subsec:Discretised-MHD-model}. 

With axisymmetry, the magnetic field can be represented in divergence-free
form with the poloidal flux (per radian) function $\psi(r,z)$ and
toroidal function $f(r,z)=rB_{\phi}$ as 
\begin{equation}
\mathbf{B=\nabla\psi\times\nabla\phi+}f\nabla\phi\label{eq:480}
\end{equation}
In this case, the current density is
\begin{equation}
\mathbf{J}=\frac{1}{\mu_{0}}\nabla\times\mathbf{B}=\frac{1}{\mu_{0}}\left(-\Delta^{*}\psi\,\nabla\phi+\nabla f\times\nabla\phi\right)\label{eq:480.1}
\end{equation}
and the Lorentz force is 
\begin{equation}
\mathbf{J}\times\mathbf{B}=-\frac{1}{\mu_{0}r^{2}}\bigg(\Delta^{*}\psi\,\nabla\psi+f\,\nabla f\bigg)+\frac{\mathbf{B}\cdot\nabla f}{\mu_{0}r}\widehat{\boldsymbol{\phi}}\label{eq:480.2}
\end{equation}
Using an axisymmetric representation of magnetic field, the magnetic
induction equation \ref{eq:479.3} can be re-written as 
\begin{eqnarray}
\dot{\psi} & = & -\mathbf{v}\cdot\nabla\psi+\eta\Delta^{*}\psi\label{eq:480.3}\\
\dot{f} & = & r^{2}\,\nabla\cdot\left(-\left(\frac{f}{r^{2}}\mathbf{v}\right)+\omega\mathbf{B}+\frac{\eta}{r^{2}}\nabla f\right)\label{eq:480.4}
\end{eqnarray}
 where $\omega=v_{\phi}/r$ is the plasma fluid angular speed.

\subsubsection{Conservation properties (continuous equations)\label{subsec:Conservation-properties-(continu}}

The system described by equations \ref{eq: 478}-\ref{eq:479.2},
\ref{eq:480.3}, \ref{eq:480.4} has a number of exact conservation
laws expressed by corresponding continuity equations. The obvious
one is conservation of particles, expressed by equation \ref{eq: 478},
with the boundary condition $\mathbf{v}_{\perp}|_{\Gamma}=0$, corresponding
to impermeable walls, where the subscript $\perp$ implies the component
perpendicular to the boundary. Another conserved quantity is toroidal
flux, defined as 
\[
\Phi=\int B_{\phi}\,dr\,dz=\int\frac{f}{r}\,dr\,dz=\frac{1}{2\pi}\int\frac{f}{r^{2}}2\pi r\,dr\,dz=\frac{1}{2\pi}\int\frac{f}{r^{2}}\,dV
\]
Equation \ref{eq:480.4} implies that the rate of change of system
toroidal flux is 
\begin{equation}
\dot{\Phi}=\frac{1}{2\pi}\int\nabla\cdot\left(-\frac{f}{r^{2}}\mathbf{v}+\omega\mathbf{B}+\frac{\eta}{r^{2}}\nabla f\right)\,dV=\frac{1}{2\pi}\int\left(-\frac{f}{r^{2}}\mathbf{v}+\omega\mathbf{B}+\frac{\eta}{r^{2}}\nabla f\right)\cdot d\boldsymbol{\Gamma}\label{eq:480.5}
\end{equation}
With appropriate boundary conditions, for example $\mathbf{v}|_{\Gamma}=\mathbf{0}$
and $\left(\nabla_{\perp}f\right)|_{\Gamma}=0$, or $\mathbf{v}_{\perp}|_{\Gamma}=\mathbf{0}$
and $\left(B_{\theta\perp}\,v_{\phi}+\eta\left(\frac{\nabla_{\perp}\,f}{r}\right)\right)|_{\Gamma}=0$,
toroidal flux is conserved. Note that the electric field at the system
boundary is 

\begin{align*}
\mathbf{E}|_{\Gamma} & =\left(-\mathbf{v\times}\mathbf{B}+\eta\nabla\times\mathbf{B}\right)|_{\Gamma}\\
 & =\left(-\mathbf{v_{\phi}\times}\mathbf{B}_{\theta}-\mathbf{v_{\theta}\times}\mathbf{B}_{\phi}-\mathbf{v_{\theta}\times}\mathbf{B}_{\theta}+\eta\left(-\Delta^{*}\psi\,\nabla\phi+\nabla f\times\nabla\phi\right)\right)|_{\Gamma} & \mbox{\mbox{(use eqn. \ref{eq:480.1})}}\\
 & =\left(-\mathbf{v_{\phi}\times}\left(\mathbf{B}_{\theta\parallel}+\mathbf{B}_{\theta\perp}\right)-B_{\phi}\,\left(\mathbf{v_{\theta\perp}}+\mathbf{v_{\theta\parallel}}\right)\times\widehat{\boldsymbol{\phi}}+(\mathbf{v}\cdot\nabla\psi)\nabla\phi+\eta\left(-\Delta^{*}\psi\,\nabla\phi+\frac{\nabla f}{r}\mathbf{\times}\widehat{\boldsymbol{\phi}}\right)\right)|_{\Gamma}
\end{align*}
Here, the subscript $\parallel$ implies the component parallel to
the boundary. Hence, the poloidal component of the electric field
parallel to the boundary is 
\begin{equation}
E_{\theta\parallel}|_{\Gamma}=\left(B_{\theta\perp}\,v_{\phi}-B_{\phi}\,v_{\theta\perp}+\eta\frac{\nabla_{\perp}f}{r}\right)|_{\Gamma}\label{eq:480.51}
\end{equation}
Therefore, either of the sets of boundary conditions listed above
for toroidal flux conservation correspond to having the poloidal component
of the  electric field perpendicular to the boundary, the condition,
in the case of azimuthal symmetry, for perfectly electrically conducting
walls.

Conservation of angular momentum is established after noting that
in the axisymmetric case 
\[
(\mathbf{v}\cdot\nabla\mathbf{v})_{\phi}=v_{r}\frac{\partial v_{\phi}}{\partial r}+\frac{v_{r}v_{\phi}}{r}+v_{z}\frac{\partial v_{\phi}}{\partial z}=\frac{v_{r}}{r}\frac{\partial(rv_{\phi})}{\partial r}+\frac{v_{z}}{r}\frac{\partial(rv_{\phi})}{\partial z}=\frac{1}{r}\mathbf{v}\cdot\nabla(rv_{\phi})
\]
and that the $\phi$ coordinate of the divergence of the viscous stress
tensor, which is transpose symmetric, is 
\[
(\nabla\cdot\underline{\boldsymbol{\pi}})_{\phi}=\frac{1}{r}\frac{\partial}{\partial r}(r\pi_{r\phi})+\frac{\partial}{\partial z}(\pi_{z\phi})+\frac{\pi_{r\phi}}{r}=\frac{1}{r^{2}}\frac{\partial}{\partial r}(r^{2}\pi_{r\phi})+\frac{1}{r^{2}}\frac{\partial}{\partial z}(r^{2}\pi_{z\phi})=\frac{1}{r}\nabla\cdot(r\pi_{r\phi}\hat{\mathbf{r}}+r\pi_{z\phi}\hat{\mathbf{z}})
\]
The continuity equation for angular momentum density is then 
\begin{eqnarray}
 &  & \frac{\partial}{\partial t}(\rho v_{\phi}r)=m_{i}\dot{n}v_{\phi}r+\rho\dot{v}_{\phi}r\nonumber \\
 & = & -\biggl[m_{i}v_{\phi}r\nabla\cdot(n\mathbf{v})+\rho\mathbf{v}\cdot\nabla(rv_{\phi})\biggr]-\nabla\cdot(r\pi_{r\phi}\hat{\mathbf{r}}+r\pi_{z\phi}\hat{\mathbf{z}})+\frac{1}{\mu_{0}}\mathbf{B}\cdot\nabla f\label{eq:481}\\
 & = & -\nabla\cdot\bigg(\rho v_{\phi}r\mathbf{v}+r\pi_{r\phi}\hat{\mathbf{r}}+r\pi_{z\phi}\hat{\mathbf{z}}-\frac{f}{\mu_{0}}\mathbf{B}\bigg)\nonumber 
\end{eqnarray}
Hence, the rate of change of total system angular momentum is 
\[
\dot{P}_{\phi}=\frac{\partial}{\partial t}\left(\int\left(\rho rv_{\phi}\right)dV\right)=-\int\left(\rho v_{\phi}r\mathbf{v}+r\pi_{r\phi}\hat{\mathbf{r}}+r\pi_{z\phi}\hat{\mathbf{z}}-\frac{f}{\mu_{0}}\mathbf{B}\right)\cdot d\boldsymbol{\Gamma}
\]
With appropriate boundary conditions, total system angular momentum
is conserved. One set of appropriate boundary conditions is as follows.
The first term here vanishes with the boundary conditions for impermeable
walls, $\mathbf{v}_{\perp}|_{\Gamma}=\mathbf{0}$. Referring to equation
\ref{eq:238-1}, it can be seen that with the boundary condition $\left(\nabla_{\perp}\omega\right)|_{\Gamma}=0$,
the second and third terms vanish too. Note that angular momentum
conservation is not physical when there is friction between the fluids
and the boundary wall. The boundary condition $\left(\nabla_{\perp}\omega\right)|_{\Gamma}=0$
is physical only in unusual cases; for example, there may be no viscosity
at the wall due to rotation of the wall. The last term is zero with
the boundary condition $\left(\nabla_{\parallel}\psi\right)|_{\Gamma}=0$
($i.e.,$ $\psi|_{\Gamma}=C$ (constant)), corresponding to having
the magnetic field parallel to the boundary.

To prove energy conservation, we note that $-\mathbf{v}\cdot\nabla\mathbf{v}=-\nabla\left(v^{2}/2\right)+\mathbf{v}\times(\nabla\times\mathbf{v})$.
Expressing the total system energy as the sum of kinetic, thermal
and magnetic energy, the continuity equation for energy density is
\begin{eqnarray}
\dot{u}_{Total} & = & \dot{u}_{K}+\dot{u}_{Th}+\dot{u}_{M}\nonumber \\
 & = & \frac{\partial}{\partial t}\left(\frac{\rho\mathbf{v}^{2}}{2}+\frac{p}{\gamma-1}+\frac{1}{2\mu_{0}r^{2}}\left(\left(\nabla\psi\right){}^{2}+f^{2}\right)\right)\nonumber \\
 & = & \frac{m_{i}\dot{n}\mathbf{v}^{2}}{2}+\rho\mathbf{v}\cdot\dot{\mathbf{v}}+\frac{\dot{p}}{\gamma-1}+\frac{1}{\mu_{0}r^{2}}\left(\nabla\psi\cdot\left(\nabla\dot{\psi}\right)+f\dot{f}\right)\nonumber \\
 & = & -\left[\frac{m_{i}\mathbf{v}^{2}}{2}\nabla\cdot(n\mathbf{v})+\rho\mathbf{v}\cdot\nabla\frac{\mathbf{v}^{2}}{2}\right]\nonumber \\
 &  & -\left[\mathbf{v}\cdot\nabla p+\frac{\mathbf{v}\cdot\nabla p+\gamma p\nabla\cdot\mathbf{v}}{\gamma-1}\right]\nonumber \\
 &  & -\left[\frac{1}{\mu_{0}r^{2}}\Delta^{*}\psi(\mathbf{v}\cdot\nabla\psi)+\frac{1}{\mu_{0}r^{2}}\nabla\psi\cdot\nabla(\mathbf{v}\cdot\nabla\psi)\right]\nonumber \\
 &  & -\left[\frac{1}{\mu_{0}}\frac{f}{r^{2}}(\mathbf{v}\cdot\nabla f)+\frac{1}{\mu_{0}}f\,\nabla\cdot\bigg(\frac{f}{r^{2}}\mathbf{v}\bigg)\right]\nonumber \\
 &  & +\left[\frac{1}{\mu_{0}}\omega\mathbf{B}\cdot\nabla f+\frac{1}{\mu_{0}}f\nabla\cdot\bigg(\omega\mathbf{B}\bigg)\right]\nonumber \\
 &  & +\left[\frac{\eta(\Delta^{*}\psi)^{2}}{\mu_{0}r^{2}}+\frac{1}{\mu_{0}r^{2}}\nabla\psi\cdot\nabla(\eta\Delta^{*}\psi)\right]\nonumber \\
 &  & +\left[\frac{\eta(\nabla f)^{2}}{\mu_{0}r^{2}}+\frac{f}{\mu_{0}}\nabla\cdot\bigg(\frac{\eta}{r^{2}}\nabla f\bigg)\right]\label{eq:481.1}\\
 &  & -\biggl[\nabla\cdot(\mathbf{q}_{i}+\mathbf{q}_{e})\biggr]\nonumber \\
 &  & -\biggl[\mathbf{v}\cdot(\nabla\cdot\underline{\boldsymbol{\pi}})+\underline{\boldsymbol{\pi}}:\nabla\mathbf{v}\biggr]\nonumber \\
 & = & -\nabla\cdot\biggl(\frac{m_{i}n\mathbf{v}^{2}}{2}\mathbf{v}+\frac{\gamma p}{\gamma-1}\mathbf{v}+\frac{(\mathbf{v}\cdot\nabla\psi)}{\mu_{0}r^{2}}\nabla\psi\nonumber \\
 &  & +\frac{f^{2}}{\mu_{0}r^{2}}\mathbf{v}-\frac{1}{\mu_{0}}\omega f\,\mathbf{B}-\frac{\eta\Delta^{*}\psi}{\mu_{0}r^{2}}\nabla\psi-\frac{\eta f}{\mu_{0}r^{2}}\nabla f+\mathbf{q}_{i}+\mathbf{q}_{e}+\underline{\boldsymbol{\pi}}\cdot\mathbf{v}\biggr)\nonumber 
\end{eqnarray}
Note that the five magnetic terms involving $\mu_{0}$ in the final
full-divergence expression constitute the Poynting flux. Taking the
integral over the final expression for $\dot{u}_{total}$ over the
system volume, and applying Gauss's theorem, it can be seen how total
system energy is conserved with appropriate boundary conditions, for
example $\mathbf{v}|_{\Gamma}=\mathbf{0},\,\mathbf{q}_{i\perp}|_{\Gamma}=\mathbf{q}_{e\perp}|_{\Gamma}=\mathbf{0},\,(\nabla_{\perp}\psi)|_{\Gamma}=0\mbox{ and }(\nabla_{\perp}f)|_{\Gamma}=0$,
or $\mathbf{v}_{\perp}|_{\Gamma}=\mathbf{0},\,\mathbf{q}_{i\perp}|_{\Gamma}=\mathbf{q}_{e\perp}|_{\Gamma}=\mathbf{0},\,(\nabla_{\perp}\psi)|_{\Gamma}=0\mbox{ and }\left(B_{\theta\perp}\,v_{\phi}+\eta\left(\frac{\nabla_{\perp}\,f}{r}\right)\right)|_{\Gamma}=0$.
Referring to equation \ref{eq:480.51}, note that either of these
sets of boundary conditions also eliminate energy loss associated
with Poynting flux. Note that the terms $+\frac{(\mathbf{v}\cdot\nabla\psi)}{\mu_{0}r^{2}}\nabla\psi$
and $-\frac{\eta\Delta^{*}\psi}{\mu_{0}r^{2}}\nabla\psi$ in the final
full-divergence expression in equation \ref{eq:481.1} constitute
the part of the Poynting flux that arises due to the toroidal component
of the electric field at the boundary. The boundary condition $(\nabla_{\perp}\psi)|_{\Gamma}=0$
or $(\Delta^{*}\psi)|_{\Gamma}=0$ is required to eliminate this contribution.
Referring to equation \ref{eq:480.3}, it can be seen that the combination
of boundary conditions $\mathbf{v}_{\perp}|_{\Gamma}=\mathbf{0}$
and $\psi|_{\Gamma}=0$ automatically leads to the boundary condition
$(\Delta^{*}\psi)|_{\Gamma}=0$. Therefore, the boundary condition
$(\nabla_{\perp}\psi)|_{\Gamma}=0$ included above in the list of
requirements for maintenance of system energy conservation may be
replaced with the requirement $\psi|_{\Gamma}=0$. 

\subsection{Discretised MHD model\label{subsec:Discretised-MHD-model}}

In this section, the full set of discretised equations for a two-temperature
MHD model is presented. The model has been constructed so as to preserve
the global conservation laws inherent to the original continuous system
of equations. In order to preserve the conservation laws of the system
in discrete form, the pairs of terms which constitute full divergences,
denoted with square brackets in equations \ref{eq:481} and \ref{eq:481.1},
have been discretised in compatible way, $i.e.,$ the appropriate
corresponding discrete operators are used in these terms, and with
particular boundary conditions, the operator properties lead to the
cancellation of each pair. The discretisation may not be unique -
several different schemes can be obtained with the same conservation
properties. More details on the construction of the scheme are given
in \cite{thesis}. The resulting set of discretised equations is:
\begin{align}
\dot{\overline{n}} & =-\overline{\overline{\nabla}}\cdot(\overline{n}\,\,\overline{\mathbf{v}})+\widehat{\overline{\nabla}}\cdot\left(\widehat{\zeta}\,\,\left(\overline{\widehat{\nabla}}\,\,\overline{n}\right)\right)\nonumber \\
\dot{\overline{v}}_{r} & =\underbrace{-\overline{\overline{Dr}}*\left(\frac{\overline{v}^{2}}{2}\right)-\overline{v}_{z}\,\,(\overline{\overline{Dz}}*\overline{v}_{r}-\overline{\overline{Dr}}*\overline{v}_{z})+\overline{v}_{\phi}\,\,(\overline{\overline{Dr}}*(\overline{r}\,\,\overline{v}_{\phi}))\,/\,\overline{r})}_{\mathclap{{-(\mathbf{v\cdot\nabla}\mathbf{v})_{r}=\left(-\nabla(v^{2}/2)+\mathbf{v}\times(\nabla\times\mathbf{v})\right)_{r}}}}-\underbrace{\left(\overline{\overline{Dr}}*\overline{p}\right)\,/\,\overline{\rho}}_{{\frac{1}{\rho}(\nabla p)_{r}}}-\underbrace{\overline{\Pi}_{r}\,/\,\overline{\rho}}_{\frac{1}{\rho}(\nabla\cdot\underline{\boldsymbol{\pi}})_{r}}\nonumber \\
 & +\underbrace{\left(\underset{}{-(\overline{\overline{Dr}}*\overline{\psi})}\,\,\left(\overline{\overline{\Delta^{^{*}}}}\,\,\overline{\psi}\right)-\left(\overline{f}\,\,\left(\overline{\overline{Dr}}*\overline{f}\right)\right)\right)\,/\,\left(\mu_{0}\,\overline{r}{}^{2}\,\,\overline{\rho}\right)}_{\frac{1}{\rho}(\mathbf{J\times}\mathbf{B})_{r}}+\overline{f}_{\zeta r}\,/\,\overline{\rho}\nonumber \\
\dot{\overline{v}}_{\phi} & =\underbrace{-\overline{\mathbf{v}}\cdot\left(\overline{\overline{\nabla}}\,\,(\overline{r}\,\,\overline{v}_{\phi})\right)\,/\,\overline{r}}_{\mathclap{{-(\mathbf{v\cdot\nabla}\mathbf{v})_{\phi}=\left(-\nabla(v^{2}/2)+\mathbf{v}\times(\nabla\times\mathbf{v})\right)_{\phi}}}}\,\,\,\,\,\,\,\,\,\,-\underbrace{\overline{\Pi}_{\phi}\,/\,\overline{\rho}}_{\frac{1}{\rho}(\nabla\cdot\underline{\boldsymbol{\pi}})_{\phi}}+\underbrace{\left(\widehat{\overline{W}}*\left(\widehat{\mathbf{B}}_{\theta}\cdot\left(\overline{\widehat{\nabla}}\,\,\overline{f}\right)\right)\right)\,/\,\left(\mu_{0}\,\overline{r}\,\,\overline{\rho}\right)}_{\frac{1}{\rho}(\mathbf{J\times}\mathbf{B})_{\phi}}+\overline{f}_{\zeta\phi}\,/\,\overline{\rho}\nonumber \\
\dot{\overline{v}}_{z} & =\underbrace{-\overline{\overline{Dz}}*\left(\frac{\overline{v}^{2}}{2}\right)+\overline{v}_{r}\,\,(\overline{\overline{Dz}}*\overline{v}_{r}-\overline{\overline{Dr}}*\overline{v}_{z})+\overline{v}_{\phi}\,\,(\overline{\overline{Dz}}*(\overline{r}\,\,\overline{v}_{\phi}))\,/\,\overline{r})}_{\mathclap{{-(\mathbf{v\cdot\nabla}\mathbf{v})_{z}=\left(-\nabla(v^{2}/2)+\mathbf{v}\times(\nabla\times\mathbf{v})\right)_{z}}}}-\underbrace{\left(\overline{\overline{Dz}}*\overline{p}\right)\,/\,\overline{\rho}}_{{\frac{1}{\rho}(\nabla p)_{z}}}-\underbrace{\overline{\Pi}_{z}\,/\,\overline{\rho}}_{\frac{1}{\rho}(\nabla\cdot\underline{\boldsymbol{\pi}})_{z}}\nonumber \\
 & +\underbrace{\left(-\left(\overline{\overline{Dz}}*\overline{\psi}\right)\,\,\left(\overline{\overline{\Delta^{^{*}}}}\,\,\overline{\psi}\right)-\overline{f}\,\,\left(\overline{\overline{Dz}}*\overline{f}\right)\right)\,/\,\left(\mu_{0}\,\overline{r}{}^{2}\,\,\overline{\rho}\right)}_{\frac{1}{\rho}(\mathbf{J\times}\mathbf{B})_{z}}+\overline{f}_{\zeta z}\,/\,\overline{\rho}\nonumber \\
\dot{\overline{p}}_{i} & =-\overline{\mathbf{v}}\cdot\left(\overline{\overline{\nabla}}\,\,\overline{p}_{i}\right)-\gamma\,\overline{p}_{i}\,\,\left(\overline{\overline{\nabla}}\cdot\overline{\mathbf{v}}\right)+(\gamma-1)\,\left[-\widehat{\overline{\nabla}}\cdot\widehat{\mathbf{q}}_{i}+\overline{Q}_{ie}+\underbrace{\overline{Q}_{\pi}}_{-\underline{\boldsymbol{\pi}}:\nabla\mathbf{v}}+\overline{Q}_{\zeta}\right]\nonumber \\
\dot{\overline{p}}_{e} & =-\overline{\mathbf{v}}\cdot\left(\overline{\overline{\nabla}}\,\,\overline{p}_{e}\right)-\gamma\,\overline{p}_{e}\,\,\left(\overline{\overline{\nabla}}\cdot\overline{\mathbf{v}}\right)+(\gamma-1)\,\biggl[-\widehat{\overline{\nabla}}\cdot\widehat{\mathbf{q}}_{e}-\overline{Q}_{ie}\nonumber \\
 & +\underset{\eta'J_{\phi}^{2}}{\underbrace{\left(\overline{\eta}/\mu_{0}\right)\,\,\left(\left(\overline{\overline{\Delta^{^{*}}}}\,\,\overline{\psi}\right)\,/\,\overline{r}\right)^{2}}}+\underset{\eta'J_{\theta}^{2}}{\underbrace{\widehat{\overline{W}}*\left(\left(\widehat{\mathbf{\eta}}/\mu_{0}\right)\,\,\left(\left(\overline{\widehat{\nabla}}\,\,\overline{f}\right)\,/\,\widehat{r}\right)^{2}\right)}}\biggr]\nonumber \\
\dot{\overline{\psi}} & =-\overline{\mathbf{v}}\cdot\left(\overline{\overline{\nabla}}\,\,\overline{\psi}\right)+\overline{\eta}\,\,\left(\overline{\overline{\Delta^{^{*}}}}\,\,\overline{\psi}\right)\nonumber \\
\dot{\overline{f}} & =\underline{r}^{2}\,\,\left[-\overline{\overline{\nabla}}\cdot\left(\overline{f}\,\,\overline{\mathbf{v}}\,/\,\overline{r}^{2}\right)+\widehat{\overline{\nabla}}\cdot\left(\widehat{\mathbf{B}}_{\theta}\,\,\widehat{\mathbf{\omega}}\right)+\widehat{\overline{\nabla}}\cdot\left(\widehat{\mathbf{\eta}}\,\,\left(\overline{\widehat{\nabla}}\,\,\overline{f}\right)\,/\,\widehat{r}^{2}\right)\right]\label{eq:517.3}
\end{align}
Here, $\widehat{\mathbf{\eta}}=<\overline{\eta}>^{e}$, and $\widehat{\mathbf{\omega}}=<\overline{v}_{\phi}\,/\,\overline{r}>^{e}$
(equation \ref{eq:502.31}). $\widehat{\overline{W}}$ is the volume-averaging
operator (equation \ref{eq:515.1}), and $\widehat{\mathbf{B}}_{\theta}=\left(-\left(\overline{\widehat{Dz}}*\overline{\psi}\right)\hat{\mathbf{r}}+\left(\overline{\widehat{Dr}}*\overline{\psi}\right)\hat{\mathbf{z}}\right)\,/\,\widehat{r}$.
These discretised equations are written in a form which is a direct
analogue of their continuous representations with inclusion of several
extra terms. Namely, the continuity equation has an artificial density
diffusion term, which is required for density smoothing and avoiding
negative density regions. The density diffusion coefficient $\zeta\,[\mbox{m}^{2}/\mbox{s}]$
can be spatially held constant or can be increased where density gradients
are high or density values approach zero. The components of force
per volume vector $\overline{\mathbf{f}}_{\zeta}=(\overline{f}_{\zeta r},\,\overline{f}_{\zeta\phi},\,\overline{f}_{\zeta z})^{T}$
in the velocity equations and the heating term $\overline{Q}_{\zeta}$
in the ion pressure equation are included to cancel the effect of
artificial density diffusion on the total system momentum and energy.
The expressions for these terms will be defined in section \ref{subsec:Maintenance-of-momentum}.

To close the system of equations, we need to specify the forms of
the viscous terms $\overline{\boldsymbol{\varPi}}$ and $\overline{Q}_{\pi}$,
the species heat exchange term $\overline{Q}_{ie}$, the resistive
diffusion coefficient $\overline{\eta}$, and the heat flux density
terms $\widehat{\mathbf{q}}_{i}$ and $\widehat{\mathbf{q}}_{e}$.
For the viscous terms, we use, for simplicity, the unmagnetised version
of the viscous stress tensor in axisymmetric cylindrical coordinates,
which is given as \cite{Cloutman}: \\
\begin{equation}
\underline{\boldsymbol{\pi}}=-\mu\left(\begin{array}{ccc}
2\frac{\partial v_{r}}{\partial r}-\frac{2}{3}\nabla\cdot\mathbf{v} & r\frac{\partial}{\partial r}\left(\frac{v_{\phi}}{r}\right) & \frac{\partial v_{z}}{\partial r}+\frac{\partial v_{r}}{\partial z}\\
r\frac{\partial}{\partial r}\left(\frac{v_{\phi}}{r}\right) & 2\frac{v_{r}}{r}-\frac{2}{3}\nabla\cdot\mathbf{v} & \frac{\partial v_{\phi}}{\partial z}\\
\frac{\partial v_{z}}{\partial r}+\frac{\partial v_{r}}{\partial z} & \frac{\partial v_{\phi}}{\partial z} & 2\frac{\partial v_{z}}{\partial z}-\frac{2}{3}\nabla\cdot\mathbf{v}
\end{array}\right)\label{eq:238-1}
\end{equation}
where $\mu\,${[}kg$\,$m$^{-1}\mbox{s}{}^{-1}${]} is dynamic viscosity.
The components of $\overline{\boldsymbol{\varPi}}$ represent the
discrete forms of the components of $\nabla\cdot\boldsymbol{\underline{\pi}}$,
and are defined as 
\begin{align}
\overline{\Pi}_{r} & =\biggl[\underset{}{-2\widehat{\overline{Dr}}*\left(\widehat{\mathbf{\mu}}\,\,\widehat{r}\,\,\left(\overline{\widehat{Dr}}*\overline{v}_{r}\right)\right)}-\underset{}{\widehat{\overline{Dz}}*\left(\widehat{\mathbf{\mu}}\,\,\widehat{r}\,\,\left(\overline{\widehat{Dr}}*\overline{v}_{z}+\overline{\widehat{Dz}}*\overline{v}_{r}\right)\right)}\biggr]\,/\,\overline{r}\nonumber \\
 & \underset{}{+\frac{2}{3}\left(\widehat{\overline{Dr}}*\left(\widehat{\mathbf{\mu}}\,\,\left(\overline{\widehat{\nabla}}\cdot\overline{\mathbf{v}}\right)\right)\right)}+\underset{}{2\,\overline{\mu}\,\,\overline{v}_{r}\,/\,\overline{r}^{2}}\nonumber \\
\overline{\Pi}_{\phi} & =-\left(\widehat{\overline{\nabla}}\cdot\left(\widehat{\mathbf{\mu}}\,\,\widehat{r}^{2}\,\,\left(\overline{\widehat{\nabla}}\,\,\overline{\omega}\right)\right)\right)\,/\,\overline{r}\nonumber \\
\overline{\Pi}_{z} & =\biggl[\underset{}{-2\widehat{\overline{Dz}}*\left(\widehat{\mathbf{\mu}}\,\,\widehat{r}\,\,\left(\overline{\widehat{Dz}}*\overline{v}_{z}\right)\right)}-\underset{}{\widehat{\overline{Dr}}*\left(\widehat{\mathbf{\mu}}\,\,\widehat{r}\,\,\left(\overline{\widehat{Dr}}*\overline{v}_{z}+\overline{\widehat{Dz}}*\overline{v}_{r}\right)\right)}\biggr]\,/\,\overline{r}\nonumber \\
 & \underset{}{+\frac{2}{3}\left(\widehat{\overline{Dz}}*\left(\widehat{\mathbf{\mu}}\,\,\left(\overline{\widehat{\nabla}}\cdot\overline{\mathbf{v}}\right)\right)\right)}\label{eq:517.4}
\end{align}
Here, $\widehat{\mathbf{\mu}}=<\overline{\mu}>^{e}$ (equation \ref{eq:502.31}).
The contraction (inner product) of two second order tensors is defined
as $\underline{\boldsymbol{T}}:\underline{\boldsymbol{U}}=\underset{i}{\Sigma}\underset{j}{\Sigma}T_{ij}U_{ij}$,
and $\nabla\mathbf{v}$, in cylindrical coordinates with azimuthal
symmetry, is

\[
\nabla\mathbf{v}=\left(\begin{array}{ccc}
\frac{\partial v_{r}}{\partial r} & -\frac{v_{\phi}}{r} & \frac{\partial v_{r}}{\partial z}\\
\\
\frac{\partial v_{\phi}}{\partial r} & \frac{v_{r}}{r} & \frac{\partial v_{\phi}}{\partial z}\\
\\
\frac{\partial v_{z}}{\partial r} & 0 & \frac{\partial v_{z}}{\partial z}
\end{array}\right)
\]
In the expression for $\dot{\overline{p}}_{i}$, $\overline{Q}_{\pi}$
is the discrete form of $-\underline{\boldsymbol{\pi}}:\nabla\mathbf{v}$:
\begin{align}
\overline{Q}_{\pi} & =\widehat{\overline{W}}*\biggl[\widehat{\mathbf{\mu}}\,\,\biggl\{2\left(\overline{\widehat{Dr}}*\overline{v}_{r}\right)^{2}+2\left(\overline{\widehat{Dz}}*\overline{v}_{z}\right)^{2}+\left(\widehat{r}\,\,\left(\overline{\widehat{\nabla}}\,\,\overline{\omega}\right)\right)^{2}\nonumber \\
 & +\left(\overline{\widehat{Dr}}*\overline{v}_{z}+\overline{\widehat{Dz}}*\overline{v}_{r}\right)^{2}-\frac{2}{3}\left(\overline{\widehat{\nabla}}\cdot\overline{\mathbf{v}}\right)^{2}\biggr\}\biggr]+2\,\overline{\mu}\,\,\left(\overline{v}_{r}\,/\,\overline{r}\right)^{2}\label{eq:517.5}
\end{align}
The representations in equations \ref{eq:517.4} and \ref{eq:517.5}
are consistent with the discrete form of energy conservation, as will
be shown in section \ref{subsec:Energy-conservation}. $\widehat{\mathbf{\mu}}\,[\mbox{kg}\,\mbox{m}^{-1}\mbox{s}{}^{-1}]=\widehat{\rho}\,\,\widehat{\nu}$
specifies dynamic viscosity at element centroids, where $\widehat{\nu}$
{[}m$^{2}\mbox{/s}${]} specifies kinematic viscosity. A certain minimum
level of artificial viscosity, which depends on the time-step, simulation
type, and mesh resolution, is required for numerical stability. Simulations
presented in this work, which involve CT formation and consequent
extreme plasma acceleration and steep gradients in the velocity fields,
are run with dynamic viscosity set to constant $\mu=\rho_{0}\nu_{0}$,
for typical code input $\nu_{0}\sim700${[}m$^{2}\mbox{/s}${]}, where
$\rho_{0}\sim6\times10^{-6}${[}kg/m$^{3}${]} is a typical representative
mass density. 

The resistive diffusion coefficient is based on the Spitzer formula

\[
\overline{\eta}=\frac{m_{e}}{1.96\,e^{2}\,\mu_{0}}\,/\,\left(Z\,\overline{n}\,\,\overline{\tau}_{ei}\right)=\left(418\,Z\left(\overline{T}_{e}\,[\mbox{eV}]\right)^{-\frac{3}{2}}\right)\,[\mbox{m}^{2}\mbox{/s}]
\]
where the electron-ion collision time is 
\[
\overline{\tau}_{ei}=\left(\frac{6\sqrt{2}\pi^{1.5}\epsilon_{0}^{2}\sqrt{m_{e}}}{\varLambda\,e^{4}Z^{2}}\right)\,\left(\overline{T}_{e}\,[\mbox{J}]\right){}^{\frac{3}{2}}\,/\,\overline{n}=\left(3.45\times10^{10}\,\frac{\left(\overline{T}_{e}\,[\mbox{eV}]\right)^{\frac{3}{2}}\,/\,\overline{n}\,[\mbox{m}^{-3}]}{Z^{2}}\right)\,[\mbox{s}]
\]
Here, $\Lambda=10$ is the Coulomb logarithm. In general, to limit
the timestep to acceptably high values, simulations are run with an
upper limit on $\eta$ of around 5000$[\mbox{m}^{2}\mbox{/s}]$.

The heat exchange term $\overline{Q}_{ie}$ gives the rate at which
energy is imparted from the electrons to the ions due to collisions
between ion and electron fluids: 
\begin{align*}
\overline{Q}_{ie} & =3(m_{e}/m_{i})\,Z\,\overline{n}\,\,\left(\overline{T}_{e}[\mbox{J}]-\overline{T}_{i}[\mbox{J}]\right)\,/\,\overline{\tau}_{ei}\\
 & =\left(7.6\times10^{-33}\,Z^{3}\,\left(\overline{T}_{e}[\mbox{eV}]-\overline{T}_{i}[\mbox{eV}]\right)\,\,\left(\overline{T}_{e}\,[\mbox{eV}]\right)^{-\frac{3}{2}}\,\,\left(\overline{n}\,[\mbox{m}^{-3}]\right)^{2}/\mu_{i}\right)\,[\mbox{W/\ensuremath{\mbox{m}^{3}}}]
\end{align*}
Here, $\mu_{i}$ is the ion mass in units of proton mass.

To include the effect of significantly enhanced thermal diffusion
parallel to the magnetic field, we include anisotropy in the model
for thermal diffusion. The species heat flux density can be expressed
as
\begin{align*}
\mathbf{q}_{\alpha} & =-\left(\kappa_{\parallel\alpha}\nabla_{\parallel}T_{\alpha}+\kappa{}_{\perp\alpha}\nabla_{\perp}T_{\alpha}\right)\\
 & =-\left(\left(\kappa_{\parallel\alpha}-\kappa_{\perp\alpha}\right)\nabla_{\parallel}T_{\alpha}+\kappa_{\perp\alpha}\nabla T_{\alpha}\right)
\end{align*}
where $\kappa_{\parallel\alpha}$ and $\kappa_{\perp\alpha}$ {[}(m-s)$^{-1}${]}
are the thermal conductivities for species $\alpha$, pertaining to
thermal diffusion parallel and perpendicular to the magnetic field.
With azimuthal symmetry, the toroidal component of $\nabla_{\parallel}T_{\alpha}$
can be dropped as it will not make a finite contribution to $\nabla\cdot\mathbf{q}_{\alpha}$
, so that the discrete forms of $\mathbf{q}_{\alpha}$ may be expressed
as 
\begin{align}
\widehat{\mathbf{q}}_{\alpha} & =-\left\{ \left(\widehat{\kappa}_{\parallel\alpha}-\widehat{\kappa}_{\perp\alpha}\right)\,\,\left(\widehat{\mathbf{B}}_{\theta}\,\,\left(\widehat{\mathbf{B}}_{\theta}\cdot\left(\overline{\widehat{\nabla}}\,\,\overline{T}_{\alpha}\right)\right)\,/\,\widehat{B}^{2}\right)+\widehat{\kappa}_{\perp\alpha}\,\,\left(\overline{\widehat{\nabla}}\,\,\overline{T}_{\alpha}\right)\right\} \label{eq:517.7}
\end{align}
Physically, the thermal conductivities vary with local conditions.
The simulations presented in this work were run with constant conductivities,
of the order $\kappa_{\parallel\alpha}=n_{0}\chi_{\parallel\alpha}$
and $\kappa_{\perp\alpha}=n_{0}\chi_{\perp\alpha}$, where $n_{0}$
{[}m$^{-3}${]} is a typical representative number density, and the
thermal diffusion coefficients $\chi_{\parallel\alpha},\,\chi_{\perp\alpha}$
$[\mbox{m}^{2}\mbox{/s}]$ are held constant. 

The code uses explicit time-stepping, with options of forward Euler,
Runge-Kutta 2 or Runge-Kutta 4 schemes. Although extremely low or
high diffusion coefficient values can be used while maintaining numerical
stability if the timestep is reduced, lower and upper limits need
to be imposed on diffusion coefficients in order to achieve numerical
stability in combination with acceptably short simulation runtimes.
This is consistent with the standard stability criteria \cite{JARDIN}
that timestep is limited by the constraints $dt\lesssim D_{min}/v_{max}^{2}$
and $dt\lesssim h_{e}^{2}/D_{max}$, where $D_{min}$ and $D_{max}$$[\mbox{m}^{2}\mbox{/s}]$
are the minimum and maximum diffusion coefficients, and $v_{max}$
is the maximum speed associated with the system.

\subsection{Conservation properties (discrete form of equations)}

Here, we demonstrate that the discretised system of equations \ref{eq:517.3}
has the same global conservation laws as the original continuous system. 

\subsubsection{Particle count conservation }

Integrating the discrete form of the continuity equation \ref{eq:517.3}
over volume, we obtain the rate of change of the total number of particles
in the system:

\begin{align*}
\dot{N} & =\overline{dV}^{T}*\dot{\overline{n}}\\
 & =\overline{dV}^{T}*\left\{ -\overline{\overline{\nabla}}\cdot\left(\overline{n}\,\,\overline{\mathbf{v}}\right)+\widehat{\overline{\nabla}}\cdot\left(\widehat{\zeta}\,\,\left(\overline{\widehat{\nabla}}\,\,\overline{n}\right)\right)\right\} 
\end{align*}
With boundary conditions $\overline{\mathbf{v}}_{\perp}|_{\Gamma}=\mathbf{0}$,
or $\overline{\mathbf{v}}|_{\Gamma}=\mathbf{0}$, identity \ref{eq:511.041}
determines that the first term is zero. The second term is always
zero, due to property \ref{eq:515.04}. Hence, total particle count
is conserved. Note that if no boundary conditions are explicitly applied
to density, $\overline{n}$ will automatically evolve to satisfy the
natural boundary condition $\left(\zeta\,\nabla_{\perp}n\right)|_{\Gamma}=0$,
as a consequence of the properties of the element-to-node divergence
operation (section \ref{subsec:Drn}).\\
\\

\subsubsection{Toroidal flux conservation\label{subsec:Toroidal-flux-conservation}}

Analogous to equation \ref{eq:480.5}, the discrete expression for
the rate of change of system toroidal flux follows from the discrete
expression for $\dot{\overline{f}}$ in equation \ref{eq:517.3} 
\begin{align*}
\dot{\Phi} & =\frac{1}{2\pi}\overline{dV}^{T}*\left\{ \dot{\overline{f}}\,/\,\overline{r}^{2}\right\} \\
 & =\frac{1}{2\pi}\overline{dV}^{T}*\left\{ -\overline{\overline{\nabla}}\cdot\left(\overline{f}\,\,\overline{\mathbf{v}}\,/\,\overline{r}^{2}\right)+\widehat{\overline{\nabla}}\cdot\left(\widehat{\mathbf{B}}_{\theta}\,\,\widehat{\mathbf{\omega}}\right)+\widehat{\overline{\nabla}}\cdot\left(\widehat{\mathbf{\eta}}\,\,\left(\overline{\widehat{\nabla}}\,\,\overline{f}\right)\,/\,\widehat{r}^{2}\right)\right\} 
\end{align*}
Once again, with boundary conditions $\overline{\mathbf{v}}_{\perp}|_{\Gamma}=\mathbf{0}$,
identity \ref{eq:511.041} determines that the first term is zero.
The second and third terms are zero, due to property \ref{eq:515.04}.
Hence, system toroidal flux is conserved. Note that in this case $\overline{f}$
will automatically evolve to satisfy the natural boundary condition
$\left(\mathbf{B}_{\theta\perp}\omega+\eta\left(\nabla_{\perp}\,f\,\right)/r^{2}\right)|_{\Gamma}=0$
if no boundary conditions are explicitly imposed on $\overline{f}$
. In combination with the boundary condition $\overline{\mathbf{v}}_{\perp}|_{\Gamma}=\mathbf{0}$,
this corresponds to having the poloidal component of the electric
field perpendicular to the boundary, the condition, in the case of
azimuthal symmetry, for a perfectly electrically conducting boundary
(see equation \ref{eq:480.51}). 

Simulations that include CT formation and compression are run with
the boundary conditions $\overline{v}_{\beta}|_{\Gamma}=0$ applied
explicitly to each velocity component, so that the automatically imposed
boundary condition for $f$, corresponding to the physical case of
a perfectly electrically conducting boundary, becomes $\left(\nabla_{\perp}\,f\,\right)|_{\Gamma}=0$.

As described in section \ref{sec:PHIconservation-with}, when part
of the computational boundary is modelled as an electrical insulator,
special care must be taken to define the explicitly applied boundary
condition for $f$ along the insulator boundary, in order to maintain
global toroidal flux conservation. In such cases, only the conducting
boundary regions are allowed to retain the naturally imposed boundary
conditions for $f$. 

\subsubsection{Angular momentum conservation\label{subsec:Angular-momentum-conservation} }

The rate of change of system angular momentum is $\dot{P}_{\phi}=\frac{\partial}{\partial t}\left(\int\left(\rho rv_{\phi}\right)dV\right)=m_{i}\int\left(rv_{\phi}\dot{n}+nr\dot{v}_{\phi}\right)dV$.
Here, it will be shown that system angular momentum is conserved for
the discrete model when the terms corresponding to density diffusion
in the discrete expressions for $\dot{n}$ and $\dot{v}_{\phi}$ are
neglected. In section \ref{subsec:Maintenance-of-momentum}, it will
be shown that angular momentum conservation can be maintained even
with the inclusion of these terms. In discrete form, the rate of change
of angular momentum is{\footnotesize{}
\begin{align}
\dot{P}_{\phi} & =m_{i}\,\overline{dV}^{T}*\left\{ -\left(\overline{r}\,\,\overline{v}_{\phi}\right)\,\,\left(\overline{\overline{\nabla}}\cdot(\overline{n}\,\,\overline{\mathbf{v}})\right)+\left(\overline{n}\,\,\overline{r}\right)\,\,\left\{ -\overline{\mathbf{v}}\cdot\left(\overline{\overline{\nabla}}\,\,(\overline{r}\,\,\overline{v}_{\phi})\right)\,/\,\overline{r}-\overline{\Pi}_{\phi}\,/\,\overline{\rho}+\widehat{\overline{W}}*\left(\widehat{\mathbf{B}}_{\theta}\cdot\left(\overline{\widehat{\nabla}}\,\,\overline{f}\right)\right)\,/\,\left(\mu_{0}\,\overline{r}\,\,\overline{\rho}\right)\right\} \right\} \label{eq:517.8}
\end{align}
}The first two terms here can be simplified to $-m_{i}\,\overline{dV}^{T}*\left\{ \left(\overline{r}\,\,\overline{v}_{\phi}\right)\,\,\left(\overline{\overline{\nabla}}\cdot(\overline{n}\,\,\overline{\mathbf{v}})\right)+\left(\overline{n}\,\,\overline{\mathbf{v}}\right)\cdot\left(\overline{\overline{\nabla}}\,\,(\overline{r}\,\,\overline{v}_{\phi})\right)\right\} $.
With boundary conditions $\overline{\mathbf{v}}_{\perp}|_{\Gamma}=\mathbf{0}$,
this combination vanishes due to identity \ref{eq:511.04}. With reference
to the definition for $\overline{\Pi}_{\phi}$ (equation \ref{eq:517.4}),
if no boundary conditions are explicitly applied to $\overline{v}_{\phi}$,
the properties of the element-to-node divergence operation will automatically
impose the natural boundary conditions $\left(\nabla_{\perp}\,\omega\,\right)|_{\Gamma}=0$,
and the third term in equation \ref{eq:517.8} will vanish due to
identity \ref{eq:515.04}. With reference to equation \ref{eq:516.1},
the fourth term in equation \ref{eq:517.8} can be expressed as 
\begin{align*}
 & \frac{1}{\mu_{0}}\,\widehat{dV}^{T}*\left\{ \widehat{\mathbf{B}}_{\theta}\cdot\left(\overline{\widehat{\nabla}}\,\,\overline{f}\right)\right\} \\
 & =\frac{2\pi}{\mu_{0}}\left(\widehat{s}\,\,\widehat{r}\right)^{T}*\left\{ \left(-\left(\overline{\widehat{Dz}}*\overline{\psi}\right)\,\,\left(\overline{\widehat{Dr}}*\overline{f}\right)+\left(\overline{\widehat{Dr}}*\overline{\psi}\right)\,\,\left(\overline{\widehat{Dz}}*\overline{f}\right)\right)\,/\,\widehat{r}\right\} \\
 & =\frac{2\pi}{\mu_{0}}\,\widehat{s}^{T}*\left\{ -\left(\overline{\widehat{Dz}}*\overline{\psi}\right)\,\,\left(\overline{\widehat{Dr}}*\overline{f}\right)+\left(\overline{\widehat{Dr}}*\overline{\psi}\right)\,\,\left(\overline{\widehat{Dz}}*\overline{f}\right)\right\} \\
 & =\frac{2\pi}{\mu_{0}}\left(-\left(\overline{\widehat{Dr}}*\overline{f}\right)^{T}*\left(\widehat{\widehat{S}}*\left(\overline{\widehat{Dz}}*\overline{\psi}\right)\right)+\left(\overline{\widehat{Dz}}*\overline{f}\right)^{T}*\left(\widehat{\widehat{S}}*\left(\overline{\widehat{Dr}}*\overline{\psi}\right)\right)\right)\\
 & \frac{2\pi}{\mu_{0}}\,\overline{f}^{T}*\left(\left(-\overline{\widehat{Dr}}^{T}*\widehat{\widehat{S}}*\overline{\widehat{Dz}}+\overline{\widehat{Dz}}^{T}*\widehat{\widehat{S}}*\overline{\widehat{Dr}}\right)*\overline{\psi}\right)\mbox{}
\end{align*}
and vanishes due to identity \ref{eq:505.003}, when the boundary
condition $\psi|_{\Gamma}=0$ is applied. Thus for angular momentum
conservation, no boundary conditions are applied on $v_{\phi}$, and
$\psi|_{\Gamma}$ must be set to zero.
\begin{figure}[H]
\centering{}\subfloat{\includegraphics[width=14cm,height=7cm]{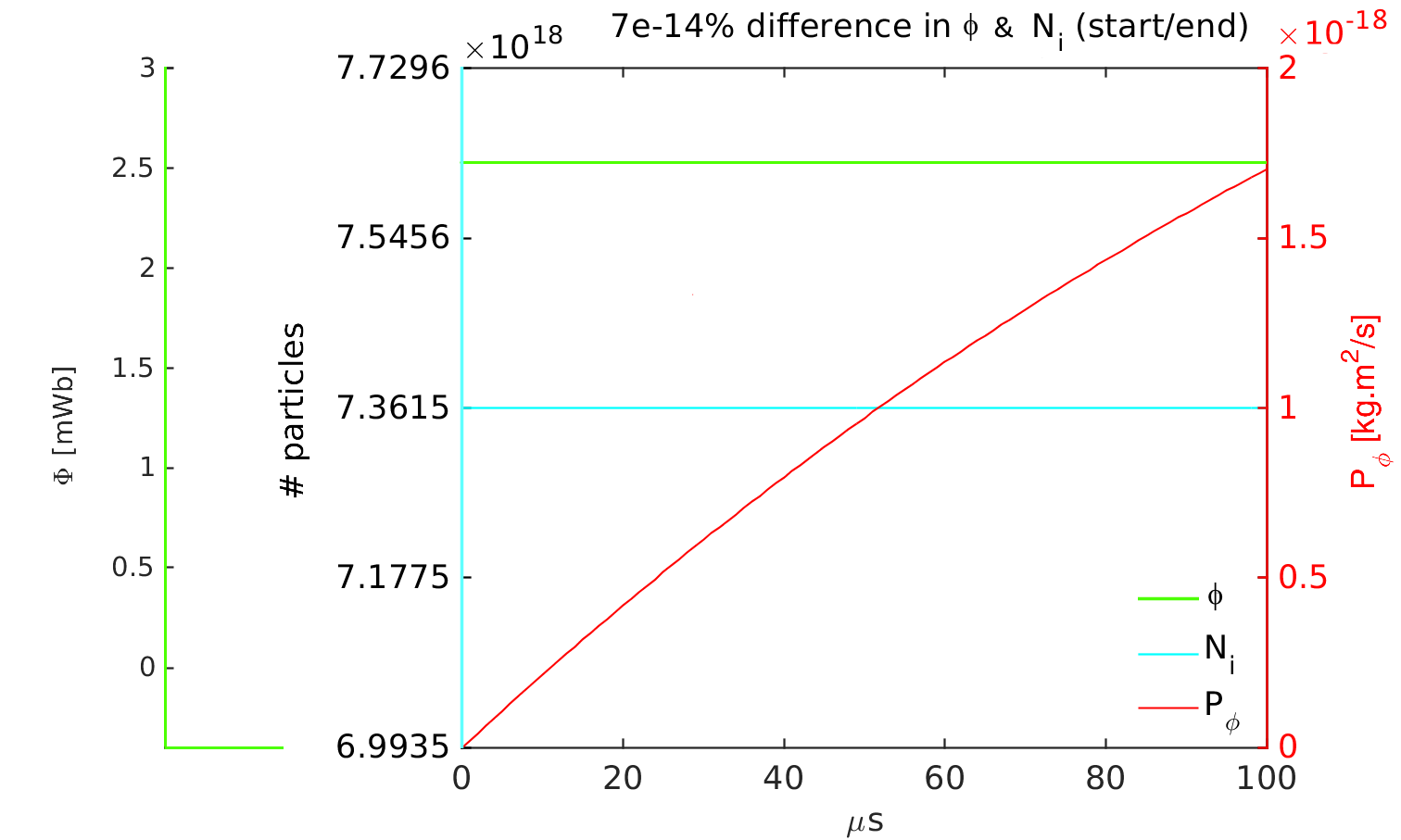}}\caption{\label{fig:Ni_phi_Pp cons}Illustration of particle, toroidal flux,
and angular momentum conservation }
\end{figure}
As shown in figure \ref{fig:Ni_phi_Pp cons}, the total number of
particles ($i.e.,$ \# ions, $N_{i})$, toroidal flux, and angular
momentum, calculated as the integrals over the computational domain,
is conserved to numerical precision for a $100\upmu$s simulation
that started from a Grad-Shafranov equilibrium, with explicitly applied
boundary conditions $v_{r}|_{\Gamma}=v_{z}|_{\Gamma}=0$, and $\psi|_{\Gamma}=0$.
As mentioned in section \ref{subsec:Conservation-properties-(continu},
angular momentum is not physical when there is friction between the
fluids and the boundary wall. The boundary condition $\left(\nabla_{\perp}\omega\right)|_{\Gamma}=0$
is physical only in unusual cases; for example, there may be no viscosity
at the wall due to rotation of the wall. Angular momentum is not conserved
in the simulations of CT formation and magnetic compression that will
be presented in section \ref{sec:Simulation-results,-and}, because
in those simulations, all velocity components are set to zero at the
boundary. 

\subsubsection{Energy conservation\label{subsec:Energy-conservation} }

Analogous to the continuous expression \ref{eq:481.1}, the discrete
expression for total system energy is

{\small{}
\begin{align}
\dot{U}_{Total} & =\dot{U}_{K}+\dot{U}_{Th}+\dot{U}_{M}\nonumber \\
 & =\frac{\partial}{\partial t}\,\left(\overline{dV}^{T}*\left\{ \frac{\overline{\rho}\,\,\overline{\mathbf{v}}^{2}}{2}+\frac{\overline{p}}{\gamma-1}+\frac{1}{2\mu_{0}}\,\overline{f}{}^{2}\,/\,\overline{r}{}^{2}\right\} +\widehat{dV}^{T}*\left\{ \frac{1}{2\mu_{0}}\left(\overline{\widehat{\nabla}}\,\,\overline{\psi}\right)^{2}\,/\,\widehat{r}{}^{2}\right\} \right)\nonumber \\
 & =\overline{dV}^{T}*\left\{ \frac{m_{i}\,\dot{\overline{n}}\,\,\overline{\mathbf{v}}^{2}}{2}+\overline{\rho}\,\,\overline{\mathbf{v}}\cdot\dot{\overline{\mathbf{v}}}+\frac{\dot{\overline{p}}}{\gamma-1}+\frac{1}{\mu_{0}}\,\overline{f}\,\,\dot{\overline{f}}\,/\,\overline{r}{}^{2}\right\} +\widehat{dV}^{T}*\left\{ \frac{1}{\mu_{0}}\,\left(\overline{\widehat{\nabla}}\,\,\overline{\psi}\right)\cdot\left(\overline{\widehat{\nabla}}\,\,\dot{\overline{\psi}}\right)\,/\,\widehat{r}{}^{2}\right\} \nonumber \\
 & =-\left[\overline{dV}^{T}*\left\{ \frac{1}{2}\overline{v}^{2}\,\,\left(\overline{\overline{\nabla}}\cdot(\overline{\rho}\,\,\overline{\mathbf{v}})\right)+\frac{1}{2}(\overline{\rho}\,\,\overline{\mathbf{v}})\cdot\left(\overline{\overline{\nabla}}\,\,\overline{v}^{2}\right)\right\} \right]\nonumber \\
 & \,\,\,-\left[\overline{dV}^{T}*\left\{ \overline{\mathbf{v}}\cdot\left(\overline{\overline{\nabla}}\,\,\overline{p}\right)+\frac{1}{\gamma-1}\left(\overline{\mathbf{v}}\cdot\left(\overline{\overline{\nabla}}\,\,\overline{p}\right)+\gamma\,\overline{p}\,\,\left(\overline{\overline{\nabla}}\cdot\overline{\mathbf{v}}\right)\right)\right\} \right]\nonumber \\
 & \,\,\,-\left[\frac{1}{\mu_{0}}\,\overline{dV}^{T}*\left\{ \left(\overline{\mathbf{v}}\cdot\left(\overline{\overline{\nabla}}\,\,\overline{\psi}\right)\right)\,\,\left(\overline{\overline{\Delta^{^{*}}}}\,\,\overline{\psi}\right)\,/\,\overline{r}{}^{2}\right\} +\widehat{dV}^{T}*\left\{ \left(\overline{\widehat{\nabla}}\,\,\overline{\psi}\right)\cdot\left(\overline{\widehat{\nabla}}\,\,\left(\overline{\mathbf{v}}\cdot\left(\overline{\overline{\nabla}}\,\,\overline{\psi}\right)\right)\right)\,/\,\widehat{r}{}^{2}\right\} \right]\nonumber \\
 & \,\,\,-\left[\frac{1}{\mu_{0}}\,\overline{dV}^{T}*\left\{ \left(\overline{f}\,\,\overline{\mathbf{v}}\,/\,\overline{r}^{2}\right)\cdot\left(\overline{\overline{\nabla}}\,\,\overline{f}\right)+\overline{f}\,\,\left(\overline{\overline{\nabla}}\cdot\left(\overline{f}\,\,\overline{\mathbf{v}}\,/\,\overline{r}^{2}\right)\right)\right\} \right]\nonumber \\
 & \,\,\,+\left[\frac{1}{\mu_{0}}\,\overline{dV}^{T}*\left\{ \overline{\omega}\,\,\left(\widehat{\overline{W}}*\left(\widehat{\mathbf{B}}_{\theta}\cdot\left(\overline{\widehat{\nabla}}\,\,\overline{f}\right)\right)\right)+\overline{f}\,\,\left(\widehat{\overline{\nabla}}\cdot\left(\widehat{\mathbf{B}}_{\theta}\,\,\widehat{\mathbf{\omega}}\right)\right)\right\} \right]\nonumber \\
 & \,\,\,+\left[\frac{1}{\mu_{0}}\,\overline{dV}^{T}*\left\{ \overline{\eta}\,\,\left(\left(\overline{\overline{\Delta^{^{*}}}}\,\,\overline{\psi}\right)\,/\,\overline{r}\right)^{2}\right\} +\widehat{dV}^{T}*\left\{ \left(\overline{\widehat{\nabla}}\,\,\overline{\psi}\right)\cdot\left(\overline{\widehat{\nabla}}\,\,\left(\overline{\eta}\,\,\left(\overline{\overline{\Delta^{^{*}}}}\,\,\overline{\psi}\right)\right)\right)\,/\,\widehat{r}{}^{2}\right\} \right]\nonumber \\
 & \,\,\,+\left[\frac{1}{\mu_{0}}\,\overline{dV}^{T}*\left\{ \widehat{\overline{W}}*\left(\left(\widehat{\mathbf{\eta}}\,\,\left(\overline{\widehat{\nabla}}\,\,\overline{f}\right)\,/\,\widehat{r}^{2}\right)\cdot\left(\overline{\widehat{\nabla}}\,\,\overline{f}\right)\right)+\overline{f}\,\,\left(\widehat{\overline{\nabla}}\cdot\left(\widehat{\mathbf{\eta}}\,\,\left(\overline{\widehat{\nabla}}\,\,\overline{f}\right)\,/\,\widehat{r}^{2}\right)\right)\right\} \right]\nonumber \\
 & \,\,\,-\left[\overline{dV}^{T}*\left\{ \widehat{\overline{\nabla}}\cdot\left(\widehat{\mathbf{q}}_{i}+\widehat{\mathbf{q}}_{e}\right)\right\} \right]\nonumber \\
 & \,\,\,-\left[\overline{dV}^{T}*\left\{ \overline{\mathbf{v}}\cdot\overline{\boldsymbol{\varPi}}-\overline{Q}_{\pi}\right\} \right]\label{eq:518}
\end{align}
}Note that each set of square brackets has a counterpart in equation
\ref{eq:481.1}, the continuous form of $\dot{u}_{Total}$. The poloidal
magnetic energy is expressed in terms of the element-centered gradient
of $\overline{\psi}$, as required for consistency with the definition
of the second order $\overline{\overline{\Delta^{^{*}}}}$ differential
operator. Here, it will be shown how, using the various mimetic properties
of the differential operators, that the terms in each set of square
brackets cancel when appropriate boundary conditions are applied,
leading to total system energy conservation.

The terms in the first set of square brackets here represents the
contribution to $\dot{U}_{K}$ due to advection, and vanish, with
boundary conditions $\mathbf{v_{\perp}}|_{\Gamma}=\mathbf{0}$, due
to equation \ref{eq:511.04} (where $\overline{\mathbf{P}}=\overline{\rho}\,\,\overline{\mathbf{v}}$
and $\overline{U}=\overline{v}^{2}$). The terms in the second of
square brackets represents the contribution to $\dot{U}_{Th}$ from
compressional heating and, with the same boundary conditions, also
vanish due to equation \ref{eq:511.04} (with $\overline{\mathbf{P}}=\overline{\mathbf{v}}$
and $\overline{U}=\overline{p}$). 

The terms in the third, fourth and fifth sets of square brackets in
equation \ref{eq:518} represent the contribution to $\dot{U}_{K}$
that arises from the discrete forms of the components of $(\mathbf{J}\times\mathbf{B})$
in combination with the ideal (non-resistive) part of $\dot{U}_{M}$.
Using identity \ref{eq:515.72} to expand the $\overline{\overline{\Delta^{^{*}}}}$
operator, the terms in the third set of square brackets can be expressed
as 

\[
-\frac{1}{\mu_{0}}\left[\overline{dV}^{T}*\left\{ \left(\overline{\mathbf{v}}\cdot\left(\overline{\overline{\nabla}}\,\,\overline{\psi}\right)\right)\,\,\left(\widehat{\overline{\nabla}}\cdot\left(\left(\overline{\widehat{\nabla}}\,\,\overline{\psi}\right)\,/\,\widehat{r}^{2}\right)\right)\right\} +\widehat{dV}^{T}*\left\{ \left(\left(\overline{\widehat{\nabla}}\,\,\overline{\psi}\right)\,/\,\widehat{r}^{2}\right)\cdot\left(\overline{\widehat{\nabla}}\,\,\left(\overline{\mathbf{v}}\cdot\left(\overline{\overline{\nabla}}\,\,\overline{\psi}\right)\right)\right)\right\} \right]
\]
These terms cancel due to identity \ref{eq:515.031}, where $\widehat{\mathbf{P}}=\left(\left(\overline{\widehat{\nabla}}\,\,\overline{\psi}\right)\,/\,\widehat{r}^{2}\right)$
and $\overline{U}=\overline{\mathbf{v}}\cdot\left(\overline{\overline{\nabla}}\,\,\overline{\psi}\right)$. 

With boundary condition $\mathbf{v\perp}|_{\Gamma}=\mathbf{0}$, the
terms in the fourth set of square brackets cancel due to identity
\ref{eq:511.04}. Using equation \ref{eq:516}, where $\overline{Q}=\overline{\omega}=\overline{v}_{\phi}\,/\,\overline{r}$,
the terms in the fifth set of square brackets can be expressed as 

\begin{align*}
\frac{1}{\mu_{0}}\,\left[\widehat{dV}^{T}*\left\{ \left(\widehat{\mathbf{B}}_{\theta}\,\,\widehat{\mathbf{\omega}}\right)\cdot\left(\overline{\widehat{\nabla}}\,\,\overline{f}\right)\right\} +\overline{dV}^{T}*\left\{ \overline{f}\,\,\left(\widehat{\overline{\nabla}}\cdot\left(\widehat{\mathbf{B}}_{\theta}\,\,\widehat{\mathbf{\omega}}\right)\right)\right\} \right]
\end{align*}
This combination cancels due to identity \ref{eq:515.031}, where
$\widehat{\mathbf{P}}=\left(\widehat{\mathbf{B}}_{\theta}\,\,\widehat{\mathbf{\omega}}\right)$
and $\overline{U}=\overline{f}$. 

The terms in the sixth set of square brackets represent the rate of
increase of thermal energy due to ohmic heating from toroidal currents
(first term), in combination with the rate of decrease of magnetic
energy associated with poloidal field, due to resistive decay of the
toroidal currents (second term). As is true in the continuous case,
these terms are balanced in the discrete case. Using identity \ref{eq:515.72},
the combination may be expressed as 

{\footnotesize{}
\begin{align*}
 & \frac{1}{\mu_{0}}\left[\overline{dV}^{T}*\left\{ \left(\overline{\eta}\,\,\left(\overline{\overline{\Delta^{^{*}}}}\,\,\overline{\psi}\right)\right)\,\,\left(\widehat{\overline{\nabla}}\cdot\left(\left(\overline{\widehat{\nabla}}\,\,\overline{\psi}\right)\,/\,\widehat{r}{}^{2}\right)\right)\right\} +\widehat{dV}^{T}*\left\{ \left(\left(\overline{\widehat{\nabla}}\,\,\overline{\psi}\right)\,/\,\widehat{r}{}^{2}\right)\cdot\left(\overline{\widehat{\nabla}}\,\,\left(\overline{\eta}\,\,\left(\overline{\overline{\Delta^{^{*}}}}\,\,\overline{\psi}\right)\right)\right)\right\} \right]
\end{align*}
}and vanishes due to identity \ref{eq:515.031}, where $\widehat{\mathbf{P}}=\left(\left(\overline{\widehat{\nabla}}\,\,\overline{\psi}\right)\,/\,\widehat{r}{}^{2}\right)$
and $\overline{U}=\overline{\eta}\,\,\left(\overline{\overline{\Delta^{^{*}}}}\,\,\overline{\psi}\right)$. 

The terms in the seventh set of square brackets in equation \ref{eq:518}
represent the rate of increase of thermal energy due to ohmic heating
from poloidal currents (first term), in combination with the rate
of decrease of magnetic energy associated with toroidal magnetic field,
due to resistive decay of the poloidal currents (second term). Again,
these terms cancel, which is physically representative. Referring
to equation \ref{eq:516.1}, the combination vanishes due to identity
\ref{eq:515.031}, where $\widehat{\mathbf{P}}=\left(\widehat{\mathbf{\eta}}\,\,\left(\overline{\widehat{\nabla}}\,\,\overline{f}\right)\,/\,\widehat{r}^{2}\right)$
and $\overline{U}=\overline{f}$. 

With reference to the definitions of the discrete forms for the thermal
flux $\widehat{\overline{\nabla}}\cdot\widehat{\mathbf{q}}_{\alpha}$
, where $\widehat{\mathbf{q}}_{\alpha}$ is defined in equation \ref{eq:517.7},
it can be seen how the terms in the eighth set of square brackets
in equation \ref{eq:518} vanish due to equation \ref{eq:515.04}.
Boundary conditions that are explicitly applied to the pressure fields
break computational-domain-energy conservation, and enable thermal
losses in accordance with the thermal conduction model and the values
explicitly applied to $n_{0},\,\left(\chi_{\parallel\alpha}\right)|_{\Gamma}$
and $\left(\chi_{\perp\alpha}\right)|_{\Gamma}$, but thermal energy
fluxes through the boundary may be systematically accounted for. 

To deal with the viscosity related terms in the final set of square
brackets in equation \ref{eq:518}, the substitutions $\overline{\mathbf{P}_{1}}=\left(\overline{v}_{r}\,/\,\overline{r}\right)\widehat{\mathbf{r}}$
, $\overline{\mathbf{P}_{2}}=\left(\overline{v}_{z}\,/\,\overline{r}\right)\widehat{\mathbf{z}}$
and $\overline{\mathbf{P}_{3}}=\left(\overline{v}_{z}\,/\,\overline{r}\right)\widehat{\mathbf{r}}+\left(\overline{v}_{r}\,/\,\overline{r}\right)\widehat{\mathbf{z}}$
are made, and the expansion of terms, using equations \ref{eq:517.4},
\ref{eq:517.5}, and \ref{eq:516.1} is

{\small{}
\begin{align*}
 & \overline{dV}^{T}*\left\{ -\overline{\mathbf{v}}\cdot\overline{\boldsymbol{\varPi}}+\overline{Q}_{\pi}\right\} \\
 & =2\left[\overline{dV}^{T}*\left\{ \overline{\mathbf{P}_{1}}\cdot\left(\widehat{\overline{\nabla}}\,\,\left(\widehat{\mathbf{\mu}}\,\,\widehat{r}^{2}\,\,\left(\overline{\widehat{\nabla}}\cdot\overline{\mathbf{P}_{1}}\right)\right)\right)\right\} +\widehat{dV}^{T}*\left\{ \widehat{\mathbf{\mu}}\,\,\widehat{r}^{2}\,\,\left(\overline{\widehat{\nabla}}\cdot\overline{\mathbf{P}_{1}}\right)^{2}\right\} \right]\\
 & +2\left[\overline{dV}^{T}*\left\{ \overline{\mathbf{P}_{2}}\cdot\left(\widehat{\overline{\nabla}}\,\,\left(\widehat{\mathbf{\mu}}\,\,\widehat{r}^{2}\,\,\left(\overline{\widehat{\nabla}}\cdot\overline{\mathbf{P}_{2}}\right)\right)\right)\right\} +\widehat{dV}^{T}*\left\{ \widehat{\mathbf{\mu}}\,\,\widehat{r}^{2}\,\,\left(\overline{\widehat{\nabla}}\cdot\overline{\mathbf{P}_{2}}\right)^{2}\right\} \right]\\
 & +\left[\overline{dV}^{T}*\left\{ \overline{\mathbf{P}_{3}}\cdot\left(\widehat{\overline{\nabla}}\,\,\left(\widehat{\mathbf{\mu}}\,\,\widehat{r}^{2}\,\,\left(\overline{\widehat{\nabla}}\cdot\overline{\mathbf{P}_{3}}\right)\right)\right)\right\} +\widehat{dV}^{T}*\left\{ \widehat{\mathbf{\mu}}\,\,\widehat{r}^{2}\,\,\left(\overline{\widehat{\nabla}}\cdot\overline{\mathbf{P}_{3}}\right)^{2}\right\} \right]\\
 & -\frac{2}{3}\left[\overline{dV}^{T}*\left\{ \overline{\mathbf{v}}\cdot\left(\widehat{\overline{\nabla}}\,\,\left(\widehat{\mathbf{\mu}}\,\,\left(\overline{\widehat{\nabla}}\cdot\overline{\mathbf{v}}\right)\right)\right)\right\} +\widehat{dV}^{T}*\left\{ \widehat{\mathbf{\mu}}\,\,\left(\overline{\widehat{\nabla}}\cdot\overline{\mathbf{v}}\right)^{2}\right\} \right]\\
 & +\left[\overline{dV}^{T}*\left\{ \underline{\mathbf{\omega}}\,\,\left(\widehat{\overline{\nabla}}\cdot\left(\widehat{\mathbf{\mu}}\,\,\widehat{r}^{2}\,\,\left(\overline{\widehat{\nabla}}\,\,\overline{\omega}\right)\right)\right)\right\} +\widehat{dV}^{T}*\left\{ \widehat{\mathbf{\mu}}\,\,\widehat{r}^{2}\,\,\left(\overline{\widehat{\nabla}}\,\,\overline{\omega}\right)^{2}\right\} \right]\\
 & -2\left[\overline{dV}^{T}*\left\{ \overline{\mu}\,\,\left(\overline{v}_{r}\,/\,\overline{r}\right)^{2}-\overline{\mu}\,\,\left(\overline{v}_{r}\,/\,\overline{r}\right)^{2}\right\} \right]
\end{align*}
}The terms in the first to fourth sets of square brackets here vanish
due to equation \ref{eq:515.03}, where, for example, for the first
set of square brackets, $\overline{\mathbf{P}}=\overline{\mathbf{P}_{1}}$
and $\widehat{U}=\widehat{\mathbf{\mu}}\,\,\widehat{r}^{2}\,\,\left(\overline{\widehat{\nabla}}\cdot\overline{\mathbf{P}_{1}}\right)$.
The terms in the fifth set of square brackets vanish due to equation
\ref{eq:515.031} where $\overline{U}=\overline{\omega}$ and $\widehat{\mathbf{P}}=\widehat{\mathbf{\mu}}\,\,\widehat{r}^{2}\,\,\left(\overline{\widehat{\nabla}}\,\,\overline{\omega}\right)$,
and the terms in the final set obviously cancel. Thus, the volume-integrated
rate of increase of thermal energy due to viscous heating is balanced
by the volume-integrated rate of decrease of kinetic energy due to
viscous dissipation, which is also a property of the physical system.
\-\-
\begin{figure}[H]
\centering{}\includegraphics[scale=0.6]{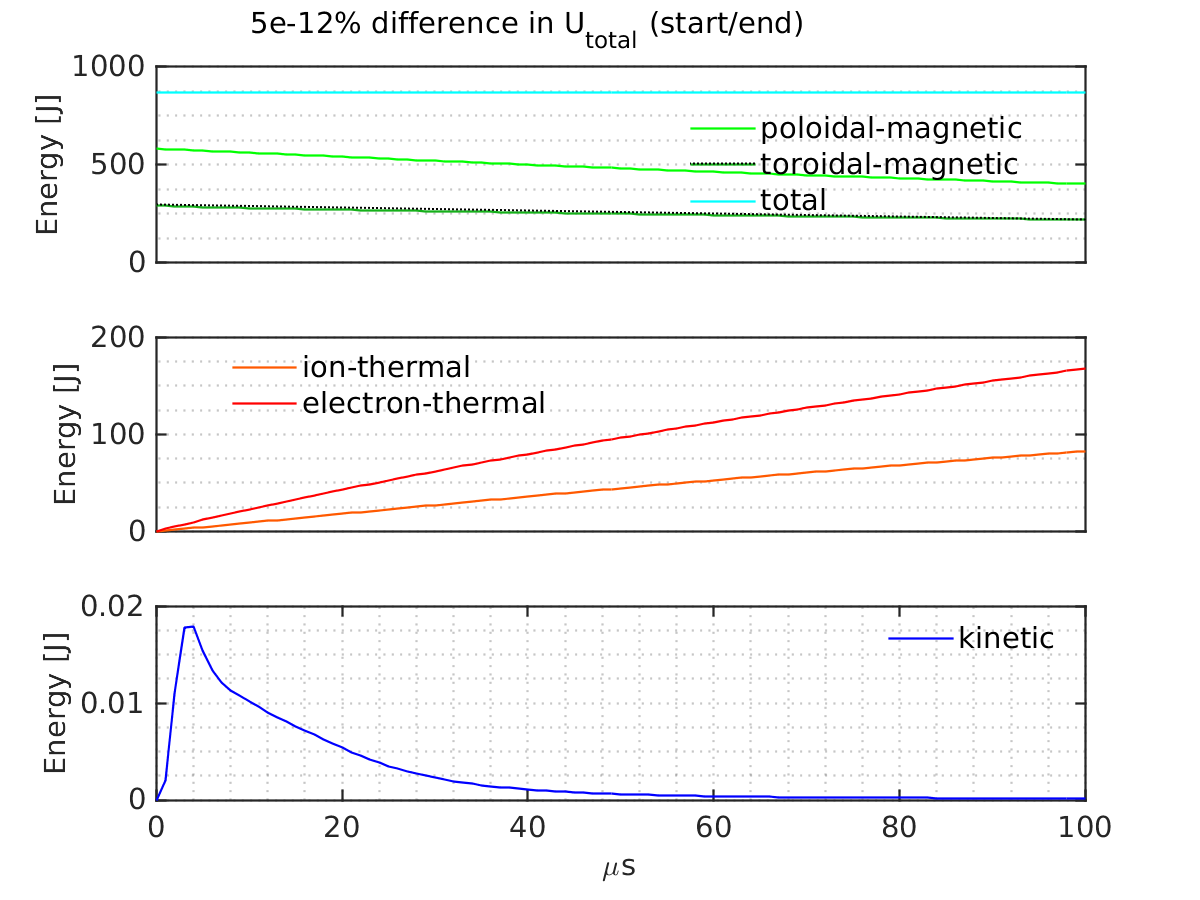}\caption{\label{fig:Energy-partition-of}Energy partitions for Grad-Shafranov
solution }
\end{figure}
Figure \ref{fig:Energy-partition-of} shows the partition of energy,
and how total energy is conserved approximately to machine precision
for a $100\upmu$s simulation with moderate time-step and mesh resolution,
using Runge-Kutta second order time-stepping, that started from a
Grad-Shafranov equilibrium. The only explicitly applied boundary conditions
are $v_{r}|_{\Gamma}=v_{z}|_{\Gamma}=\psi|_{\Gamma}=0$. In contrast,
simulations presented in section \ref{sec:Simulation-results,-and}
include magnetic levitation and magnetic compression, and are run
with explicitly applied boundary conditions for $\psi$ that are determined
by experiment. This enables inward electromagnetic energy fueling
and outward electromagnetic losses (Poynting flux). In those simulations,
boundary conditions are explicitly applied to the pressure fields
to enable thermal losses in accordance with the thermal conduction
model, and boundary conditions for $f$ are explicitly applied only
on the parts of the domain boundary that represent insulating regions. 

\section{Experiment overview\label{sec:Experiment-overview}}

\begin{figure}[H]
\centering{}\subfloat[]{\centering{}\includegraphics[scale=0.4]{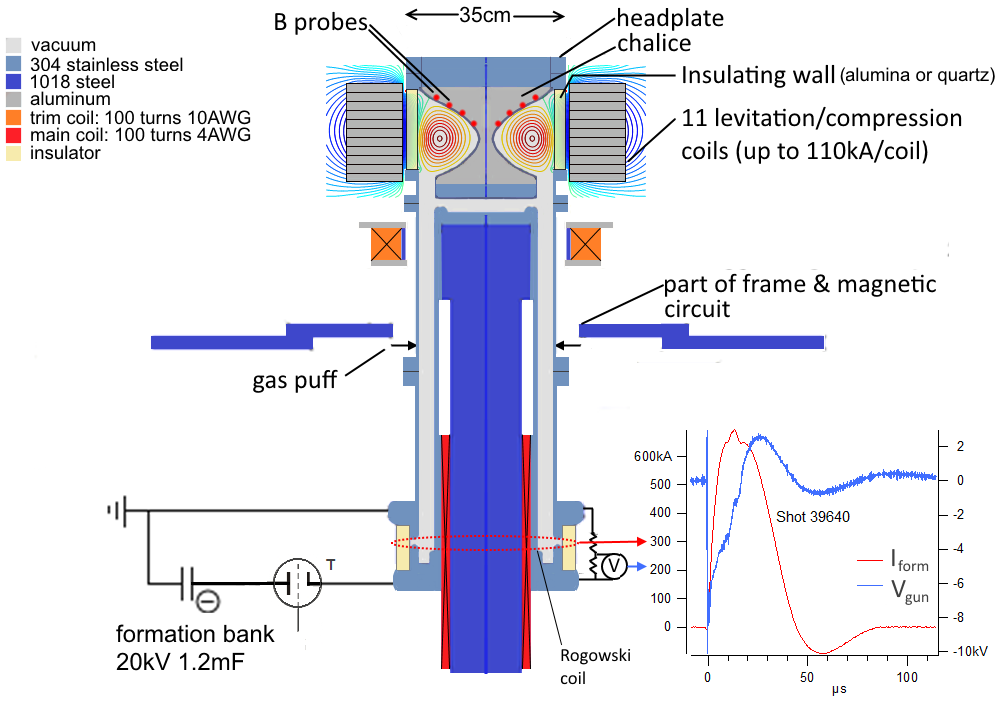}}\hfill{}\subfloat[]{\centering{}\includegraphics[scale=0.5]{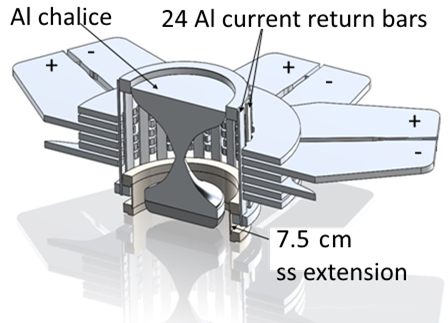}}\caption{\label{fig:Machine-Schematic}Machine schematic (a) and view of CT
confinement region surrounded by the stack of six coil plates (b) }
\end{figure}
The structure and experiment design of the SMRT magnetic compression
device can be seen in figure \ref{fig:Machine-Schematic}(a), with
CT and levitation $\psi$ (poloidal flux) contours from an equilibrium
model superimposed, and with a schematic of the formation circuit.
Measurements of formation current $(I_{form}(t))$, and voltage across
the formation electrodes $(V_{gun}(t))$ are also indicated. Note
that the principal materials used in the machine construction, and
some key components, are indicated by the color-key at the top left
of the figure. Eleven levitation/compression coils are located outside
the insulating tube around the CT containment region. Figure \ref{fig:Machine-Schematic}(b)
indicates, for the original configuration with six coils, the coils,
inner tungsten-coated aluminum chalice flux conserver, stainless steel
extension, and aluminum return current bars that carry axial current
outside the insulating wall that is depicted in figure \ref{fig:Machine-Schematic}(a).
A more thorough overview of the experiment can be found in \cite{thesis,exppaper}.
The sequence of machine operation is as follows:\\

\vspace{1cm}
\begin{tabular}[t]{lll}
(i) & $t\sim-3$s  & Main coil is energised with steady state ($\sim4\mbox{s}$ duration)\tabularnewline
 &  & current, $I_{main}\sim70$A \tabularnewline
(ii) & $t=t_{gas}\sim-400\upmu$s  & Gas is injected into vacuum \tabularnewline
(iii) & $t=t_{lev}\sim-400\upmu\mbox{s}\rightarrow-40\upmu\mbox{s}$  & Levitation banks are fired ($I_{lev}\sim180$kA)\tabularnewline
(iv) & $t=0$s  & Formation banks are fired ($I_{form}\sim700$kA)\tabularnewline
(v) & $t=t_{comp}\sim40\upmu\mbox{s}\rightarrow150\upmu\mbox{s}$  & Compression banks are fired ($I_{comp}\sim1.2$MA)\tabularnewline
\end{tabular}

\section{Model for CT formation, levitation and magnetic compression\label{sec:Model-for-CT}}

\subsection{Boundary conditions\label{subsec:Boundary-conditions}}

No boundary condition is explicitly applied to density, so that the
natural boundary condition is automatically imposed: 
\[
\left(\nabla_{\perp}n\right)|_{\Gamma}=0
\]
Simulations of CT formation, levitation and magnetic compression presented
here are run with all velocity components set to zero on the boundary
$i.e.,$ the no-slip condition for an impermeable boundary: 
\[
\mathbf{v}|_{\Gamma}=\mathbf{0}
\]
Boundary conditions for pressures follow from equations \ref{eq:479.31},
with temperatures $T_{i}$, $T_{e}$ set to approximately zero at
the boundary: 
\[
T_{i}|_{\Gamma}=T_{e}|_{\Gamma}=0.02\mbox{eV}\approx20^{o}\mbox{C}
\]
Toroidal currents in the main, and levitation/compression coils constitute
sources of poloidal flux, which can be included in the model with
the application of appropriate boundary conditions for $\psi$. After
assigning the material properties, solution frequencies, and peak
coil currents pertaining to either the main stuffing field, the levitation
field or the compression field, and running a $\mbox{\mbox{FEMM}}$
\cite{FEMM} model for the relevant configuration, the values of $\psi_{main},\,\psi_{lev},\,\mbox{and }\psi_{comp}$
at the boundary nodes are extracted to file. In the experiment, the
main (stuffing) coils are energised around three seconds before the
formation banks are fired, so the stuffing field can soak entirely
through all materials in the model and can be modelled as being steady
state, therefore the solution frequency input to $\mbox{\mbox{FEMM}}$
for obtaining $\psi_{main}$ is 0Hz. The solution frequency for obtaining
$\psi_{lev}$ is taken as a function of $|t_{lev}|$, the time over
which the levitation field is allowed to soak into the machine before
firing the formation bank. $|t_{lev}|$ is taken as the quarter-period
of the current waveform input to $\mbox{\mbox{FEMM}}$, so that the
solution frequency is $\sim5$kHz for typical values, from experiment,
of $|t_{lev}|=50\upmu\mbox{s}$. Similarly, as the rise time of the
compression current is $\sim20\upmu$s, the solution frequency for
obtaining $\psi_{comp}$ is taken as $\sim12.5$kHz. Field diffusion
times in the most electrically conducting parts of the machine are
much longer than the timescales associated with the levitation and
compression fields. Magnetic boundary conditions are applied at all
boundary nodes, but levitation boundary values are significant only
on the levitation/compression coil boundaries and on the boundaries
of the regions representing the stainless steel above and below the
insulating wall, while compression boundary values are significant
only on the coil boundaries, and are approximately zero at other boundary
regions. 
\begin{figure}[H]
\subfloat[$\widetilde{I}_{lev}(t)$]{\begin{centering}
\includegraphics[width=7cm,height=5cm]{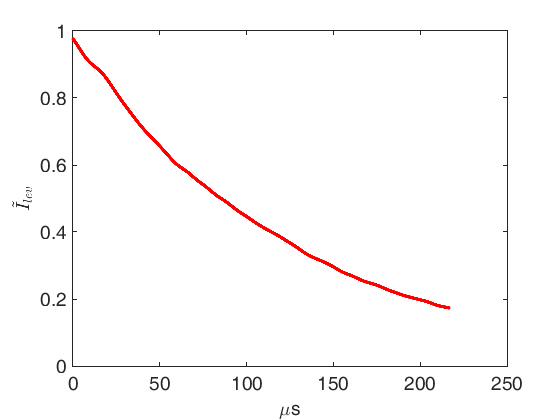}
\par\end{centering}
}\hfill{}\subfloat[$\widetilde{I}_{comp}(t)$]{\centering{}\includegraphics[width=7cm,height=5cm]{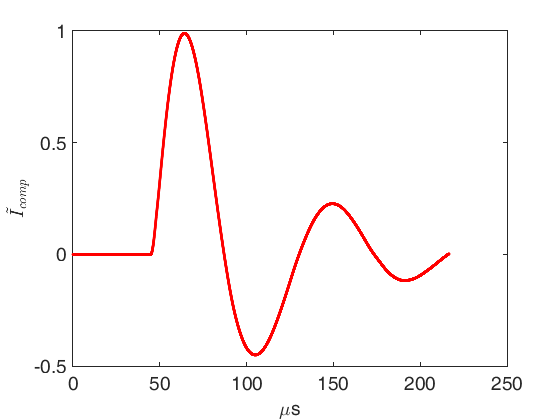}}\caption{Normalised levitation (a) and compression (b) current signals\label{fig:Ilev_comp_sim}}
\end{figure}
 The boundary conditions for $\psi_{main}(\mathbf{r})$ are held constant
over time, while boundary values obtained for $\psi_{lev}(\mathbf{r},\,t)$
and $\psi_{comp}(\mathbf{r},\,t)$, which pertain to the peak levitation/compression
currents, are scaled over time according to the experimentally measured
waveforms for $I_{lev}(t)$ and $I_{comp}$(t): 
\begin{equation}
\psi(\mathbf{r},\,t)|_{\Gamma}=\psi_{main}(\mathbf{r})|_{\Gamma}+\psi_{lev}(\mathbf{r},\,t)|_{\Gamma}+\psi_{comp}(\mathbf{r},\,t)|_{\Gamma}\label{eq:519}
\end{equation}
Figure \ref{fig:Ilev_comp_sim}(a) shows $\widetilde{I}_{lev}(t)$,
the normalised levitation current signal, as measured with Rogowski
coils during the experiment; similarly, figure \ref{fig:Ilev_comp_sim}(b)
shows $\widetilde{I}_{comp}(t)$. For this simulation, $t_{comp}=45\upmu$s,
and the simulation was run until around $220\upmu$s. Note that $\widetilde{I}_{lev}<1$
at $t=0$, because $t_{lev}=-50\upmu$s , $i.e.,$ in the experiment,
the levitation capacitor banks were fired $50\upmu$s before the formation
banks. With a rise time of $\sim40\upmu$s, the levitation current
peaks at $t\sim-10\upmu$s, when $\widetilde{I}_{lev}=1$. Typical
peak levitation currents were $\sim180$kA total, peak compression
currents were up to 1.2 MA.

No boundary conditions are explicitly applied to $f$ on boundary
regions representing electrical conductors. In combination with the
boundary conditions for velocity, the natural boundary conditions
\[
\left(\nabla_{\perp}f\right)|_{\Gamma}=0
\]
are automatically imposed, corresponding to the condition for a perfectly
electrically conducting boundaries. Since $f=rB_{\phi}=\mu_{0}I_{\theta}/2\pi$
(from Ampere's law when $\frac{\partial}{\partial\phi}=0$), the boundary
condition $\left(\nabla_{\perp}f\right)|_{\Gamma}=0$ implies that
any poloidal currents that flow into or out of the wall ($e.g.,$
radial intra-electrode formation current or crowbarred shaft current
during magnetic compression) flow perpendicular to the wall, $i.e.,$
$\mathbf{E}$ has no component parallel to the boundary. Since no
currents can flow in insulating regions, $f$ must be spatially constant
on the part of the boundary representing the insulating wall surrounding
the CT confinement region (see section \ref{sec:Vacuum-field-in}
below). As detailed below in section \ref{sec:PHIconservation-with},
the value of that constant, consistent with conservation of toroidal
flux in the combined domains, is calculated at each timestep and explicitly
applied as a boundary condition for $f$ on the interface shared by
the plasma and insulating domains. 

\subsection{Vacuum field in insulating region\label{sec:Vacuum-field-in}}

As described in \cite{thesis,exppaper}, levitated CT lifetime and
temperature, and the recurrence rate of shots in which the CT conserved
most of its poloidal flux during compression, was increased significantly
when the original six levitation/compression coil configuration was
replaced with an eleven coil configuration. It is thought that the
improvement in levitated CT lifetime and temperature was due to reduced
levels of levitation field displacement and interaction between the
plasma and insulating wall during the CT formation process, leading
to consequent reductions in sputtering, plasma impurity concentration,
and radiative cooling.

\begin{figure}[H]
\subfloat[Full mesh]{\includegraphics[scale=0.5]{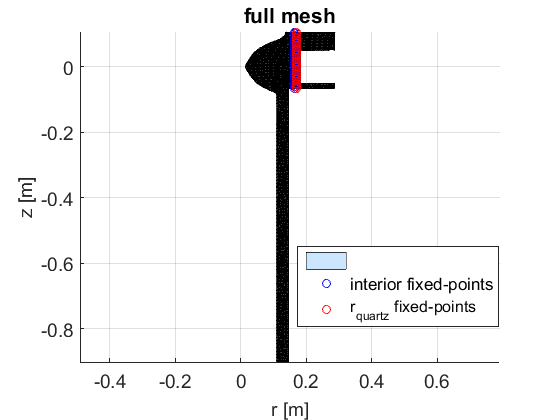}}
\hfill{}\subfloat[Full mesh (close up)]{\includegraphics[scale=0.5]{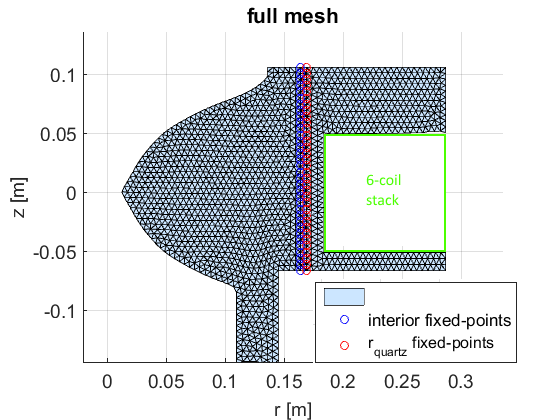}} 

\subfloat[Plasma mesh with indication of magnetic probe locations]{\includegraphics[scale=0.5]{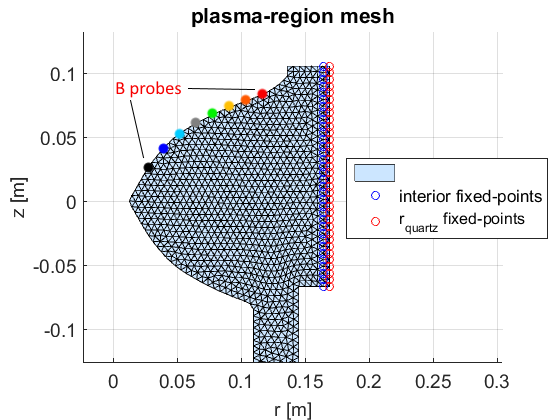}

} \hfill{}\subfloat[Insulator mesh]{\includegraphics[scale=0.5]{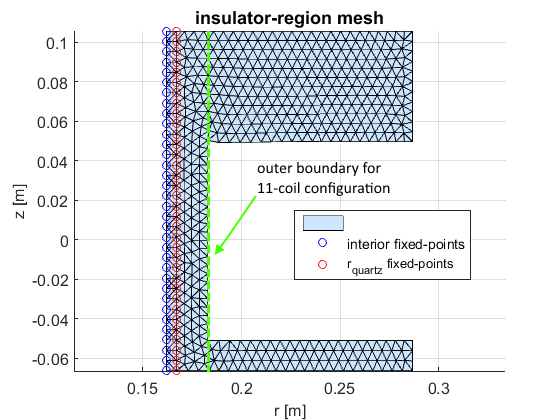}}

\caption{\label{fig:Grid-arrangement-with}Computational mesh with insulating
region, 6-coil configuration. Note that the magnetic probe locations
are indicated in (c).}
\end{figure}
In order to model the interaction of plasma with the insulating wall
during the CT formation process, a vacuum field should be solved for
in the insulating area between the inner radius of the insulating
outer wall and the levitation/compression coils. The insulating area
includes the wall itself as well as the air between the wall and the
coil-stack, and the air above and below the coil-stack. If this area
is included in the domain in which the plasma fields are solved, then
unphysical currents will be allowed to flow in insulating regions.
To solve for a vacuum field in the insulating region and couple it
to the plasma fields, the computational grid for the combined plasma
and insulating domain is split in two, as shown in figure \ref{fig:Grid-arrangement-with},
in which the original six levitation/compression coil configuration
is being modelled. The six coil stack representation is specifically
indicated in figure \ref{fig:Grid-arrangement-with}(b). For the 11-coil
configuration, in which the stack of coils extends along the entire
height of the insulating wall, the main part of the mesh (figure \ref{fig:Grid-arrangement-with}(c))
is unchanged, while the parts of the insulator mesh located above
and below the 6-coil stack are deleted, as indicated in figure \ref{fig:Grid-arrangement-with}(d).
To couple between the plasma and vacuum solutions at each timestep,
we use the two vertical rows of shared fixed mesh points indicated
in the figure to exchange boundary values along the plasma/insulating
interface. 

No currents can flow in a vacuum, so 
\begin{equation}
\Delta^{^{*}}\psi_{v}=0\label{eq:111.1}
\end{equation}
where $\psi_{v}$ represents $\psi$ in a vacuum (or insulator). Defining
$\overline{\psi}_{v}=\overline{\psi}_{v\Gamma}+\overline{\psi}_{vi}$,
where $\overline{\psi}_{v\Gamma}$ and $\overline{\psi}_{vi}$ have
the values of $\psi$ set to zero at the nodes corresponding to internal
/ boundary nodes respectively, the discrete form of equation \ref{eq:111.1}
leads to $\overline{\overline{\Delta^{^{*}}}}\,\,\overline{\psi}_{v}=\overline{\overline{\Delta^{^{*}}}}\,\,\overline{\psi}_{v\Gamma}+\overline{\overline{\Delta_{0}^{^{*}}}}\,\,\overline{\psi}_{vi}=\overline{0}$
, so that: 
\begin{equation}
\overline{\psi}_{vi}=-\left(\overline{\overline{\Delta_{0}^{^{*}}}}\right)^{-1}*\left(\overline{\overline{\Delta^{^{*}}}}\,\,\overline{\psi}_{v\Gamma}\right)\label{eq:112}
\end{equation}
The discrete operator $\overline{\overline{\Delta_{0}^{^{*}}}}$ is
designed to produce the same results as $\overline{\overline{\Delta^{^{*}}}}$
at internal nodes if the boundary values of the operand field are
set to zero. The boundary values of $\overline{a}$ are preserved
in the operation $\overline{c}=\overline{\overline{\Delta_{0}^{^{*}}}}\,\,\overline{a}$,
and do not contribute to the values calculated for elements corresponding
to internal nodes in \textbf{$\overline{c}$}. Defining $\overline{b}\,[N_{n}\times1]$
as the logical vector defining the indexes of the boundary nodes,
$\overline{\overline{\Delta_{0}^{^{*}}}}$ is defined as $\overline{\overline{\Delta_{0}^{^{*}}}}=\overline{\overline{\Delta^{^{*}}}},$
with $\overline{\overline{\Delta_{0}^{^{*}}}}(\overline{b},\,:)=\overline{\overline{\Delta_{0}^{^{*}}}}(:,\,\overline{b})=0$,
and $\overline{\overline{\Delta_{0}^{^{*}}}}(\overline{b},\,\overline{b})=1$.
$\overline{\overline{\Delta^{^{*}}}}$ and $\overline{\overline{\Delta_{0}^{^{*}}}}$
operators based on the geometries of the plasma domain and the insulating
domain were constructed. With zero initial simulated plasma density
everywhere, $\psi$ in the entire domain is initially a vacuum field,
so additional $\overline{\overline{\Delta^{^{*}}}}$ and $\overline{\overline{\Delta_{0}^{^{*}}}}$
operators based on the geometry of the combined domain are required
in order to construct the initial vacuum $\psi$ solution. For formation
simulations, initially the time-relevant boundary values for $\psi$
are applied at the boundary points of the combined mesh, according
to equation \ref{eq:519}, and the vacuum field in the combined domain
is obtained using equation \ref{eq:112}. The sequence of subsequent
steps followed at each timestep to couple the solutions is as follows:
\begin{enumerate}
\item After evolving $\psi_{plasma}$ at the nodes in the plasma domain
according to equations \ref{eq:517.3}, boundary conditions $\psi(\mathbf{r},\,t)|_{\Gamma}$
are applied separately to the boundary points of the plasma and insulator
meshes.
\item The values at the nodes along the left boundary (\textquotedbl interior
points\textquotedbl{} - blue circles in figure \ref{fig:Grid-arrangement-with}(d))
of the insulator mesh are overwritten using the equivalent values
from $\psi_{plasma}$. 
\item Equation \ref{eq:112} is used to calculate $\psi_{v}$ in the insulating
area.
\item Finally, the values at the outer right boundary of the plasma-grid
(red circles in figure \ref{fig:Grid-arrangement-with}(c)) are overwritten
using the equivalent values from the insulating region solution $\psi_{v}.$
\end{enumerate}

\subsection{Formation Simulation\label{subsec:Formation-Simulations}}

Physically, CT formation in a magnetized Marshall gun is achieved
as a result of $\mathbf{J}_{r}\times\mathbf{B}_{\phi}$ forces acting
on plasma, where $\mathbf{J}_{r}$ is the radial formation current
density across the plasma between the electrodes, and $\mathbf{B}_{\phi}$
is the toroidal field due to the axial formation current in the electrodes.
Open stuffing magnetic field lines that are resistively pinned to
the electrodes, and frozen into the conducting plasma, are advected
with the plasma by the $\mathbf{J}_{r}\times\mathbf{B}_{\phi}$ force,
into the containment region, where they reconnect to form CT closed
flux surfaces. Simulated formation is initiated with the addition
of toroidal flux below the physical locations of the gas puff valves
- initial radial formation current is assumed to occur at the z-coordinate
of the valves, where the gas density is initially highest.

When this external source of toroidal flux is included, the continuous
form of the expression for the time-rate of change of $f$ is 
\begin{equation}
\dot{f}(\mathbf{r},\,t)=r^{2}\,\nabla\cdot\bigg(-\frac{f}{r^{2}}\mathbf{v}+\omega\mathbf{B}+\frac{\eta}{r^{2}}\nabla f\bigg)+\dot{f}_{form}(z,\,t)\label{eq:144.2}
\end{equation}
In practice, the expression for $f_{form}$ is added to the existing
value for $f$ at the beginning of each timestep, so that the natural
toroidal flux conserving boundary condition ($(\nabla_{\perp}f)|_{\Gamma}=0$)
is maintained on electrically conducting boundaries. Conservation
of the system's intrinsic toroidal flux implies that the total system
flux is equal to the initial flux plus the flux added due to $f_{form}(z,\,t)$,
at each time.

Integrating Faraday's law in the poloidal plane over the area defining
the formation current path (up the outer gun electrode as far as the
region with the initially highest gas concentration at the z-coordinate
of the valves, through the plasma across the intra-electrode gap,
and down the inner electrode) we have: 
\begin{align}
 & \int\,\nabla\times\mathbf{E}_{\theta}(\mathbf{r},t)\cdot d\mathbf{S}=-\int\,\dot{\mathbf{B}}_{\phi}(\mathbf{r},t)\cdot d\mathbf{S}\nonumber \\
 & \Rightarrow\int\,\mathbf{E}_{\theta}(\mathbf{r},t)\cdot d\mathbf{l}=V(t)=-\dot{\Phi}_{form}(t)\nonumber \\
 & \Rightarrow\Phi_{form}(t)=-\intop_{0}^{t}V(t')\:dt'\label{eq:145}
\end{align}
Here, 
\begin{equation}
V(t)=V_{gun}(t)+I_{form}(t)\,R(t)\label{eq:145.001}
\end{equation}
where $V_{gun}(t)$ is the voltage measured across the formation electrodes
and $I_{form}(t)\,R(t)$ is the resistive voltage drop along the formation
current path. $\mathbf{E}_{\theta}(\mathbf{r},t)$ is the poloidal
electric field that is established through the application of $V_{gun}(t)$,
and $\dot{\mathbf{B}}_{\phi}(\mathbf{r},t)$ and $\dot{\Phi}_{form}(t)$
are the time-rates of change of the toroidal magnetic field and toroidal
flux induced in the area bounded by the path along which formation
current driven by $\mathbf{E}_{\theta}(\mathbf{r},t)$ flows. 

For the simplifying assumption of a fixed formation current path of
constant inductance, $-V(t)=\dot{\Phi}_{form}(t)=L\dot{I}_{form}(t)$,
where $I_{form}(t)$ is the formation current and $L$ is the inductance
of the formation current path. $R(t)$ in equation \ref{eq:145.001}
is a resistance that depends on the resistivity of the metal along
the formation current path, and also on the resistivity of the plasma
between the gun electrodes. In addition, it is expected that when
plasma, which is advected upwards during the formation process, occupies
the gap between the chalice and the inner gun electrode (see figure
\ref{fig:Machine-Schematic}(a)), that part of the formation current
will flow on a path up the outer electrode, through the aluminum bars
located outside the  insulating wall (see figure \ref{fig:Machine-Schematic}(b)),
down the chalice, across the plasma in the gap and down the inner
electrode. Physically, $R(t)$ includes the contribution of the resistance
of the plasma in the gap shortly after initiation of the formation
process. Note that the chalice ($i.e.,$ the inner flux conserver
indicated in figure \ref{fig:Machine-Schematic}(b)) and inner electrode
are modelled as a continuous conductor - we do not have sufficient
information, for example the voltage measured across the gap, to properly
model the effect of the gap. For simplicity, we assume here that $R$
is constant in time. Thus the expression for the voltage measured
across the formation electrodes is

\begin{align}
-V_{gun}(t) & =L\dot{I}_{form}(t)+I_{form}(t)R\nonumber \\
 & =\dot{\Phi}_{form}(t)+\frac{\Phi_{form}(t)}{\tau_{LR}}\label{eq:150.1}
\end{align}
where $\tau_{LR}=L/R$ is the $LR$ time determining the formation
current decay rate. This time constant can be estimated from the e-folding
time of the toroidal field that is experimentally measured at the
probes indicated in figure \ref{fig:Machine-Schematic}(a). The toroidal
field at the probes is a measure of the crow-barred shaft current
that flows poloidally around the machine after formation. The value
of $\tau_{LR}$ varies depending on the formation current path and
resistances of the intra-electrode plasma and the plasma in the gap
between the chalice and the gun inner electrode, so it can vary from
shot to shot depending on plasma conditions. For relatively long-lived
magnetically levitated CTs, $\tau_{LR}\sim90\upmu$s. Reduced values
of $\tau_{LR}$ can be chosen as a simulation input when shorter-lived
CTs ($e.g.,$ for simulations with relatively high thermal diffusion
coefficients, or low flux CTs associated with reduced $V_{form}$
and $I_{main}$) are being modelled. Equation \ref{eq:150.1} can
be solved as 
\begin{equation}
\Phi_{form}(t)=-e^{-\frac{t}{\tau_{LR}}}\intop_{0}^{t}V_{gun}(t')\,e^{\frac{t'}{\tau_{LR}}}\:dt'\label{eq:145.01}
\end{equation}
By definition,
\begin{equation}
\Phi_{form}(t)=\int\frac{f_{form}(z,t)}{r}dr\;dz\label{eq:146}
\end{equation}
 $f_{form}(z,t)$ may be expressed as
\begin{equation}
f_{form}(z,t)=A_{form}(t)\,g_{form}(z,t)\label{eq:147}
\end{equation}
where $A_{form}(t)$ determines the amplitude of $f_{form}(z,t)$,
and $g_{form}(z,t)$ is a geometric profile that determines where
plasma formation current flows. The smooth logistic profile
\begin{equation}
g_{form}(z,t)=\frac{e^{m\,z_{I}(t)}}{e^{m\,z_{I}(t)}+e^{m\,z}}\label{eq:145-1}
\end{equation}
defines $g_{form}(z,t)$, where $m\sim40$ defines the profile slope,
and $z_{I}(t)$ is the $z$ coordinate at which the greatest concentration
of radial formation current flows between the machine electrodes.
For simplicity, we neglect the time dependence in $g_{form}$ by replacing
$z_{I}(t)$ with $z_{gp}=-0.43\mbox{m}$, the $z$ coordinate of the
machine's gas puff valves, around which the greatest concentration
of radial formation current is expected to flow.
\begin{figure}[H]
\subfloat[]{\begin{centering}
\includegraphics[width=7cm,height=5cm]{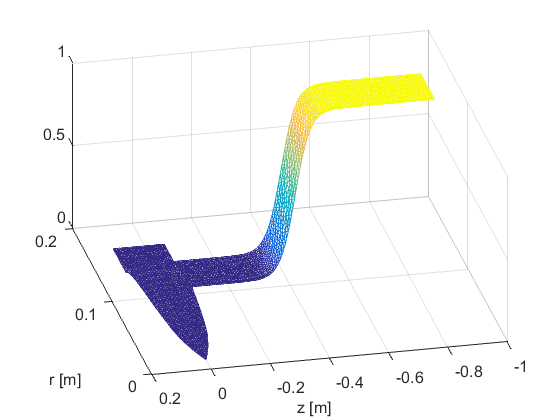}
\par\end{centering}
}\hfill{}\subfloat[]{\centering{}\includegraphics[width=7cm,height=5cm]{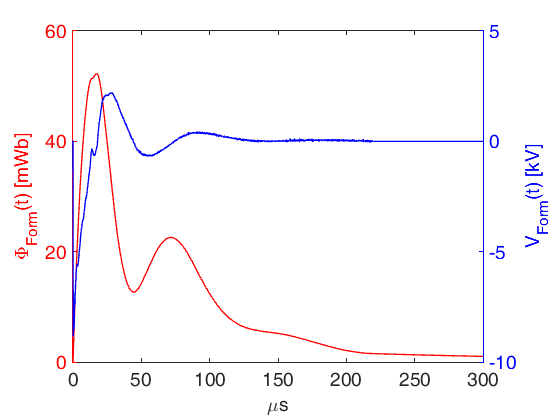}}\caption{\label{fig:gform_Phiform_Vform}Logistic function $g_{form}(z)$ profile
for formation flux input (a), and measured $V_{gun}(t)$ with calculated
$\Phi_{form}(t)$ (b)}
\end{figure}
The logistic function formation profile is shown in figure \ref{fig:gform_Phiform_Vform}(a).
As outlined in \cite{thesis}, it is possible to reproduce the measured
formation current over the bulk of the simulated formation process,
when the time dependence of $g_{form}$ is included. Figure \ref{fig:gform_Phiform_Vform}(b)
shows the experimentally measured formation voltage signal (right
axis) which is numerically integrated over time to evaluate $\Phi_{form}(t)$
(left axis) as defined by equation \ref{eq:145.01}.\\
\\
Combining equations \ref{eq:146} and \ref{eq:147}, we have:
\begin{equation}
\Phi_{form}(t)=A_{form}(t)\,\int\frac{g_{form}(z)}{r}dr\;dz\label{eq:148}
\end{equation}
Together with equation \ref{eq:145.01}, this gives: 
\begin{equation}
A_{form}(t)=\frac{-e^{-\frac{t}{\tau_{LR}}}\intop_{0}^{t}V_{gun}(t')\,e^{\frac{t'}{\tau_{LR}}}\:dt'}{\int\frac{g_{form}(z)}{r}dr\;dz}\label{eq:149}
\end{equation}
Combining equations \ref{eq:147} and \ref{eq:149} yields the expression
for $f_{form}(z,t)$, which is used in equation \ref{eq:144.2} to
simulate the formation process: 
\begin{equation}
f_{form}(z,t)=\left(-e^{-\frac{t}{\tau_{LR}}}\intop_{0}^{t}V_{gun}(t')\,e^{\frac{t'}{\tau_{LR}}}\:dt'\right)\left(\frac{g_{form}(z)}{\int\frac{g_{form}(z)}{r}dr\;dz}\right)\label{eq:150}
\end{equation}

\subsection{$\Phi$ conservation with insulating wall\label{sec:PHIconservation-with}}

Integrating toroidal field over the area of the total domain (combined
plasma and insulating wall regions) in the poloidal plane, we have:

\begin{equation}
\Phi_{tot}(t)=\overset{\underbrace{\int_{\Omega_{P}}\frac{f_{P}(\mathbf{r},t)}{r}dr\,dz+\int_{\Omega_{I}}\frac{f_{I}(t)}{r}dr\,dz}}{\Phi_{PI}(t)}+\overset{\underbrace{\int_{\Omega_{P}+\Omega_{I}}\frac{f_{form}(z,t)}{r}dr\,dz}}{\Phi_{form}(t)}\label{eq:176}
\end{equation}
The subscripts $P$ and $I$ refer to plasma and insulating wall regions.
Note that the insulating wall area includes the insulating wall only,
and not the insulating (air) regions above and below the coil stack
as depicted in figure \ref{fig:Grid-arrangement-with}. Poloidal currents
in the aluminum bars located outboard of the insulating wall (depicted
in figure \ref{fig:Machine-Schematic}(b)) add toroidal flux to the
system only in the insulating wall and plasma regions. Since no currents
can flow in insulators, $f$ must be spatially constant in the insulating
wall inside aluminum bars. The flux-conserving natural boundary condition
that is imposed on $f_{P}$ (which is evolved according to equation
\ref{eq:144.2}) must be overwritten on the interface between the
insulating wall and the plasma domain. The second term in equation
\ref{eq:176} can be re-expressed as $\int_{\Omega_{I}}\frac{f_{I}(t)}{r}dr\;dz=f_{I}(t)\,h_{I}\,\mbox{ln}(r_{out}/r_{in})$,
where $h_{I}$ is the height of the rectangular cross-section of the
insulating wall, and $r_{out}$ and $r_{in}$ are the outer and inner
radii of the wall. Initial system toroidal flux is zero, and $\Phi_{PI}(t)$
is to be conserved by design, so the system's only source of toroidal
flux is $\Phi_{form}(z,t)$:
\begin{align}
 & \Phi_{tot}(t)=\Phi_{form}(t)\Rightarrow\Phi_{PI}(t)=0\nonumber \\
 & \Rightarrow\int_{\Omega_{P}}\frac{f_{P}(\mathbf{r},t)}{r}dr\;dz+\int_{\Omega_{I}}\frac{f_{I}(t)}{r}dr\;dz=0\nonumber \\
 & \Rightarrow\stackrel[i=1]{N_{n}}{\Sigma}\left(\frac{f_{P_{i}}(t)\,s_{i}}{3r_{i}}\right)+f_{I}(t)\left(h_{I}\,\mbox{ln}(r_{out}/r_{in})\right)=0\label{eq:140}
\end{align}
Here, $N_{n}$ is the number of nodes in the non-insulating part of
the domain ($i.e.,$ the plasma domain) in which the MHD equations
are solved, $s_{i}/3$ is the area, and $r_{i}$ is the radial coordinate,
associated with node $i$. The summation is over all nodes in the
plasma domain, which includes the fixed-point nodes along the inner
radius of the insulating wall, at which we want to evaluate $f_{I}$.
Equation \ref{eq:140} can be solved for the constant $f_{I}$ if
$f_{P_{i}}$ is temporarily set to zero at the fixed-point nodes along
the inner wall of the insulator (see figure \ref{fig:Grid-arrangement-with}),
so that $f_{P}\rightarrow f_{P0}$. Equation \ref{eq:140} is modified
as:
\begin{equation}
\stackrel[i=1]{N_{n}}{\Sigma}\left(\frac{f_{P0_{i}}(t)\,s_{i}}{3r_{i}}\right)+\left(\stackrel[j=1]{N_{int}}{\Sigma}\left(\frac{s_{j}}{3r_{j}}\right)+h_{I}\,\mbox{ln}(r_{out}/r_{in})\right)f_{I}(t)=0\label{eq:141}
\end{equation}
Here, $\stackrel[j=1]{N_{int}}{\Sigma}$ implies summation over the
interface fixed-point nodes along the inner insulating wall. In equation
\ref{eq:141}, $\tilde{L}_{ins}=h_{I}\,\mbox{ln}(r_{out}/r_{in})$
is related to the inductance of the insulating wall's cross-sectional
area (recall that the inductance of a co-axial cable is $L_{coaxial}=\frac{\mu_{0}}{2\pi}l\,\mbox{ln}(r_{out}/r_{in})$,
where $l$ is the length of the coaxial cable and $r_{out}/r_{in}$
are the cable's outer and inner radii), and $\tilde{L}_{int\Delta}=\stackrel[j=1]{N_{int}}{\Sigma}\left(\frac{s_{j}}{3r_{j}}\right)$
is related to the inductance of the area of the parts of the triangular
elements in the plasma domain that are associated with the interface
fixed point nodes. The resulting expression for $f_{I}$ is:
\begin{equation}
f_{I}(t)=\:\frac{-1}{\tilde{L}_{ins}+\tilde{L}{}_{int\Delta}}\stackrel[i=1]{N_{n}}{\Sigma}\left(\frac{f_{P0_{i}}(t)s_{i}}{3r_{i}}\right)\label{eq:142}
\end{equation}
The constant $f_{I}$ is calculated at each timestep and is applied
as a boundary condition for $f$ on the interface shared by the plasma
domain and the insulating domain, resulting in conservation of total
toroidal flux in the combined domains. 
\begin{figure}[H]
\centering{}\subfloat{\centering{}\includegraphics[width=7cm,height=5cm]{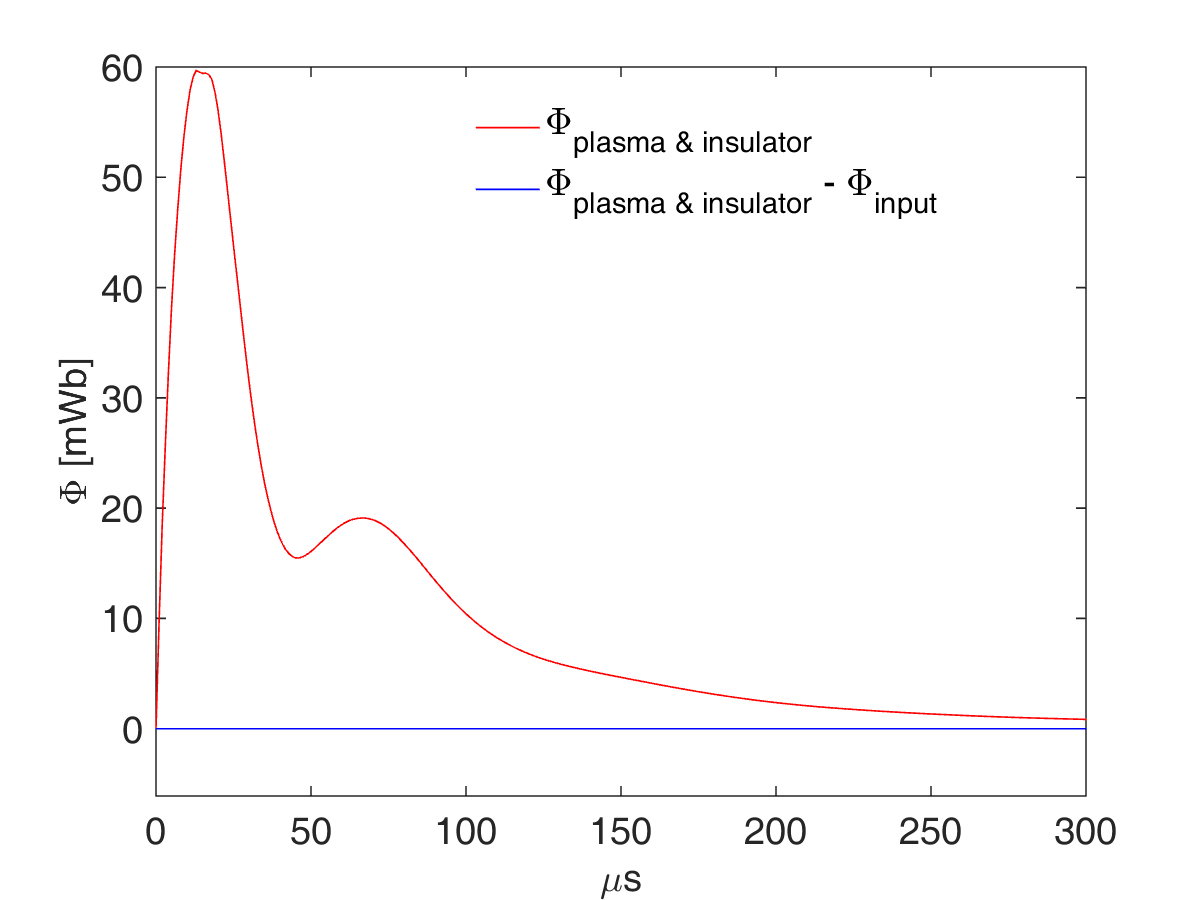}}\caption{\label{fig:phi_cons} Indication of toroidal flux conservation for
MHD simulation of CT formation }
\end{figure}
Now, with the natural boundary condition $\left(\nabla_{\perp}f\right)|_{\Gamma}=0$
imposed at all boundary points in the plasma domain (in combination
with the explicitly applied boundary conditions $\overline{\mathbf{v}}|_{\Gamma}=\mathbf{0}$),
and with equation \ref{eq:142} used to overwrite the values for $f_{P}$
on the interface, intrinsic toroidal flux conservation is ensured
in the combined computational domain. As indicated in figure \ref{fig:phi_cons}
toroidal flux is conserved over time for simulations of CT formation
- total flux less input (formation-associated) flux is zero at all
times. 

\section{Simulation results, and comparison with experiment\label{sec:Simulation-results,-and}}

In this section, results from an MHD simulation of CT formation into
a levitation field, and magnetic compression will be presented and
compared with experimental results. The models described in section
\ref{sec:Model-for-CT} have been incorporated. All simulation and
experiment results presented will correspond to the machine configuration
with eleven levitation/compression coils, with the exception that
simulations of CT formation will be compared across the two configuration
in section \ref{subsec:Simulated-plasma-wall-interactio}. Note that
code verification studies have been performed to confirm convergence
of the field solutions with decreasing element size or timestep. In
addition, the Grad-Shafranov equilibrium and $q$ profile solutions
have been benchmarked against the corresponding solutions from the
well-established Corsica code \cite{corsica code}. Verification studies,
and the methods implemented to evaluate the equilibrium and $q$ profile
solutions within the DELiTE framework, are presented in \cite{thesis}. 

\begin{figure}[H]
\subfloat[]{\raggedright{}\includegraphics[scale=0.4]{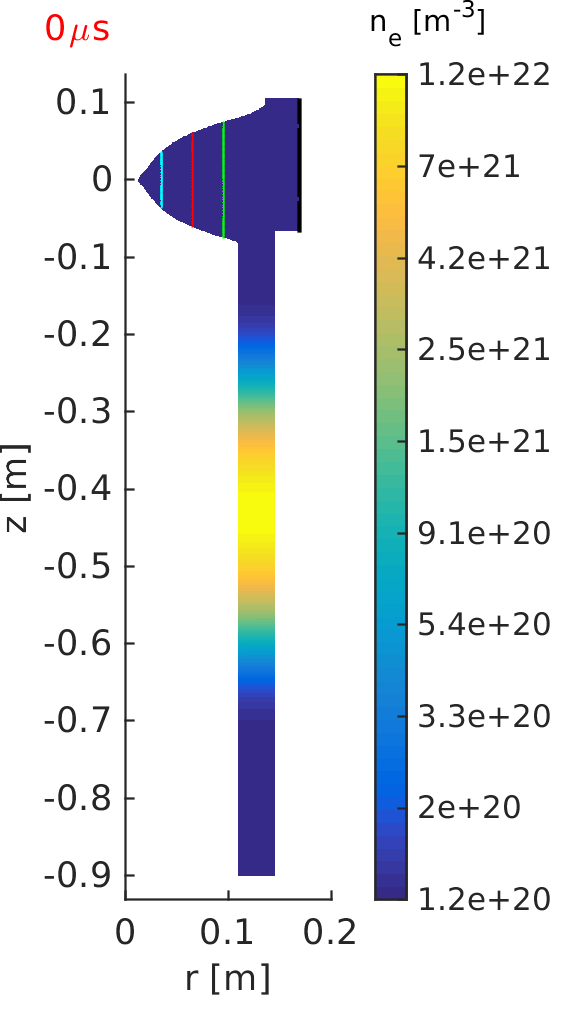}}\hfill{}\subfloat[]{\raggedright{}\includegraphics[scale=0.4]{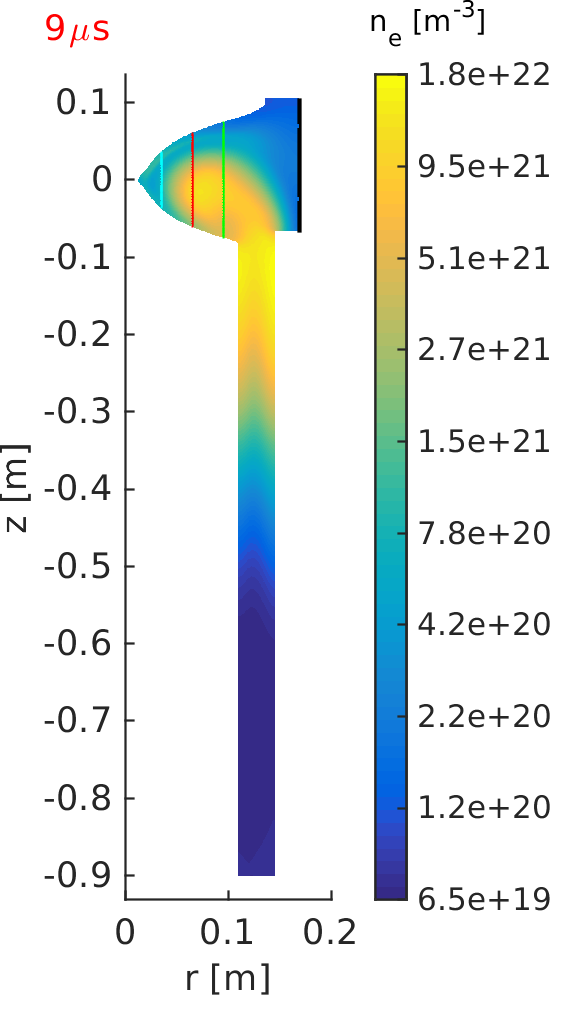}}\hfill{}\subfloat[]{\raggedright{}\includegraphics[scale=0.4]{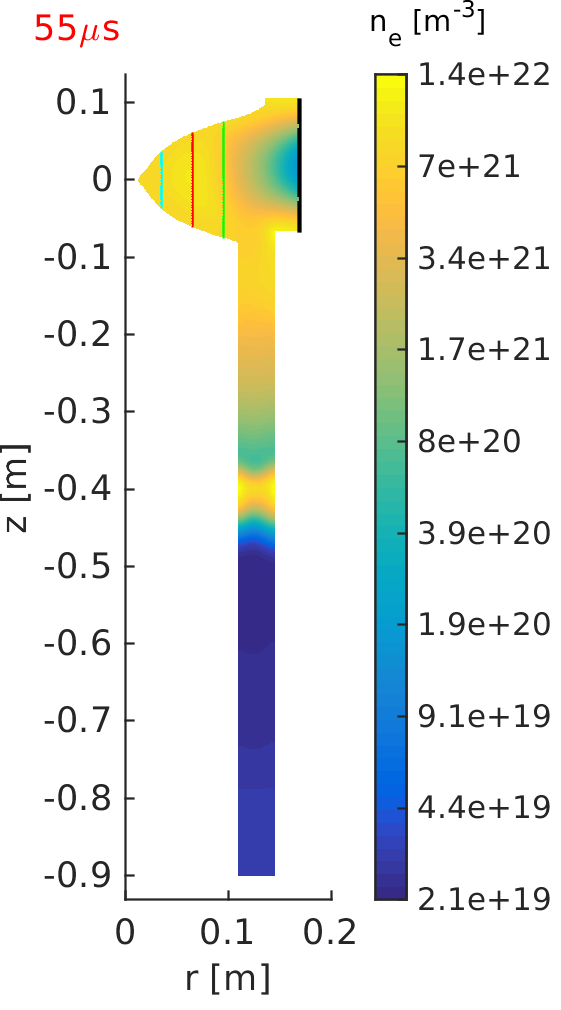}}\hfill{}\subfloat[]{\raggedright{}\includegraphics[scale=0.4]{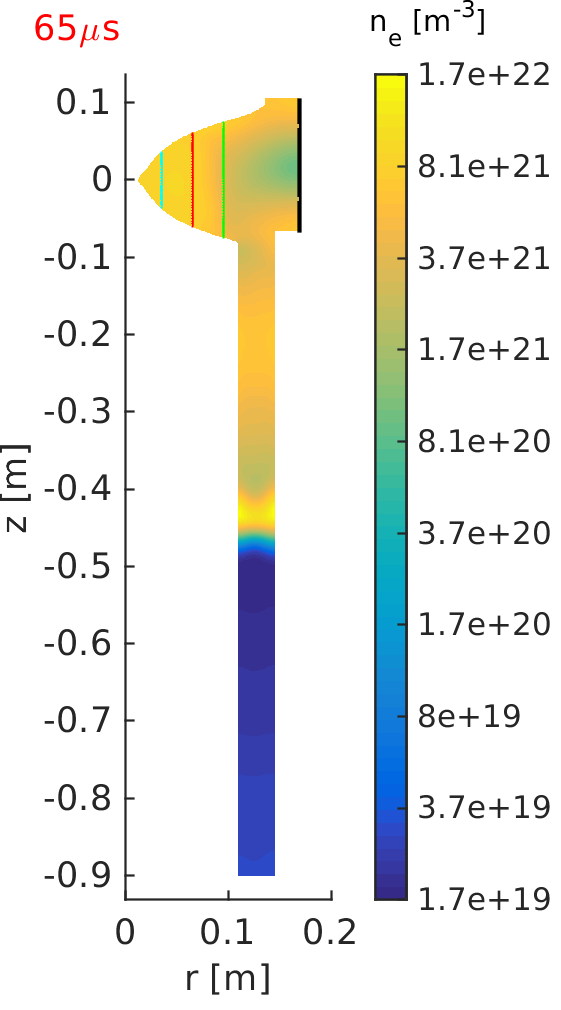}}

\caption{\label{fig: n_11coils}Electron density profiles at $0\upmu$s (a),
$9\upmu$s (b), $55\upmu$s (c), $65\upmu$s (d). The vertical blue,
red and green chords represent the lines of sight of the interferometer
measurements, along which simulated line averaged-electron density
is evaluated for comparison with experimental data. }
\end{figure}

\begin{figure}[H]
\begin{minipage}[t]{0.4\columnwidth}%
\subfloat[]{\raggedright{}\includegraphics[scale=0.4]{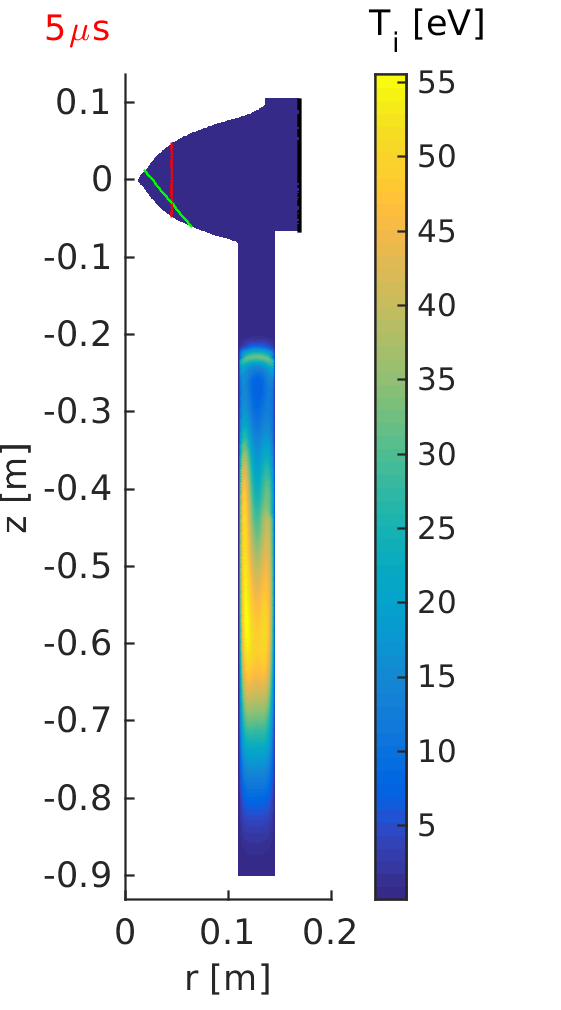}}\hfill{}\subfloat[]{\raggedright{}\includegraphics[scale=0.4]{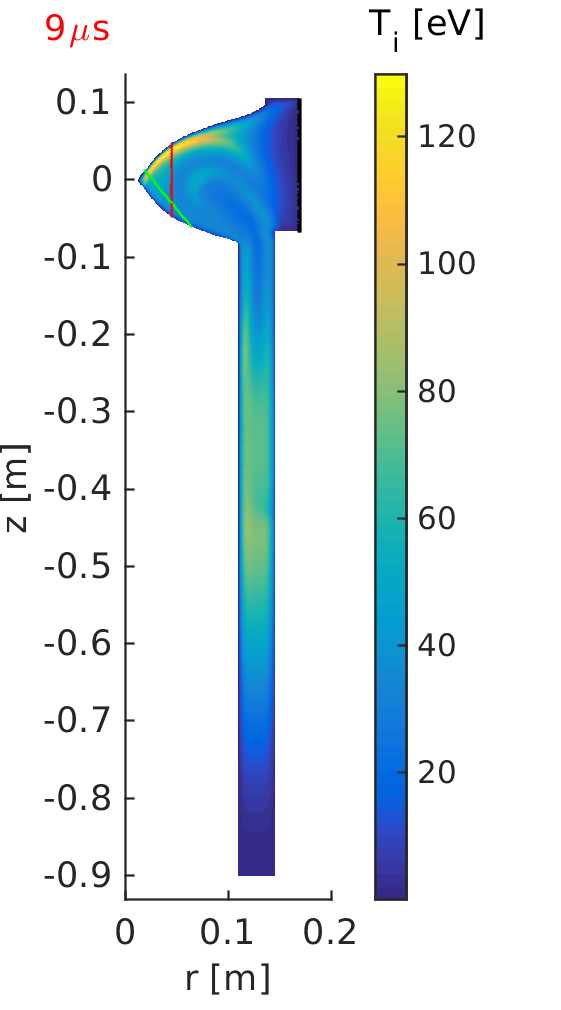}}%
\end{minipage}\hfill{}%
\begin{minipage}[t]{0.2\columnwidth}%
\subfloat[]{\raggedright{}\includegraphics[scale=0.2]{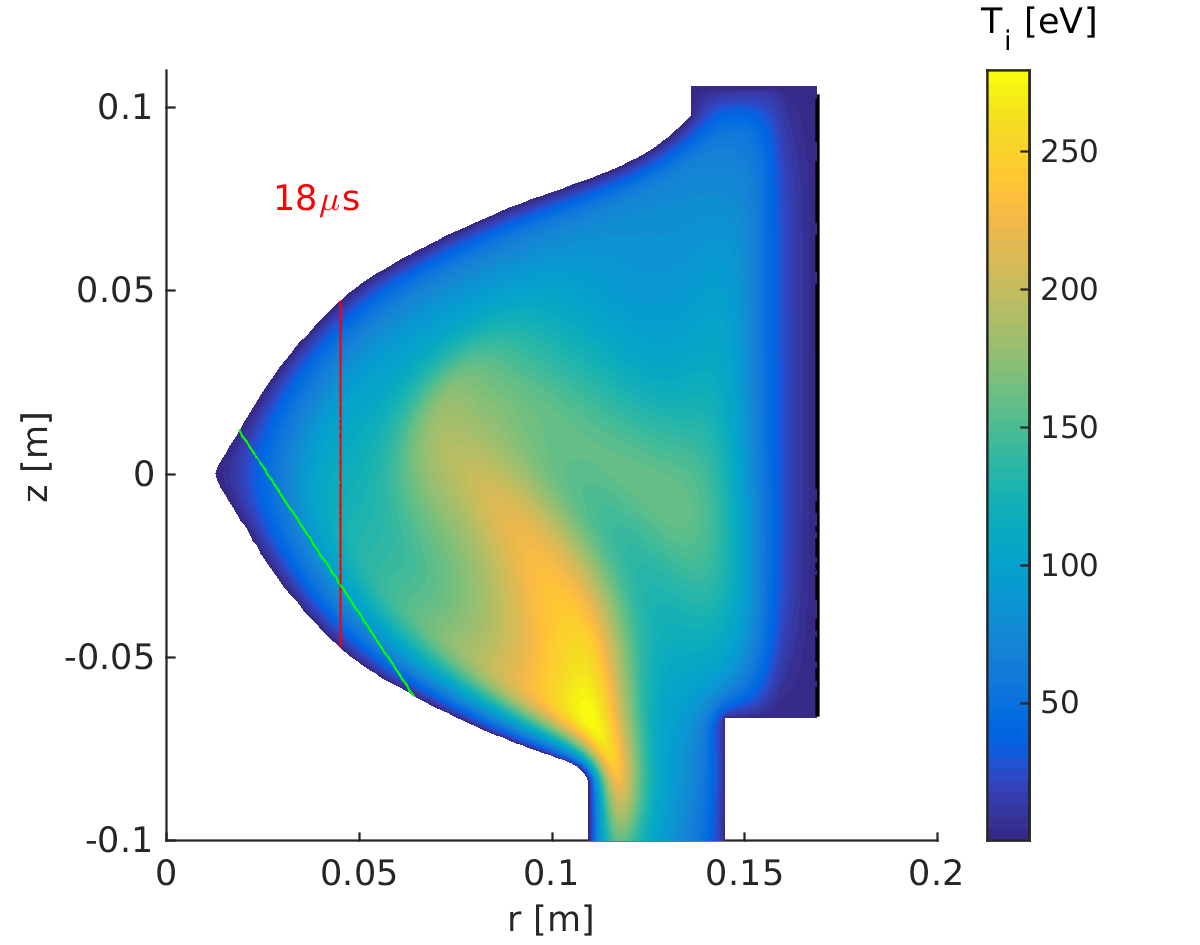}}\hfill{}\subfloat[]{\raggedright{}\includegraphics[scale=0.2]{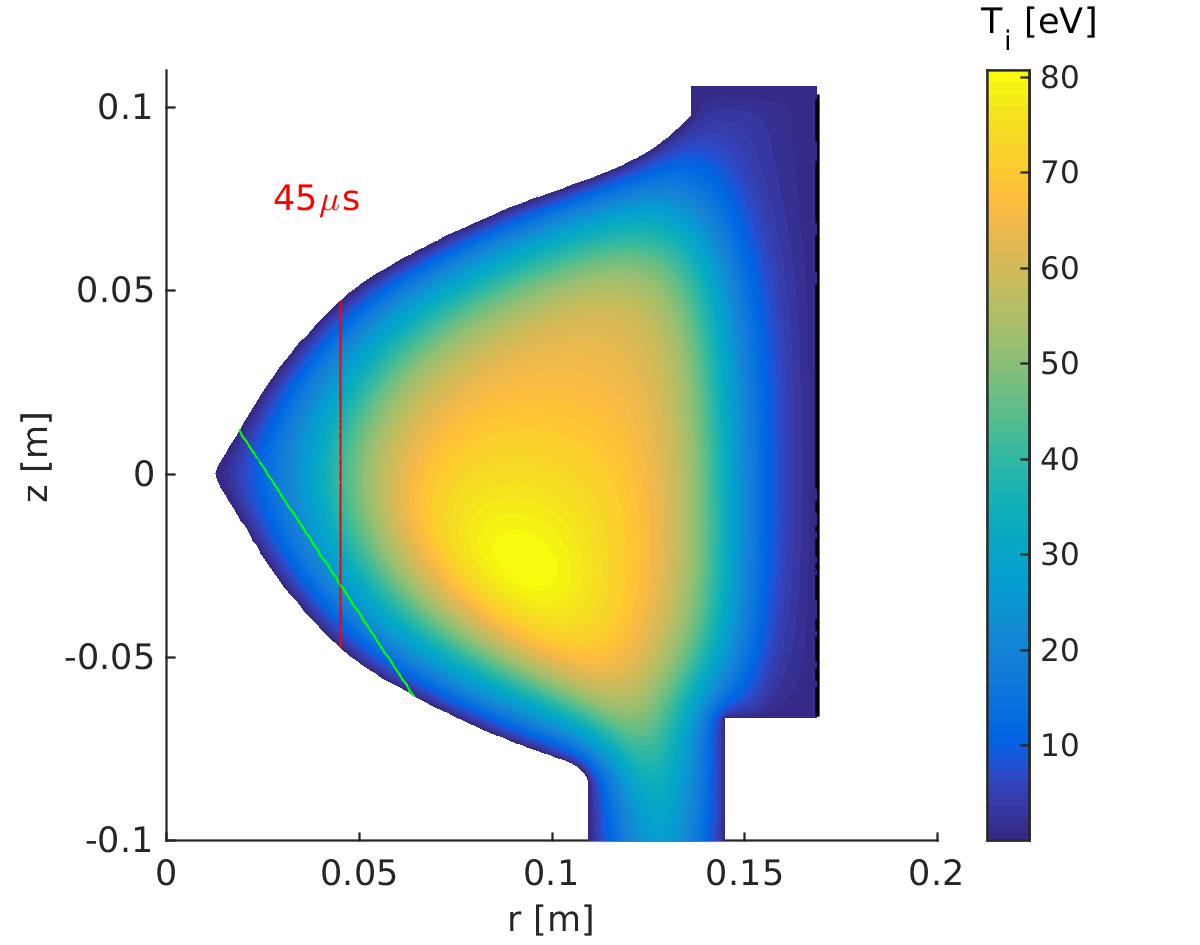}}%
\end{minipage}\hfill{}%
\begin{minipage}[t]{0.2\columnwidth}%
\subfloat[]{\raggedright{}\includegraphics[scale=0.2]{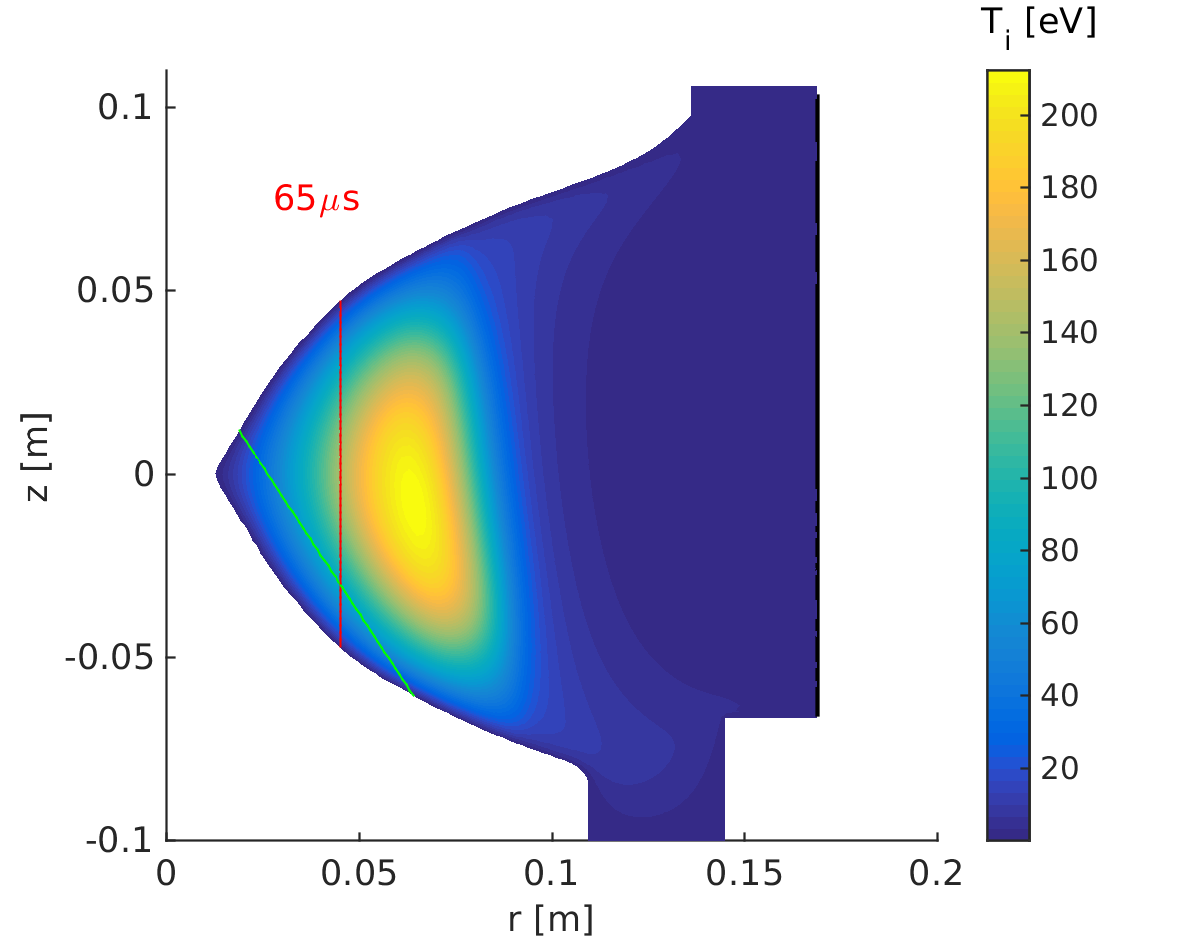}}\hfill{}\subfloat[]{\raggedright{}\includegraphics[scale=0.2]{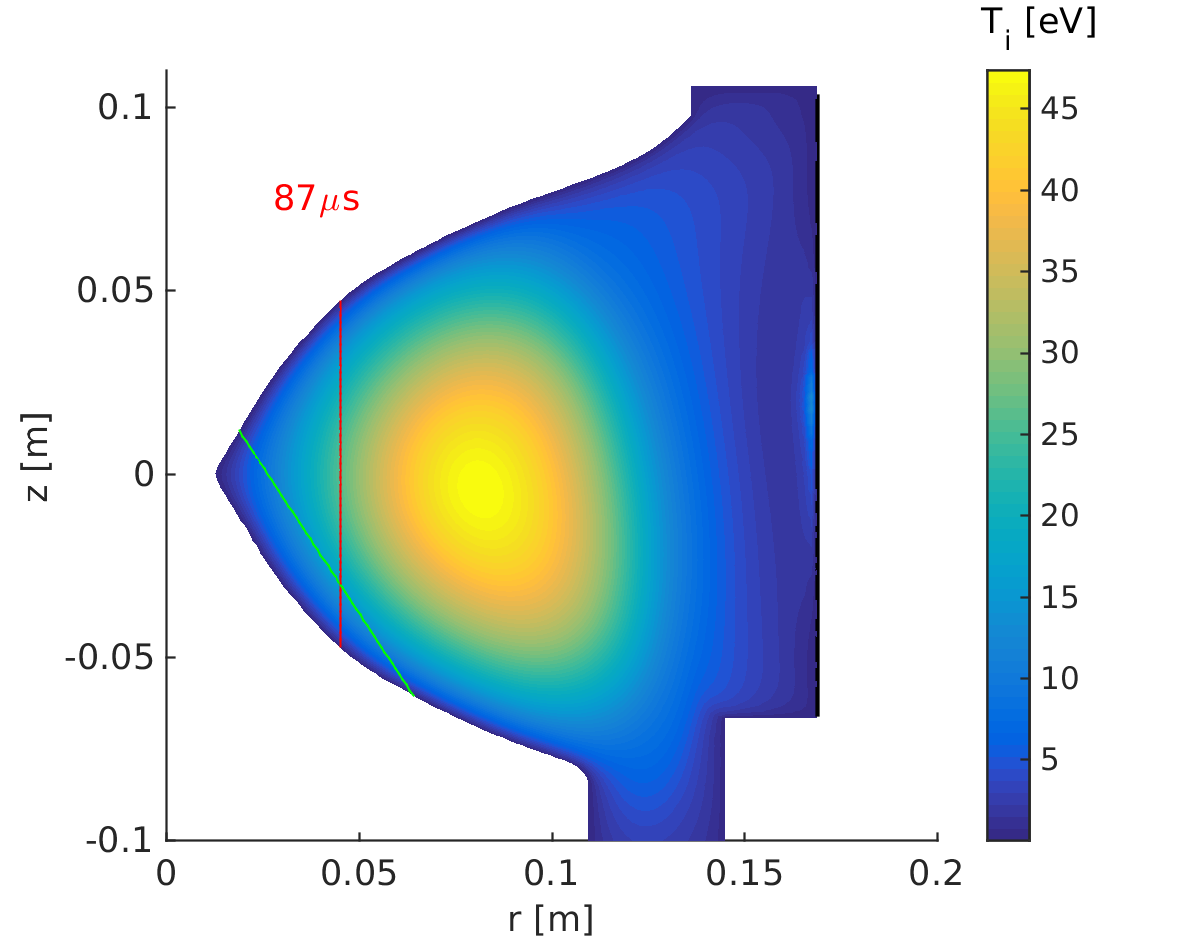}}%
\end{minipage}

\caption{\label{fig: Ti_11coils}Ion temperature profiles at $5\upmu$s (a),
$9\upmu$s (b), $18\upmu$s (c), $45\upmu$s (d), $65\upmu$s (e),
$87\upmu$s (f). The vertical red chord and diagonal green chord in
the CT containment region represent the lines of sight of the ion-Doppler
measurements, along which simulated line-averaged ion temperature
is evaluated for comparison with experimental data. }
\end{figure}

In total, there are over one hundred code input parameters. Principal
code input parameters for the simulation results presented in this
section are as follows:
\begin{itemize}
\item $h_{e}$=2mm defines the mesh resolution. 
\item $V_{form}=16$kV, $I_{main}=70$A, $V_{lev}=16$kV, $V_{comp}=18$kV:
The formation gun voltage waveform used to calculate $\Phi_{form}(t)$
in the expression for $f_{external}(z,t)=f_{form}(z,t)$ (equation
\ref{eq:145.01}) is taken from an experimental measurement for a
typical shot with $V_{form}=16$kV, where $V_{form}$ is the voltage
to which the formation capacitor banks were charged. As described
in section \ref{subsec:Boundary-conditions}, the boundary conditions
for $\psi$, pertaining to $\psi_{main},\,\psi_{lev},\,\mbox{and }\psi_{comp}$
are obtained using FEMM models. For the FEMM model used to find $\psi_{main}$,
the dc current in the main coil was set to $I_{main}=$70A. Similarly,
the currents in the levitation and compression coils in the FEMM models,
used to obtain boundary conditions for $\psi_{lev}$ and $\psi_{comp}$,
were set to the experimentally measured values corresponding to experimentally
recorded $V_{lev}=16$kV and $V_{comp}=18$kV (voltages to which the
levitation and compression capacitor banks were charged) respectively.
Note that the waveform for $V_{form}$, and the boundary values pertaining
to the peak values of $\psi_{main},\,\psi_{lev},\,\mbox{and }\psi_{comp}$,
can be linearly scaled according to the settings for $V_{form},\,I_{main,}\,V_{lev}$
and $V_{comp}$ that are required for a particular simulation.
\item $t_{comp}=45\upmu$s means that magnetic compression ($i.e.,$ superimposition
of the $\psi_{comp}$ boundary conditions, scaled by experimentally
recorded $\widetilde{I}_{comp}(t)$, on the $\psi_{main}$ and $\psi_{lev}$
boundary conditions, see section \ref{subsec:Boundary-conditions})
is started $45\upmu$s into the simulation.
\item $m_{0}=4$ implies that the ion mass is four proton masses ($m_{p}=1.67\times10^{-27}$kg),
which is relevant for modeling helium plasmas.
\item $Z=1.3$ defines the estimate for the (volume averaged) ion charge,
and determines the ratio of electron to ion number density ($n_{e}=Z\,n_{i})$.
The plasma is being modeled as a single fluid, while the energy equation
is split into ion and electron components, so $Z$ determines the
ratio between the ion and electron pressures ($p=p_{i}+p_{e}=n\,(T_{i}[\mbox{J}]+Z\,T_{e}[\mbox{J}])$,
where $n\iff n_{i}$ ). $Z$ also enters into the determination of
plasma Spitzer resistivity and the ion-electron heat exchange term
$Q_{ie}$. Fully ionized pure helium plasma would have $Z=2$, here
we assume that a proportion of the helium atoms are only singly ionized.
$Z$ can be increased to include some of the effects of the inclusion
of high ionic-charge impurities in the plasma.
\item Code inputs $n_{0}${[}m$^{-3}${]}$=9\times10^{20}${[}m$^{2}$/s{]}
and $\sigma_{n}^{2}=0.005$ {[}m$^{2}${]} determine the initial plasma
density distribution. $\sigma_{n}^{2}$ {[}m$^{2}${]} defines the
variance of the Gaussian function that defines the distribution, which
is centered around $z_{gp}=-0.43$m, the $z$ coordinate of the physical
location of the machine gas puff valves. The initial number density
distribution is 
\begin{equation}
n_{t0}(z)=n_{0}\,\left((n_{high}-n_{low})\,\widetilde{g}(z)+n_{low}\right)\label{eq:900}
\end{equation}
where $\widetilde{g}(z)=g(z)/\mbox{max}(g(z))$. Here, $g(z)=\frac{1}{\sqrt{2\pi\,\sigma_{n}^{2}}}\,exp\left(\frac{-\left(z-z_{gp}\right)^{2}}{2\sigma_{n}^{2}}\right)$.
$n_{high}$ and $n_{low}$ are typically set to 10 and 0.1 respectively.
Note that $n_{low}$ must be finite - density is not allowed to approach
too close to zero anywhere in the computational domain. The sound
speed of a neutral gas is given by $V_{s}=\sqrt{\frac{\gamma p}{\rho}}$.
At room temperature, for helium gas, this is $V_{s}\approx900$m/s.
In a typical shot, gas is puffed into the vacuum vessel through the
eight gas valves spaced toroidally on the outer electrode, around
$400\upmu$s before firing the formation capacitors. A simple estimate
of $d_{n}$, the maximum spread of the neutral gas around the valves
can be estimated using this time and the sound speed, as $d_{n}=0.36$m,
so that choosing $\sigma_{n}^{2}=0.005$m$^{2}$, corresponding to
a Gaussian base width of around 0.45m, is a reasonable estimate for
the initial spread. See figure \ref{fig: n_11coils}(a) for the initial
density profile pertaining to these input parameters. Note that for
the simulation presented here, interaction between plasma and neutral
fluid is included in the model, where the initial neutral particle
inventory distribution is also determined by a Gaussian profile centered
around the location of the machine gas valves. Results from simulated
neutral fluid dynamics, and a description of the model for plasma-neutral
interaction, won't be discussed further here; this material is presented
separately in \cite{Neut_paper}.
\item $\zeta$ is the density diffusion coefficient, note that 50 {[}m$^{2}$/s{]}
is close to the minimum value required for numerical stability for
simulations that include CT formation and magnetic compression.
\item $\nu_{num}=700$m$^{2}$/s and $\nu_{phys}=410$m$^{2}$/s: $\nu_{num}$
is the value for the viscosity diffusion coefficient that defines
$\mu$, in the viscous terms in the momentum equation (see equation
\ref{eq:517.4}), as $\mu_{num}=m_{i}n_{0}\,\nu_{num}$. To reduce
viscous heating of the ions to acceptable levels ($i.e.,$ to achieve
a match between experimentally measured and simulated ion temperatures),
a reduced value $\mu_{phys}=m_{i}n_{0}\,\nu_{phys}$ defines $\mu$
in the viscous term in the ion energy equation (see equation \ref{eq:517.5}).
At 700m$^{2}$/s, $\nu_{num}$ is close to the minimum value required
for numerical stability ($i.e,$ sufficient velocity field smoothing),
while $\nu_{phys}=410$m$^{2}$/s is chosen as a crude estimate for
the \textquotedbl physical\textquotedbl{} viscosity because it leads
to simulated ion temperature close to the levels indicated by the
ion-Doppler system. Ideally, the strategy of having different values
for $\nu_{num}$ and $\nu_{phys}$ would be avoided because it breaks
the conservation of system energy pertaining to the viscous terms
(the terms in the last set of square brackets in equation \ref{eq:518}).
However, it was found that the discrepancy does not, apart from a
reduction in ion viscous heating, significantly alter these particular
simulation results, and appears at this stage to be an acceptable
method for reducing ion heating in simulation scenarios with low density
and high plasma velocities.
\item Code input parameters pertaining to the thermal diffusion coefficients,
which are held constant in this simulation, are $\chi_{\parallel e}\sim16000\,[\mbox{m}^{2}\mbox{/s}],\,\chi_{\parallel i}\sim5000\,[\mbox{m}^{2}\mbox{/s}],\,\chi_{\perp e}\sim240\,[\mbox{m}^{2}\mbox{/s}]$,
and $\chi_{\perp i}\sim120$$[\mbox{m}^{2}\mbox{/s}]$. The high field-parallel
thermal diffusion coefficients represent rapid equilibration of temperature
along field lines, while the field-perpendicular coefficients are
chosen to match the experimentally observed resistive decay rate of
CT currents. Enhanced thermal diffusion acts as a proxy for anomalous
energy sink mechanisms, including radiative losses.
\end{itemize}
\begin{figure}[H]
\subfloat[]{\raggedright{}\includegraphics[scale=0.4]{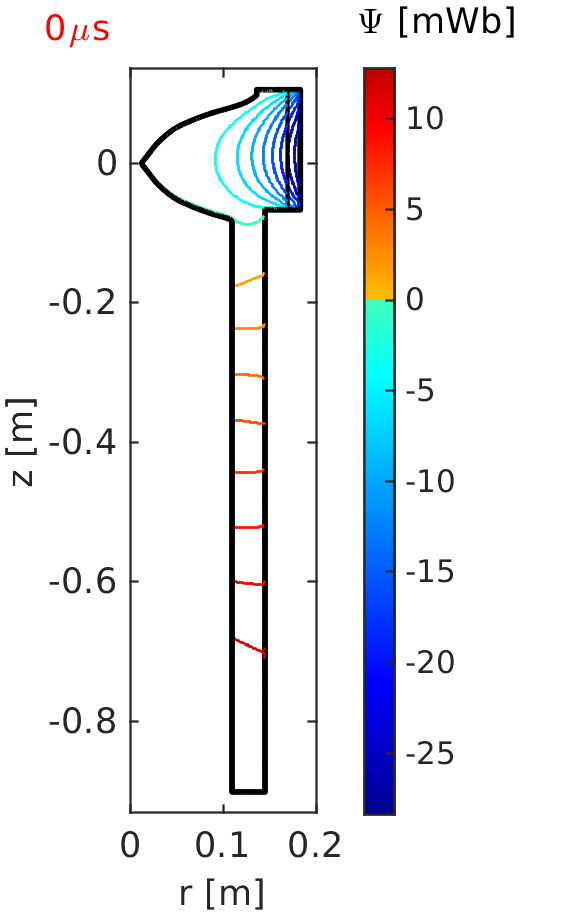}}\hfill{}\subfloat[]{\raggedright{}\includegraphics[scale=0.4]{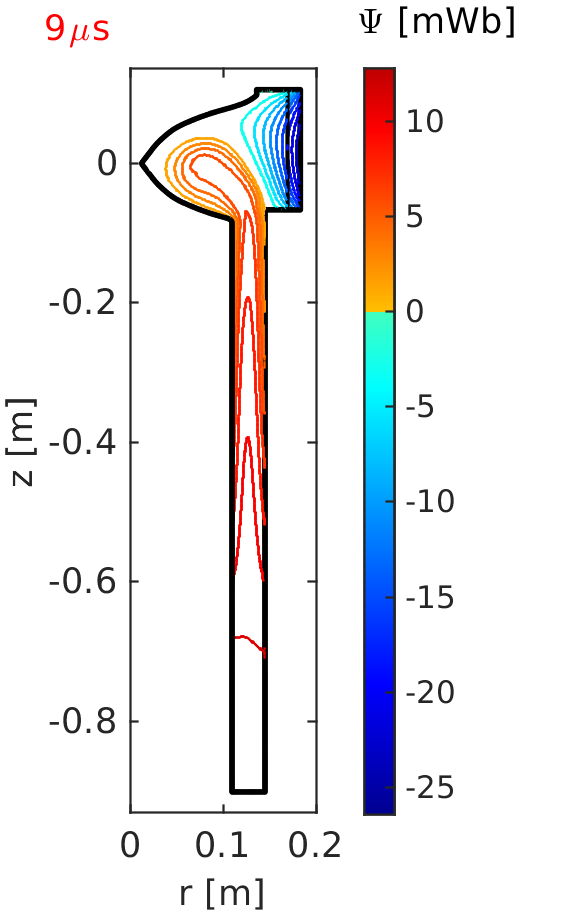}}\hfill{}\subfloat[]{\raggedright{}\includegraphics[scale=0.4]{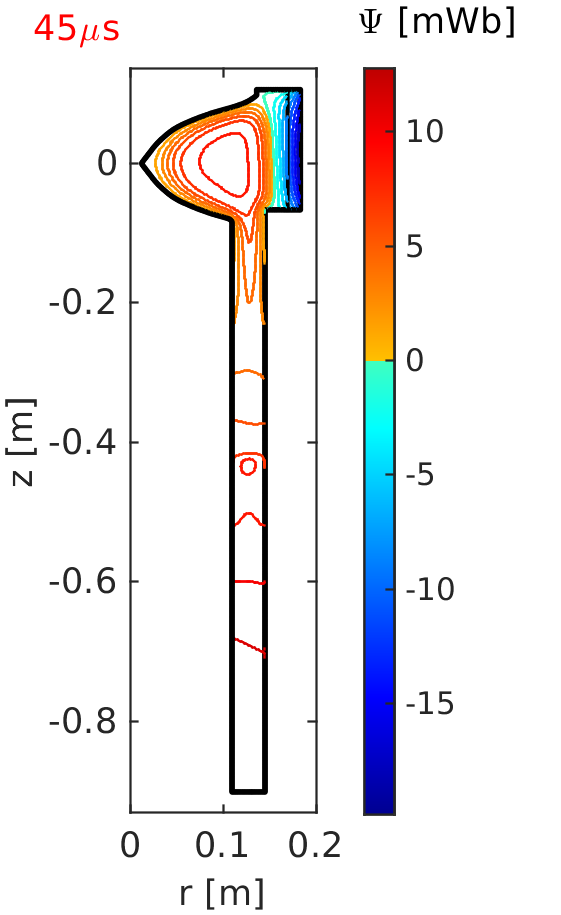}}\hfill{}\subfloat[]{\raggedright{}\includegraphics[scale=0.4]{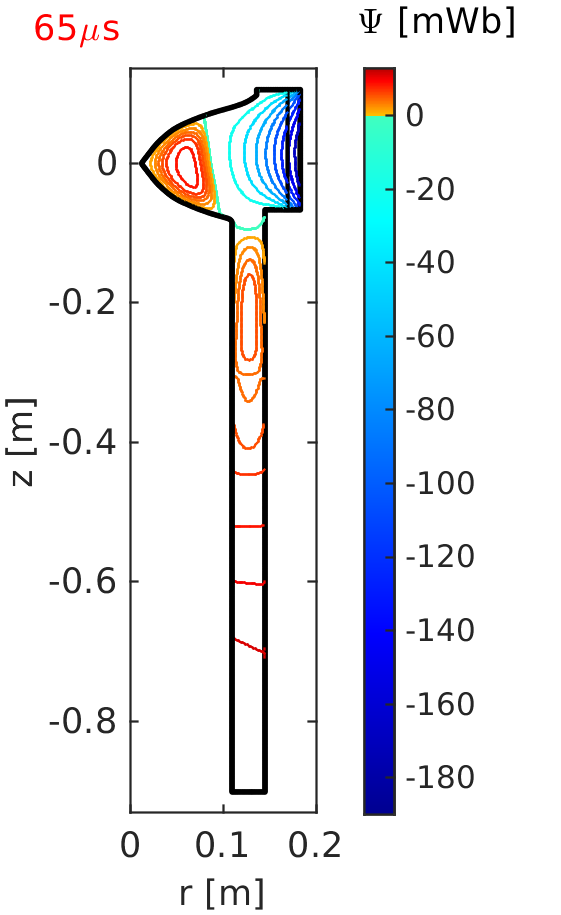}}

\caption{\label{fig: psi_11coils}Poloidal flux contours, 11-coil configuration,
at $0\upmu$s (a), $9\upmu$s (b), $45\upmu$s (c), $65\upmu$s (d).
Note that the colorbar scaling changes over time (here and in subsequent
contour plots). }
\end{figure}
Figure \ref{fig: psi_11coils} shows simulated $\psi$ contours, representing
poloidal magnetic field lines, at various times for a simulation of
CT formation, levitation and magnetic compression in the eleven coil
configuration, with code input parameters as defined above. Note that
the vertical black line at the top-right of the figures at $r\sim17\mbox{cm}$
represents the inner radius of the insulating wall. Vacuum field only
is solved for to the right of the line, and the plasma dynamics are
solved for in the remaining solution domain to the left of the line.
The inner radius of the stack of eleven levitation/compression coils
(which are not depicted here) is located at the outer edge of the
solution domain, at $r\sim18\mbox{cm}$. Simulation times are notated
in red at the top left of the figures. Note that max$(\psi)$ decreases
slowly over time as the CT decays, while min$(\psi)$ increases as
the levitation current in the external coils decays, and then drops
off rapidly as the compression current in the external coils is increased,
starting at $t_{comp}=45\upmu\mbox{s}$ in this simulation. At time
$t=0,$ the stuffing field ($\psi>0$) due to currents in the main
coil fills the vacuum below the containment region, and has soaked
well into all materials around the gun, while the levitation field
fills the containment region. As described in section \ref{subsec:Formation-Simulations},
simulated CT formation is initiated with the addition of toroidal
flux below the gas puff valves located at $z=-0.43$m; initial intra-electrode
radial formation current is assumed to flow at the z-coordinate of
the valves. Open field lines that are resistively pinned to the electrodes,
and partially frozen into the conducting plasma, have been advected
by the $\mathbf{J}_{r}\times\mathbf{B}_{\phi}$ force into the containment
region by $t=9\upmu\mbox{s}$ ($\mathbf{J}_{r}$ is the radial formation
current density across the plasma between the electrodes, and $\mathbf{B}_{\phi}$
is the toroidal field due to the axial formation current in the electrodes).
By $45\upmu\mbox{s}$, open field lines have reconnected at the entrance
to the containment region to form closed CT flux surfaces. At these
early times, open field lines remain in place surrounding the CT.
Compression starts at $45\upmu\mbox{s}$ and peak compression is at
$65\upmu\mbox{s}$. Note magnetic compression causes closed CT poloidal
field lines that extend down the gun to be pinched off at the gun
entrance, where they reconnect to form a second smaller CT. Field
lines that remain open surrounding the main CT are then also reconnectively
pinched off, forming additional closed field lines around the main
CT, while the newly reconnected open field lines below the main CT
act like a slingshot that advects the smaller CT down the gun, as
can be seen at $65\upmu\mbox{s}$. The CT expands again between after
peak compression at $65\upmu\mbox{s}$ as the compression current
in the external levitation/compression coils decreases. Note that
2D simulations, which neglect inherently three dimensional turbulent
transport and flux conversion, generally overestimate the level of
hollowness of the field profiles.\\
\\
\begin{figure}[H]
\subfloat[]{\raggedright{}\includegraphics[scale=0.4]{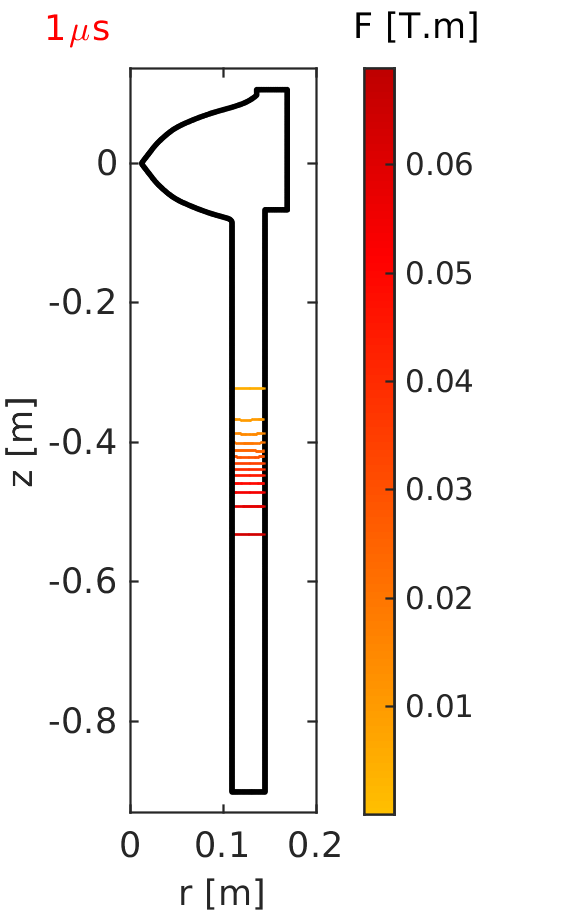}}\hfill{}\subfloat[]{\raggedright{}\includegraphics[scale=0.4]{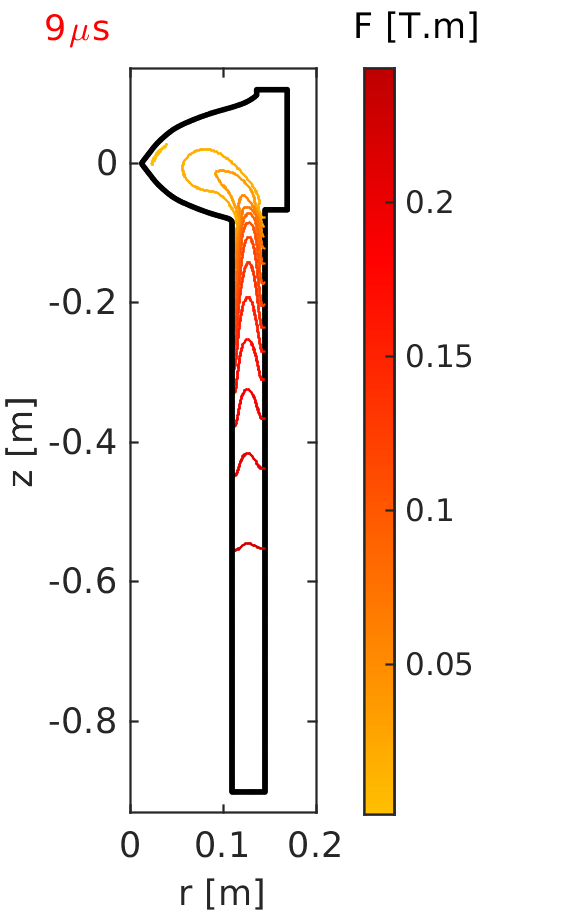}}\hfill{}\subfloat[]{\raggedright{}\includegraphics[scale=0.4]{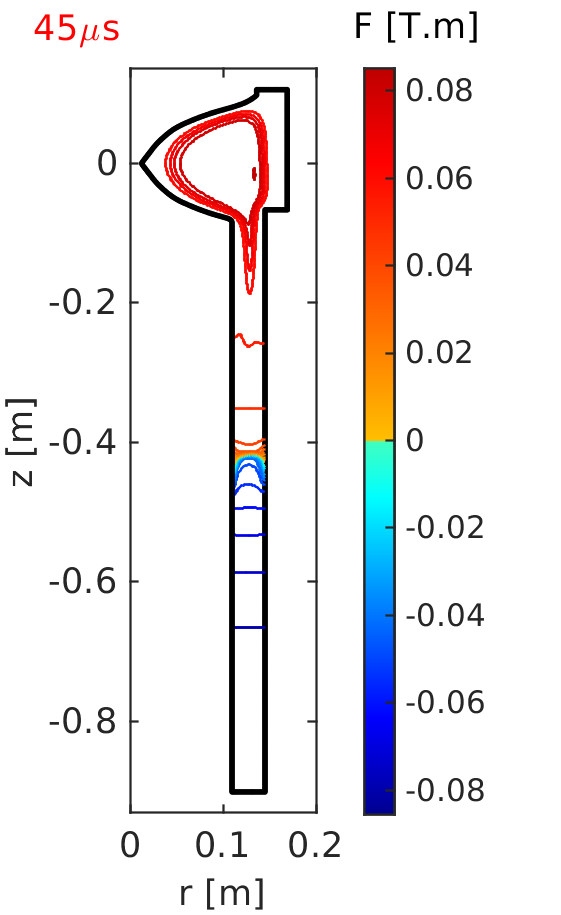}}\hfill{}\subfloat[]{\raggedright{}\includegraphics[scale=0.4]{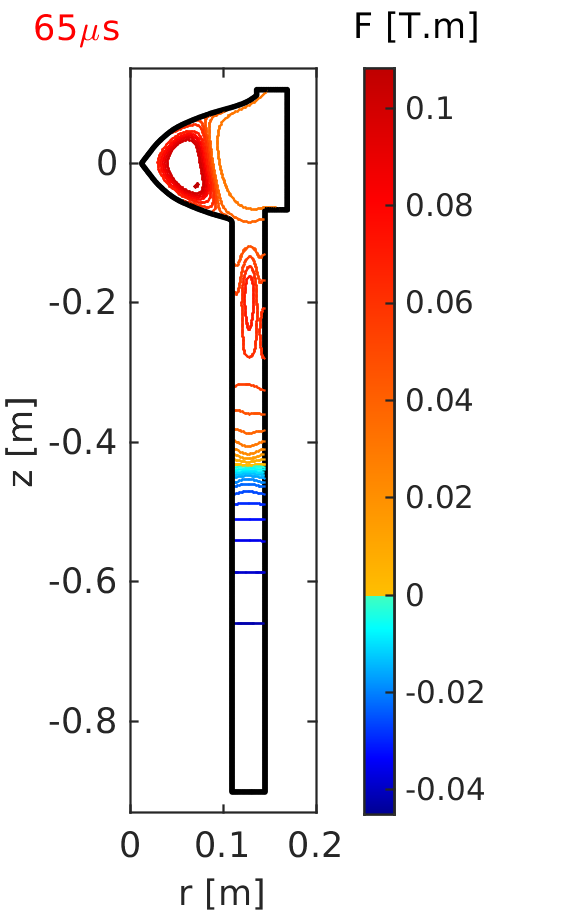}}

\caption{\label{fig: F_11coils} $f$ contours (lines of poloidal current),
at $0\upmu$s (a), $9\upmu$s (b), $45\upmu$s (c), $65\upmu$s (d)}
\end{figure}
Contours of $f$ at various times are shown in figure \ref{fig: F_11coils}.
Recall that contours of $f=rB_{\phi}$ represent lines of poloidal
current. Initially, $f$ is zero at all nodes. At $1\upmu$s, it can
seen how radial formation current between the gun electrodes, corresponding
to the axial gradient of $f$, is concentrated around the gas valve
location, at $z=-$0.43m. Closely spaced contours indicate regions
of high gradients and high poloidal currents. To simulate formation,
toroidal flux is being added to the domain according to the geometric
profile shown in figure \ref{fig:gform_Phiform_Vform}(a). At $1\upmu$s,
$f$ is constant, at its highest value in the domain, below the lowermost
contour at $z\sim-0.55$m representing $f\sim0.07$ T-m. By $9\upmu$s,
plasma has been advected upwards, and poloidal current is flowing
along open poloidal magnetic field lines that remain resistively pinned
to the electrodes down the gun. By $45\upmu$s (figure \ref{fig: F_11coils}(c))
poloidal current is flowing along closed and open poloidal magnetic
field lines, and most of the toroidal flux in the domain has been
advected with the plasma into the containment region. At this time,
as the level of toroidal flux added to the domain below the gas valve
locations, which is equal to the total toroidal flux in the domain,
falls off (see figure \ref{fig:phi_cons}), radial currents between
the electrodes further down the gun reverse direction, and $f$ becomes
negative, as consequences of toroidal flux conservation and upward
advection of toroidal flux. As shown in figure \ref{fig: F_11coils}(d),
poloidal current flows around the closed poloidal field lines of the
main CT as it is being compressed, and around the closed poloidal
field lines of the smaller pinched off CT at 65$\upmu$s. Total toroidal
flux is conserved, leading to induced internal currents between various
points on the conducting walls. In particular, the wall-to-wall currents,
that are apparent outboard of the CT at peak compression (figure \ref{fig: F_11coils}
(d)), are thought to be a representation of real physical phenomenon
which led to experimentally measured toroidal field signals that suggested
the appearance of an external kink instability at compression, as
described in detail in \cite{thesis,exppaper}.
The initial electron density profile is shown in figure \ref{fig: n_11coils}(a).
As described earlier, the initial density profile is defined with
a Gaussian profile centered around the axial coordinate of the locations
of the gas puff valves. Plasma is advected up the gun during CT formation,
and bubbles-in to the containment region, as shown in figure \ref{fig: n_11coils}(b).
In figures \ref{fig: n_11coils}(c) and (d), it can be seen how density
rises over compression around the CT core, and that a region of low
density remains outboard of the compressed CT. 
\begin{figure}[H]
\begin{minipage}[t]{0.4\columnwidth}%
\subfloat[]{\raggedright{}\includegraphics[scale=0.4]{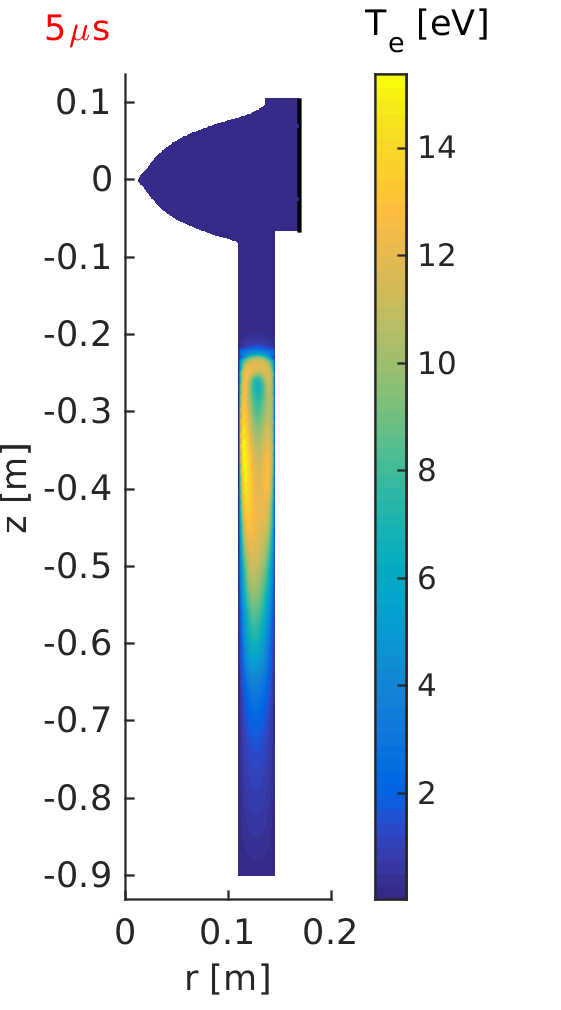}}\hfill{}\subfloat[]{\raggedright{}\includegraphics[scale=0.4]{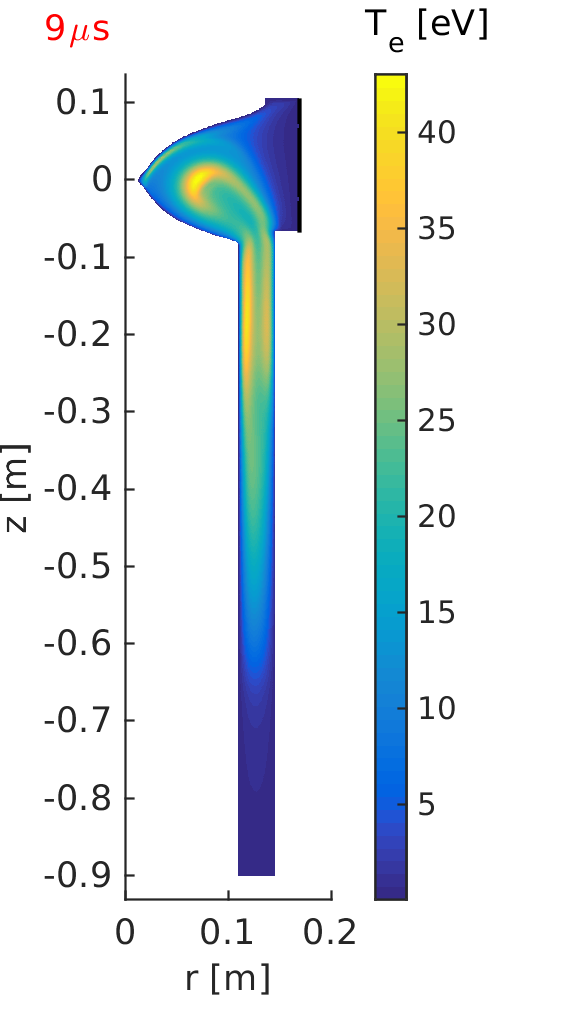}}%
\end{minipage}\hfill{}%
\begin{minipage}[t]{0.2\columnwidth}%
\subfloat[]{\raggedright{}\includegraphics[scale=0.2]{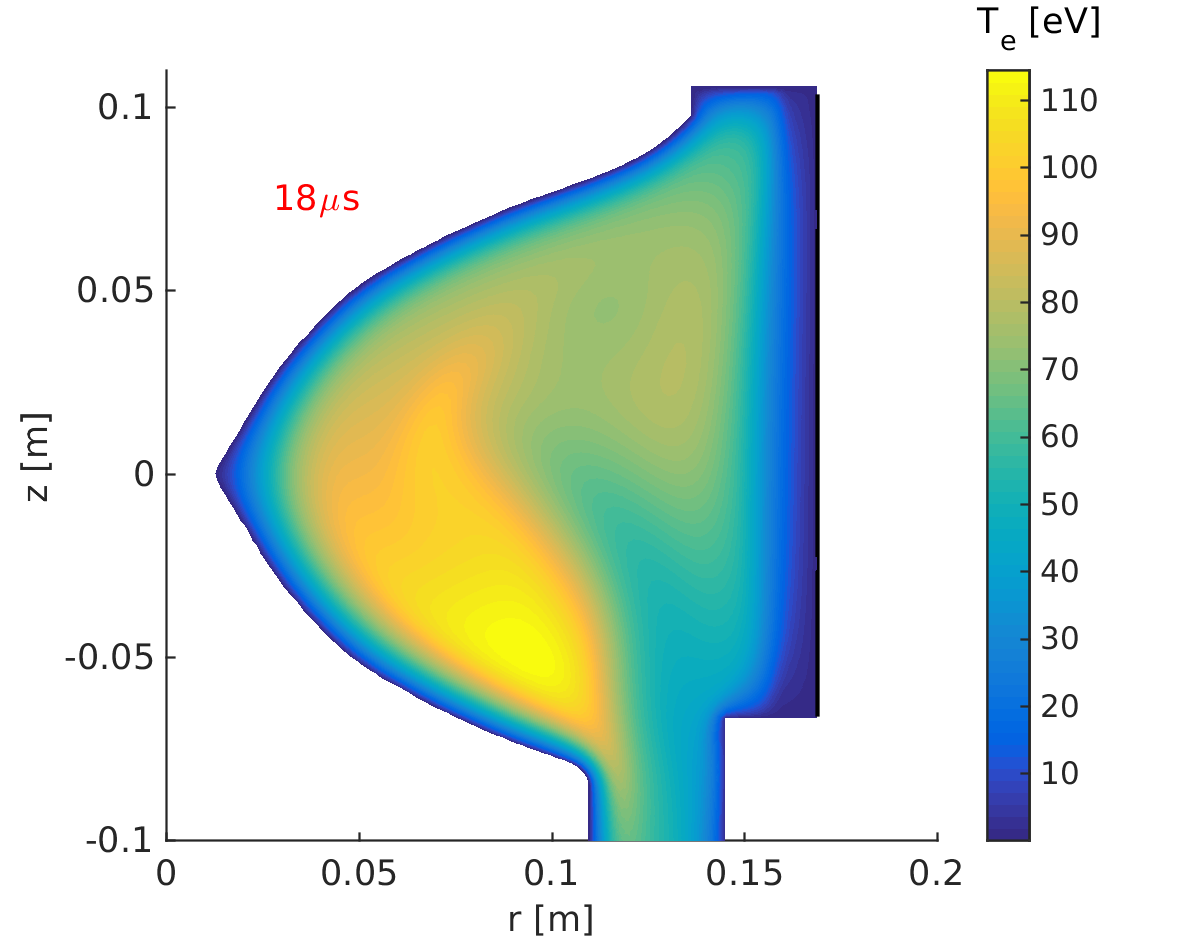}}\hfill{}\subfloat[]{\raggedright{}\includegraphics[scale=0.2]{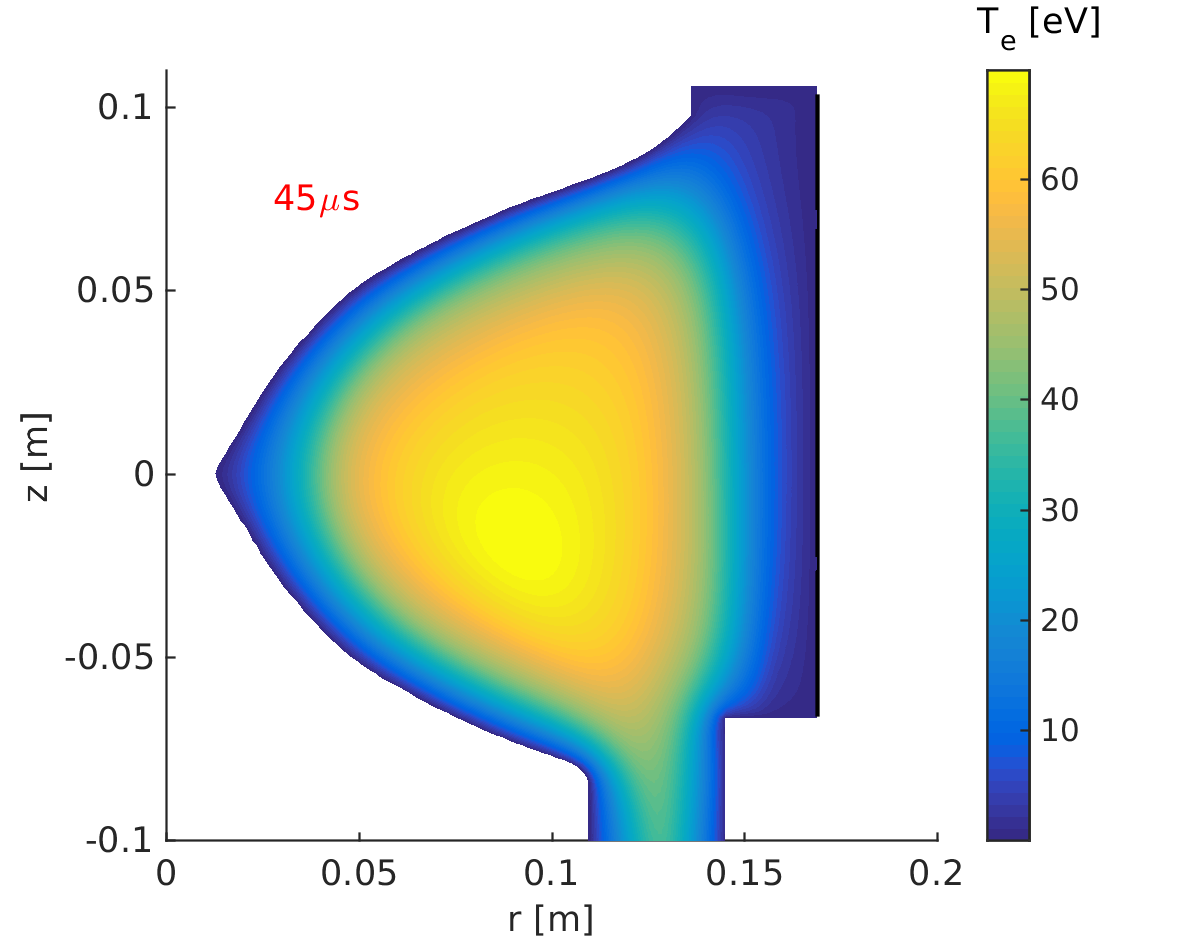}}%
\end{minipage}\hfill{}%
\begin{minipage}[t]{0.2\columnwidth}%
\subfloat[]{\raggedright{}\includegraphics[scale=0.2]{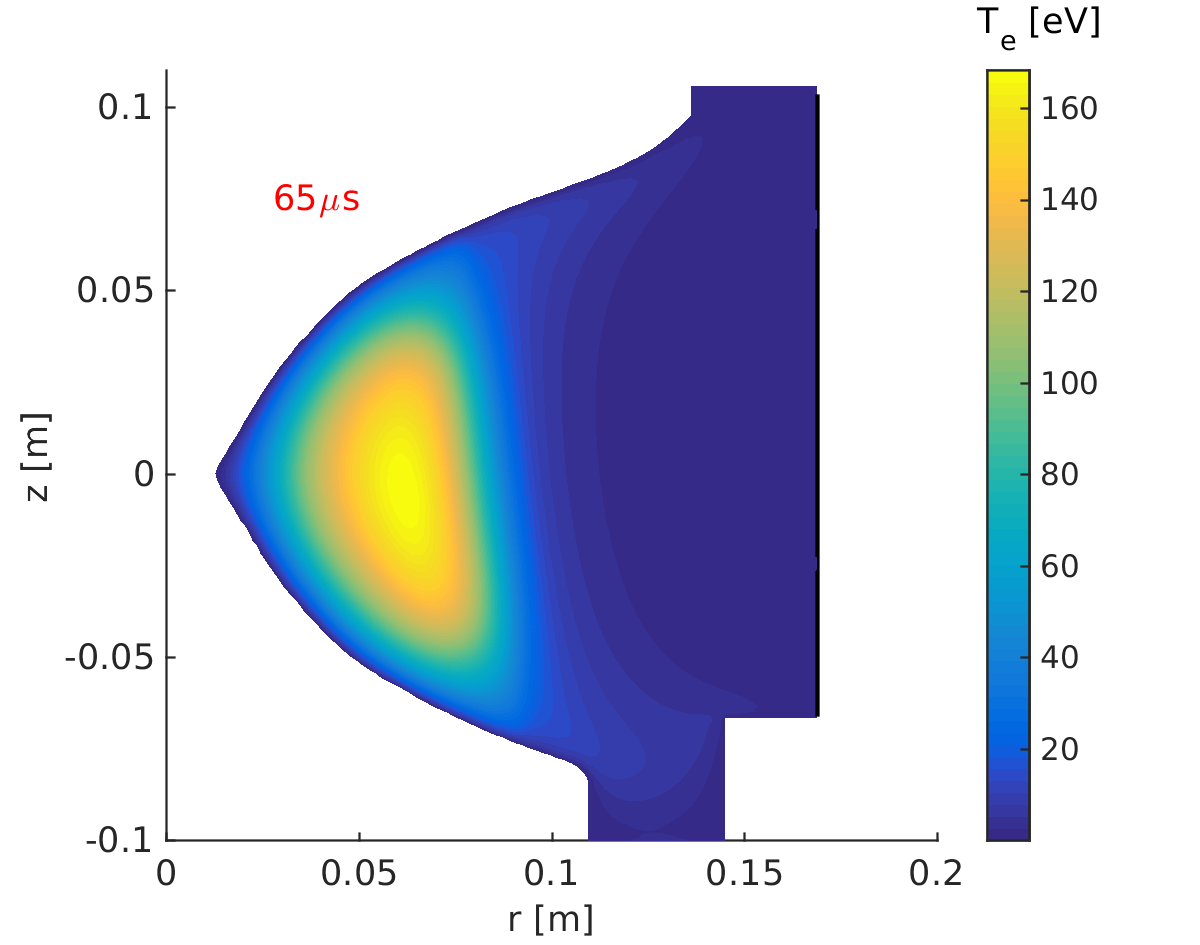}}\hfill{}\subfloat[]{\raggedright{}\includegraphics[scale=0.2]{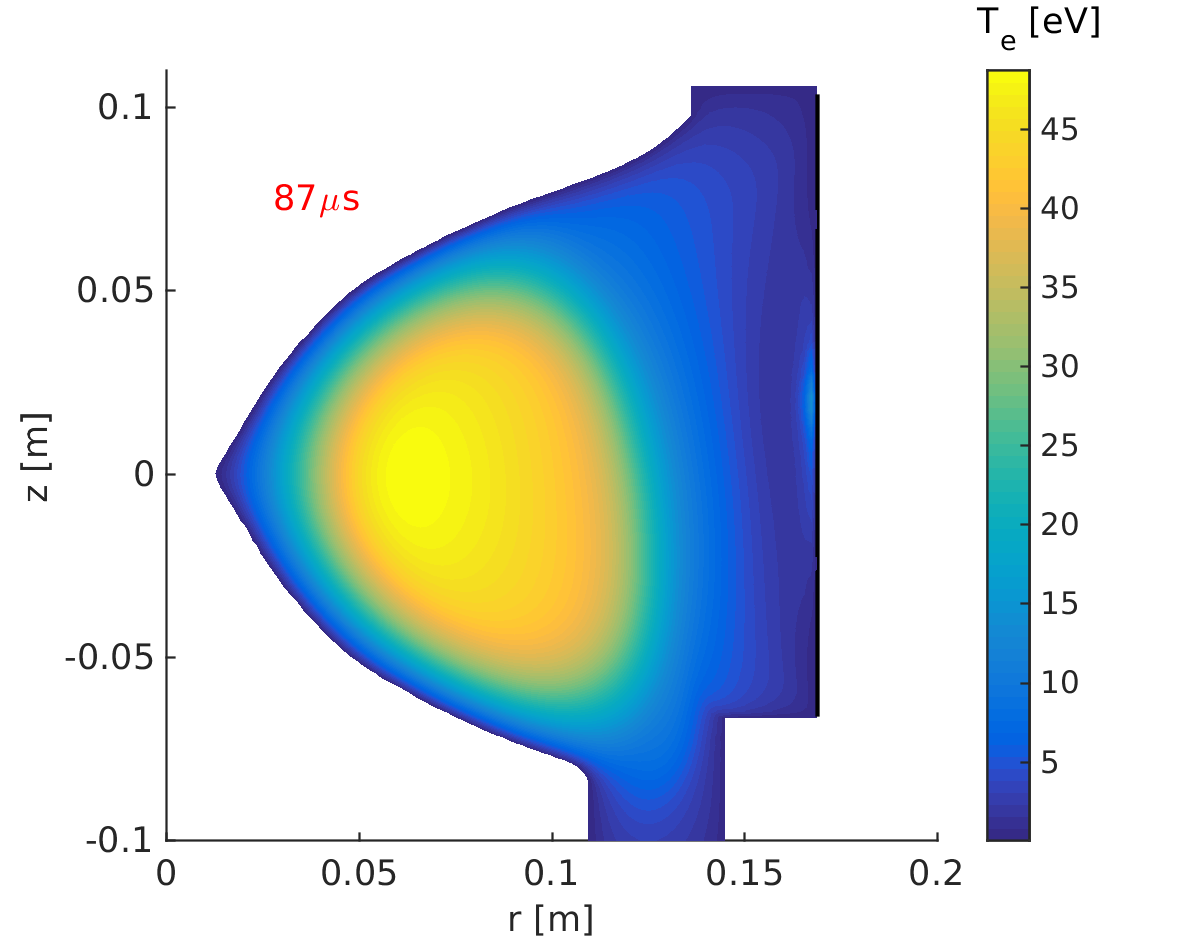}}%
\end{minipage}

\caption{\label{fig: Te_11coils}Electron temperature profiles at $5\upmu$s
(a), $9\upmu$s (b), $18\upmu$s (c), $45\upmu$s (d), $65\upmu$s
(e), $87\upmu$s (f)}
\end{figure}
Formation current ohmically heats initially cold electrons as plasma
is pushed up the gun into the containment region (figures \ref{fig: Te_11coils}(a)
and (b)). The $115$ eV attained by electrons near the containment
region entrance at $18\upmu$s (figure \ref{fig: Te_11coils}(c))
is partially due to heat exchange with ions, which, as indicated in
figure \ref{fig: Ti_11coils}(c)), have been heated to around 280eV
near the same area at this time due to viscous heating. Thermal diffusion
through the boundary causes the electron temperature around the CT
core to be reduced to around 70 eV just prior to magnetic compression
at 45$\upmu$s. Electron temperature is more than doubled, to 165
eV, at peak magnetic compression at 65$\upmu$s (figure \ref{fig: Te_11coils}(e)).
Referring to figure \ref{fig: J_11coils}, it can be seen how ohmic
heating is a principal electron heating mechanism. Compressional heating
is the main heating mechanism during magnetic compression, and is
supplemented by enhanced ohmic heating. Compression coil current falls
to zero by around 87$\upmu$s, by which time the CT has re-expanded,
and core electron temperature has fallen to around 50 eV (figure \ref{fig: Te_11coils}(f)).
Simulated ion temperatures, at the same times referred to in figure
\ref{fig: Te_11coils} for electron temperatures, are shown in figure
\ref{fig: Ti_11coils}. Over the early stages of the simulation, ion
temperature is significantly higher than electron temperature, due
to ion viscous heating during the formation process. Later, ohmic
heating of electrons is the main heating source and electron temperature
approaches ion temperature. Note that, as a result of compressional
heating in combination with heat exchange with the electrons, which
are heated ohmically by enhanced currents during compression, ion
temperature is more than doubled at peak compression, increasing from
around 80 to 210eV from 45$\upmu$s to 65$\upmu$s. The vertical red
chord and diagonal green chord in the CT containment region in figures
\ref{fig: Ti_11coils}(a) to (g) represent the lines of sight of the
ion-Doppler measurements, along which simulated line-averaged ion
temperature is evaluated for comparison with experimental data. \\
\\
\\
\begin{figure}[H]
\begin{minipage}[t]{0.4\columnwidth}%
\subfloat[]{\raggedright{}\includegraphics[scale=0.4]{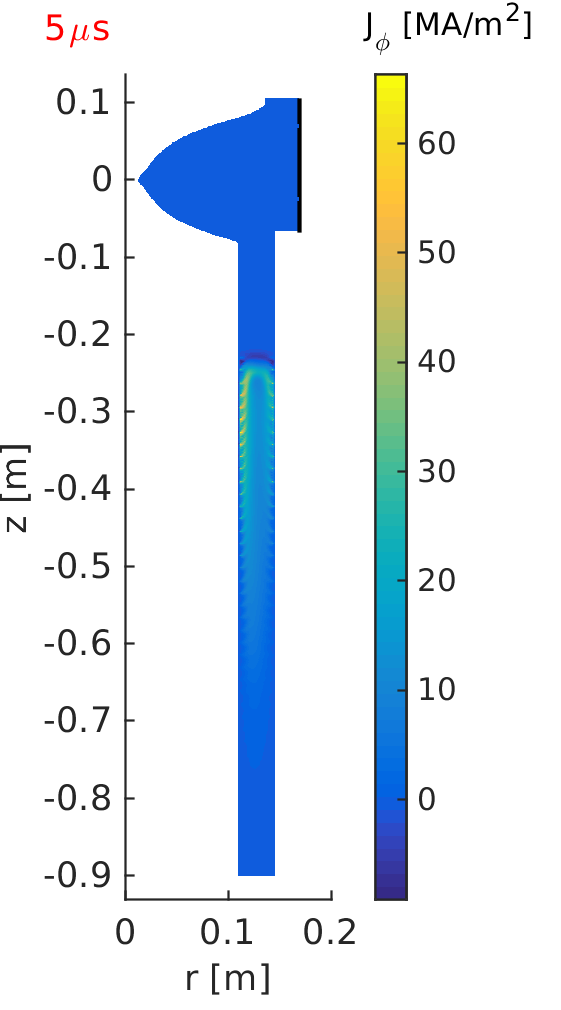}}\hfill{}\subfloat[]{\raggedright{}\includegraphics[scale=0.4]{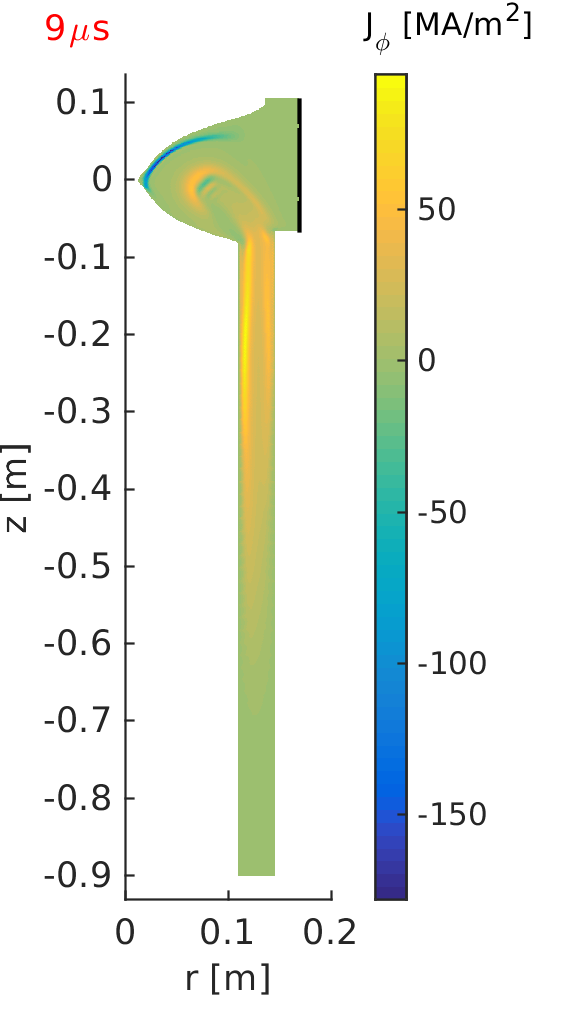}}%
\end{minipage}\hfill{}%
\begin{minipage}[t]{0.2\columnwidth}%
\subfloat[]{\raggedright{}\includegraphics[scale=0.2]{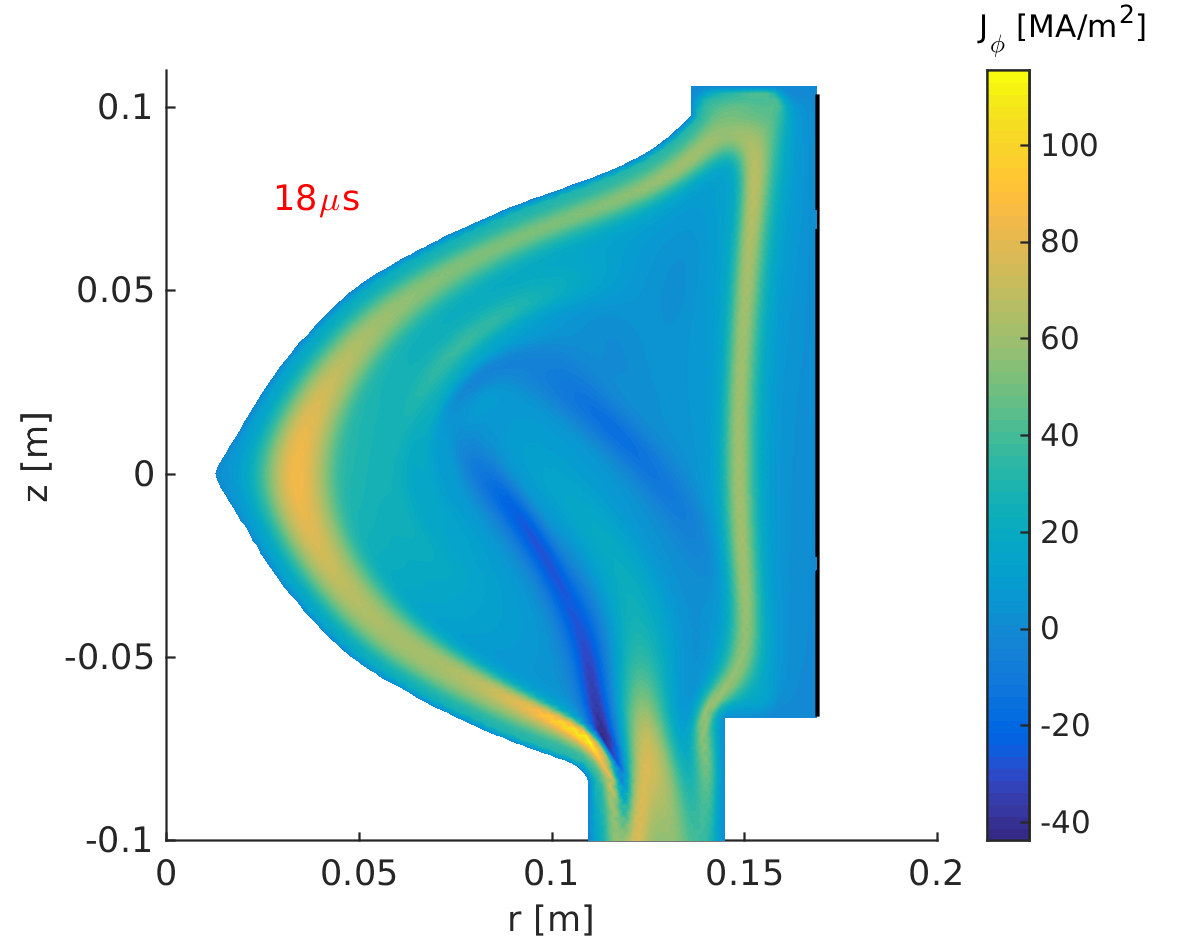}}\hfill{}\subfloat[]{\raggedright{}\includegraphics[scale=0.2]{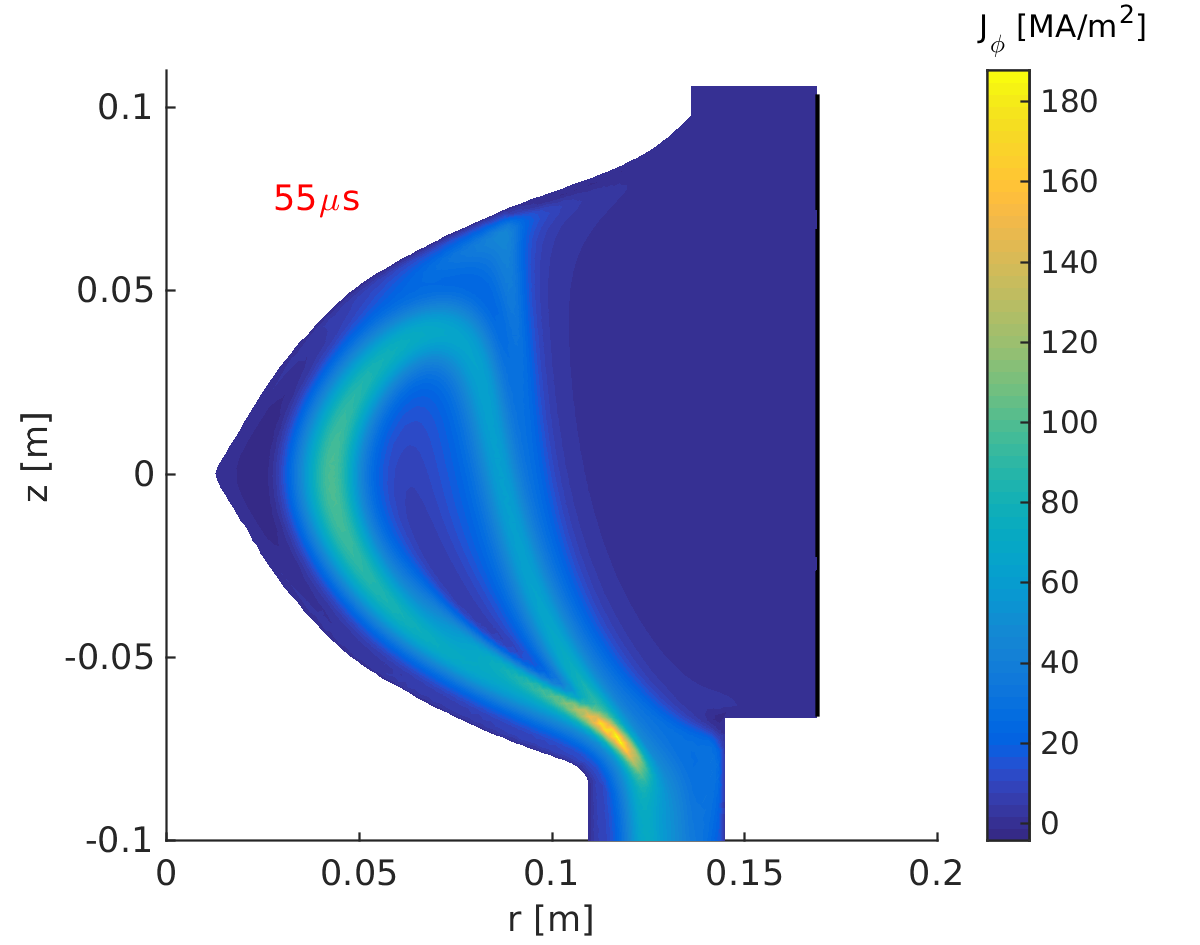}}%
\end{minipage}\hfill{}%
\begin{minipage}[t]{0.2\columnwidth}%
\subfloat[]{\raggedright{}\includegraphics[scale=0.2]{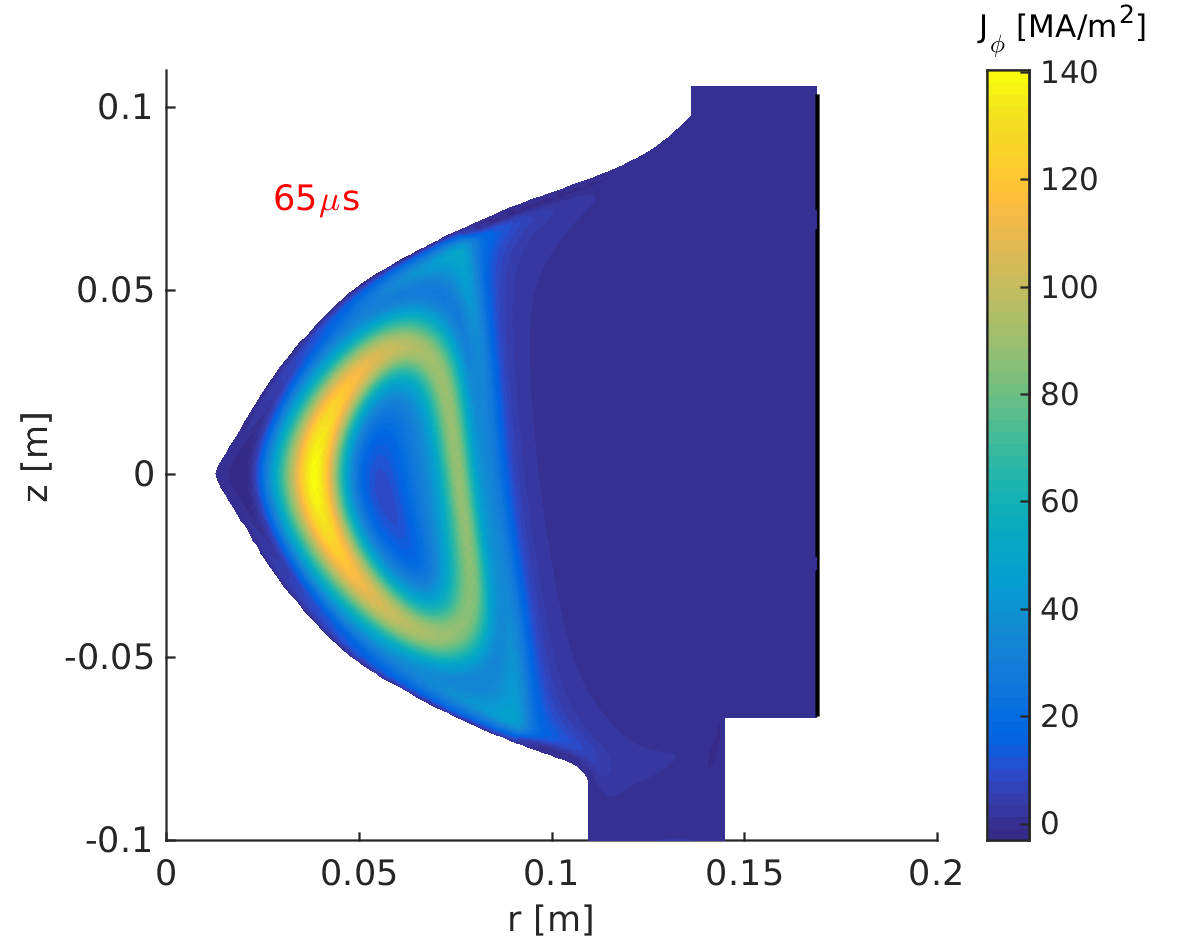}}\hfill{}\subfloat[]{\raggedright{}\includegraphics[scale=0.2]{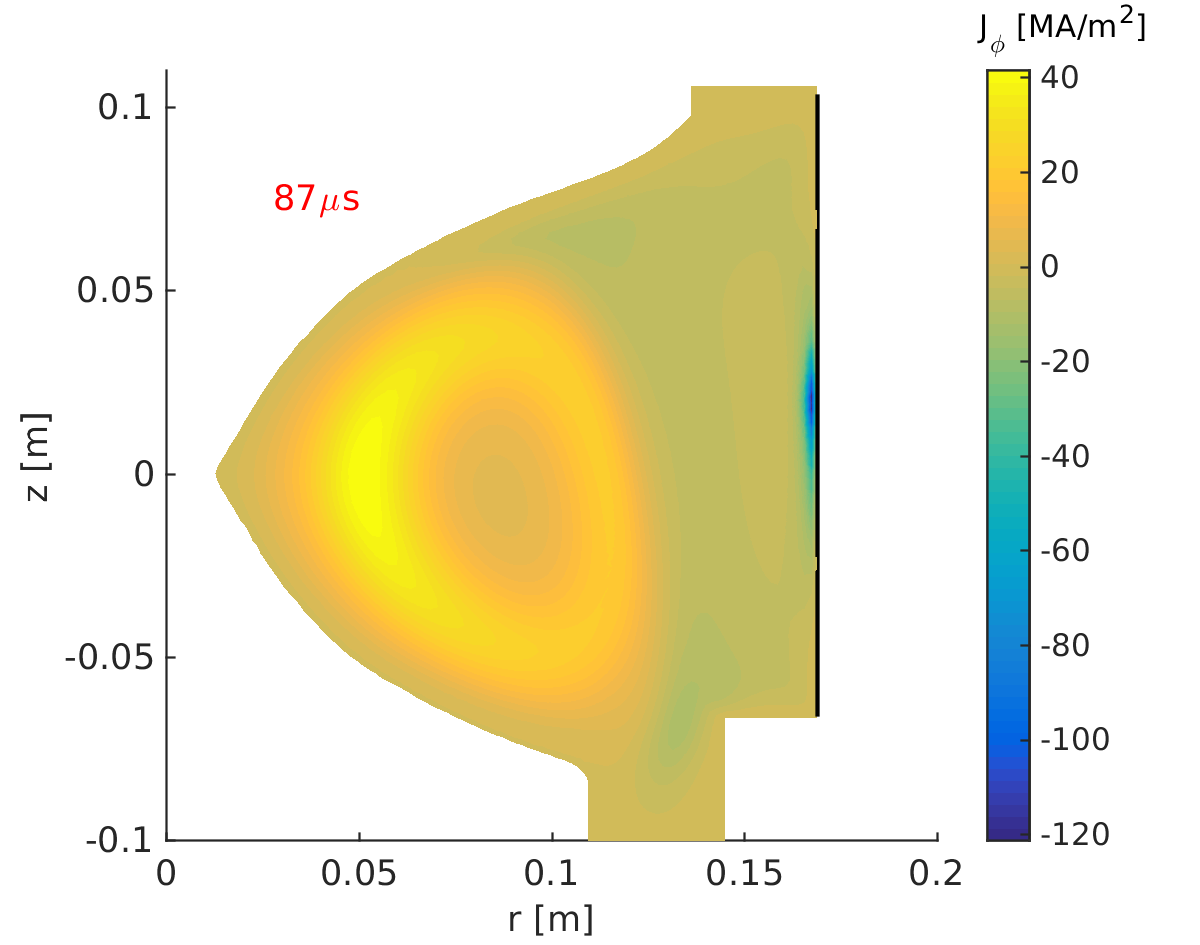}}%
\end{minipage}

\caption{\label{fig: J_11coils}Toroidal current density profiles at $5\upmu$s
(a), $9\upmu$s (b), $18\upmu$s (c), $55\upmu$s (d), $65\upmu$s
(e), $87\upmu$s (f)}
\end{figure}
 Simulated toroidal current density is shown at various times in figure
\ref{fig: J_11coils}. Comparing with figure \ref{fig: Te_11coils},
it can be seen how ohmic heating is a principal electron heating mechanism.
Note the high concentration of $J_{\phi}$ near the entrance to the
CT confinement region at 55$\upmu$s (figure \ref{fig: J_11coils}(d)),
relating to the toroidally directed current sheet present between
oppositely directed poloidal field lines during the magnetic reconnection
process that occurs as closed CT field lines that extend down the
gun, and then open field lines surrounding the CT, are pinched off
during magnetic compression. A reconnection-related current sheet
is also evident when open poloidal field lines reconnect to form closed
CT flux surfaces, for example at 18$\upmu$s (figure \ref{fig: J_11coils}(c)).
Toroidal current, and hence ohmic heating of the electrons, increases
significantly over the compression cycle, as seen in figures \ref{fig: J_11coils}(d)
and (e). 
\begin{figure}[H]
\begin{minipage}[t]{0.4\columnwidth}%
\subfloat[]{\raggedright{}\includegraphics[scale=0.4]{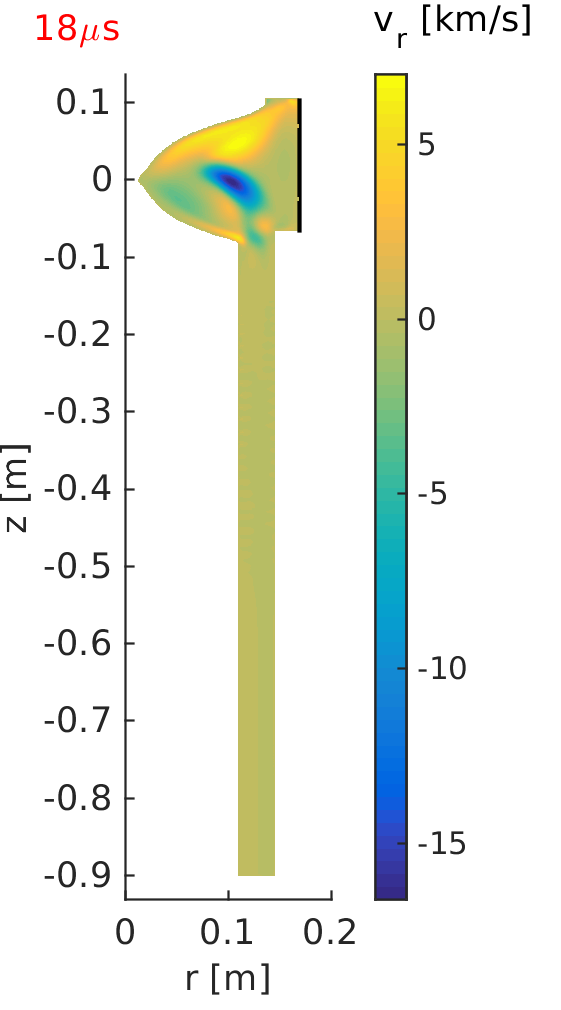}}\hfill{}\subfloat[]{\raggedright{}\includegraphics[scale=0.4]{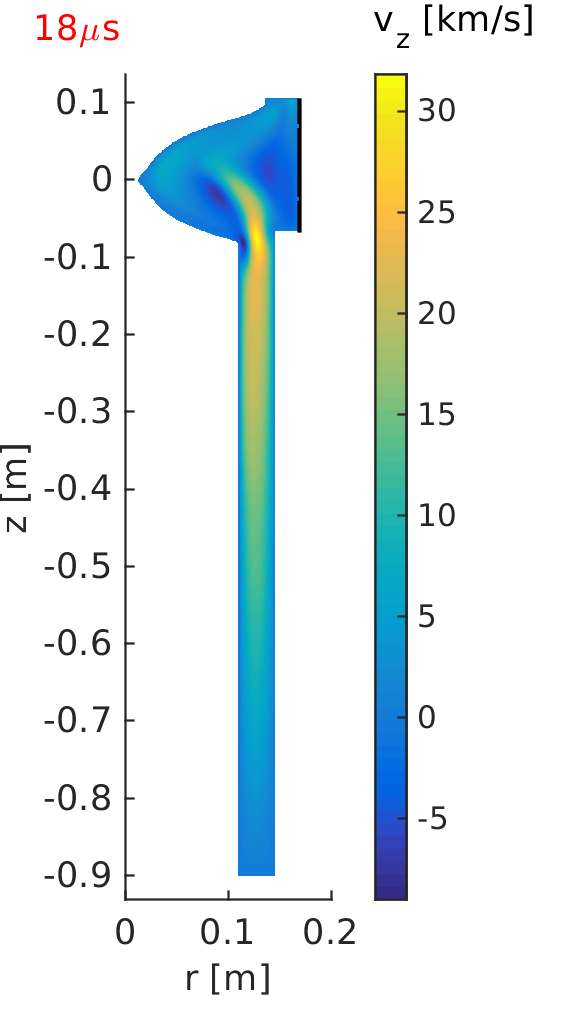}}%
\end{minipage}\hfill{}%
\begin{minipage}[t]{0.2\columnwidth}%
\subfloat[]{\raggedright{}\includegraphics[scale=0.2]{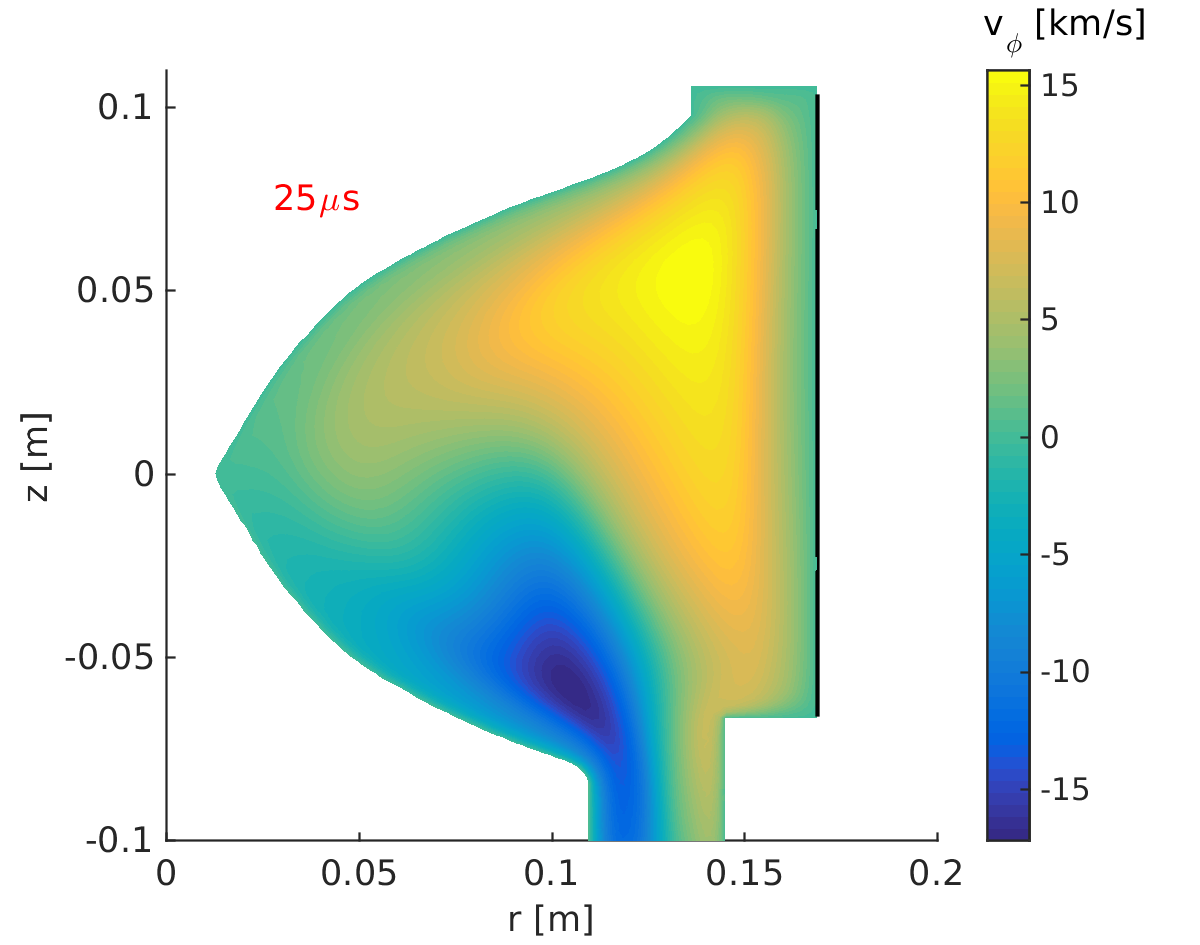}}\hfill{}\subfloat[]{\raggedright{}\includegraphics[scale=0.2]{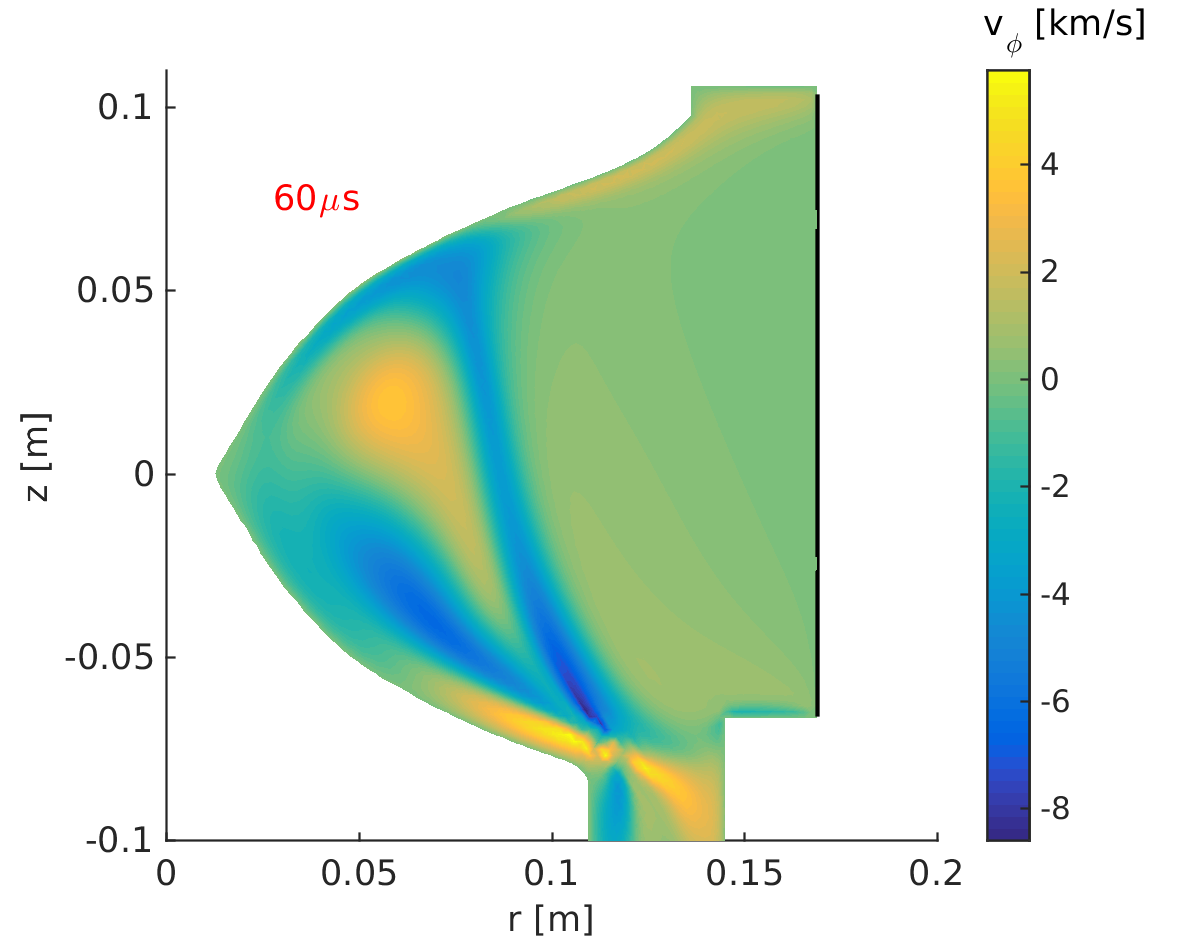}}%
\end{minipage}\hfill{}%
\begin{minipage}[t]{0.2\columnwidth}%
\subfloat[]{\raggedright{}\includegraphics[scale=0.2]{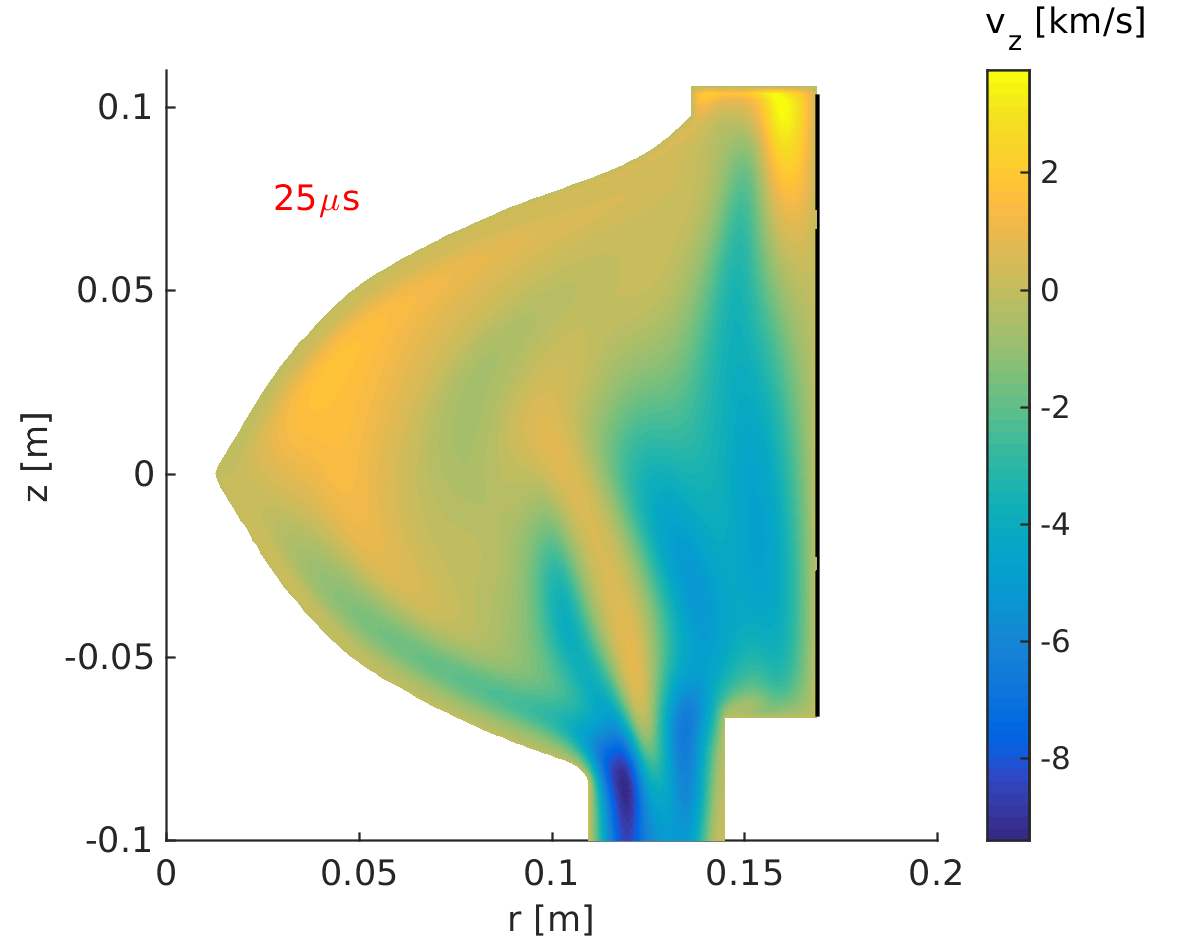}}\hfill{}\subfloat[]{\raggedright{}\includegraphics[scale=0.2]{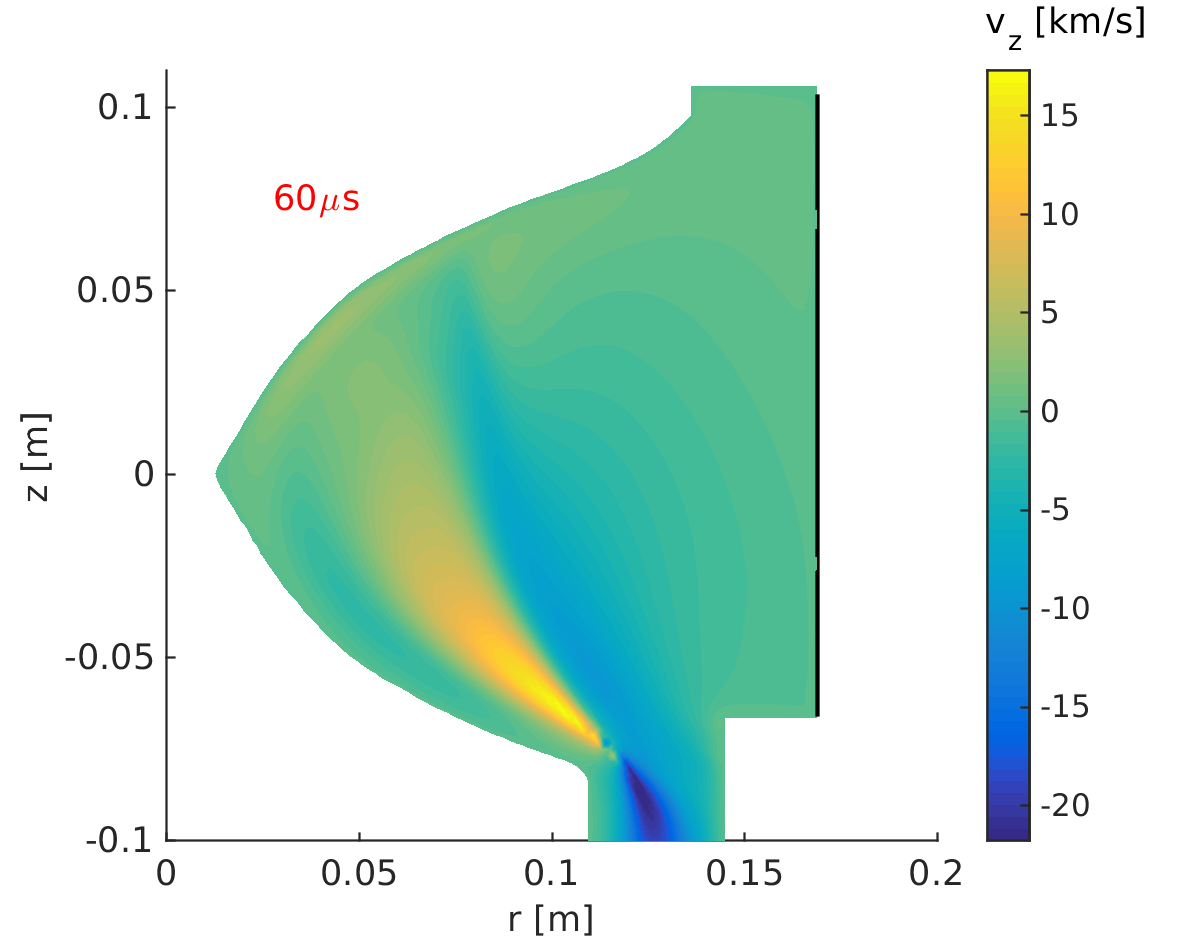}}%
\end{minipage}

\caption{\label{fig: V_11coils}Velocity component profiles: radial velocity
component at $18\upmu$s (a), azimuthal velocity component at $25\upmu$s
(c) and $60\upmu$s (d), axial velocity component at $18\upmu$s (b),
$25\upmu$s (e), and $60\upmu$s (f)}
\end{figure}
Profiles of $v_{r}$ and $v_{z}$ at $18\upmu$s are indicated in
figure \ref{fig: V_11coils}(a) and (b). The dynamics are complicated,
but some key features can be noted. Velocities, in particular $v_{z}$,
are especially high around this time during the CT formation process,
and lead to enhanced ion viscous heating ($cf$ figure \ref{fig: Ti_11coils}(c)).
Pre-compression profiles of $v_{\phi}$ and $v_{z}$ are shown at
$25\upmu$s in figures \ref{fig: V_11coils}(c) and (e).  The particularly
high downward-directed axial velocity near the entrance to the CT
containment region at $25\upmu$s in figure \ref{fig: V_11coils}(e)
is related to jets of plasma fluid associated with magnetic reconnection
during the formation of closed CT flux surfaces. Velocities are gradually
dissipated by viscous effects until $t_{comp}=45\upmu$s, when magnetic
compression is initiated and velocities rise sharply in reaction to
the various dynamics associated with compression. The structures in
the profiles of $v_{\phi}$ and $v_{z}$ near the entrance to the
CT containment region at 60$\upmu$s (figures \ref{fig: V_11coils}(d)
and (f)), just after halfway through the main compression cycle, are
related to jets of plasma fluid associated with magnetic reconnection
when open poloidal field lines surrounding the CT are pinched off.

\subsection{Comparison of simulated and experimental diagnostics\label{subsec:Comparison-of-simulated}}

\begin{figure}[H]
\subfloat[$B_{\theta}$ for shot  39475]{\raggedright{}\includegraphics[width=8cm,height=5cm]{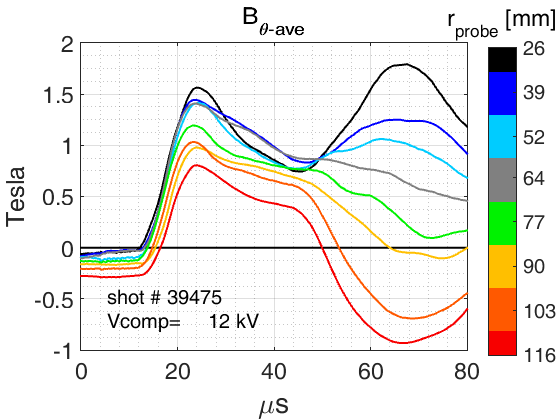}}\hfill{}\subfloat[$B_{\theta}$ for simulation  2350]{\raggedleft{}\includegraphics[width=8.5cm,height=5cm]{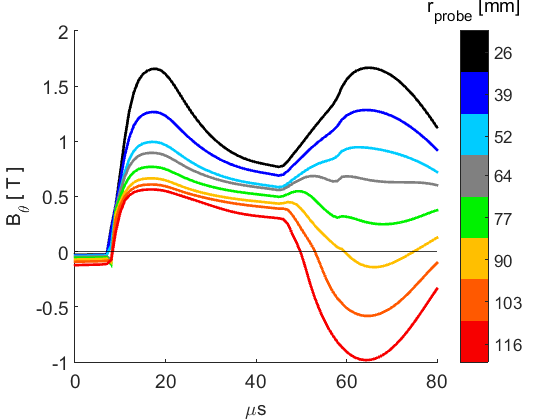}}

\caption{\label{fig:Bpol_meas_cf_sim39475}Comparison of measured (a) and simulated
(b) poloidal magnetic field at magnetic probe locations ($V_{comp}=$12kV,
11-coil configuration). Poloidal magnetic field signals are colored
by the radial coordinates of the probe locations. The magnetic probe
locations are indicated in figure \ref{fig:Grid-arrangement-with}(c).}
\end{figure}

\begin{figure}[H]
\subfloat[Measured $B_{\phi}$ ]{\raggedright{}\includegraphics[width=8.6cm,height=5cm]{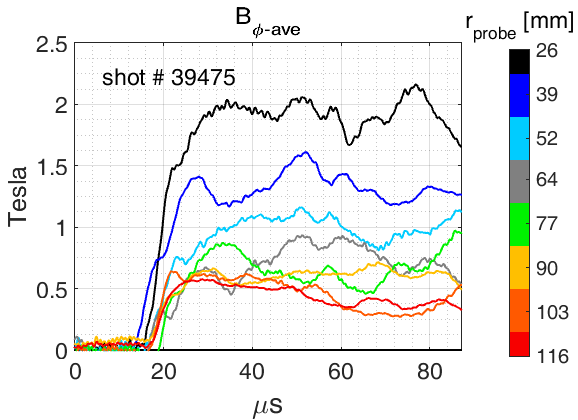}}\hfill{}\subfloat[Simulated $B_{\phi}$]{\raggedleft{}\includegraphics[width=7.5cm,height=5cm]{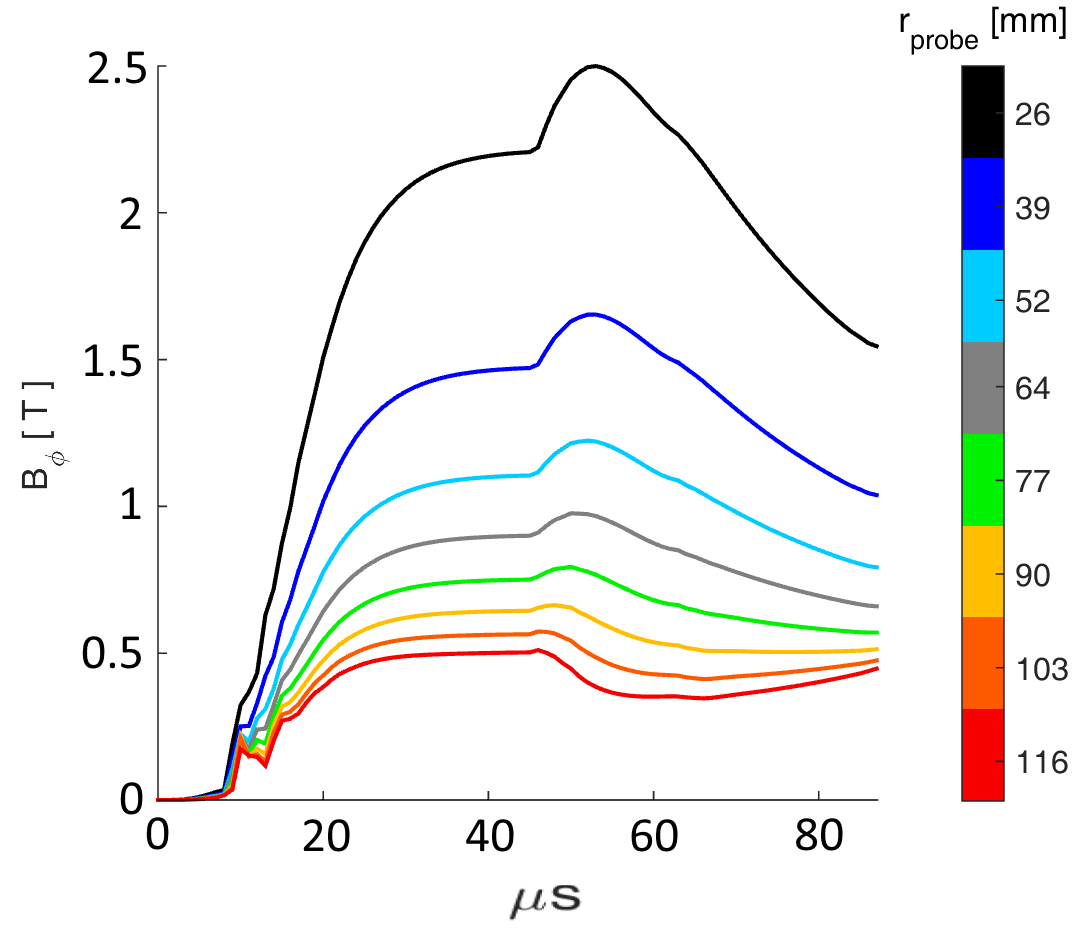}}

\caption{\label{fig:Bphi_exp_sim_comp39475_2350}Comparison of measured (a)
and simulated (b) toroidal magnetic field at magnetic probe locations
($V_{comp}=$12kV, 11-coil configuration)}
\end{figure}

\begin{figure}[H]
\subfloat[Simulated and measured $n_{e}$]{\raggedright{}\includegraphics[width=7cm,height=5cm]{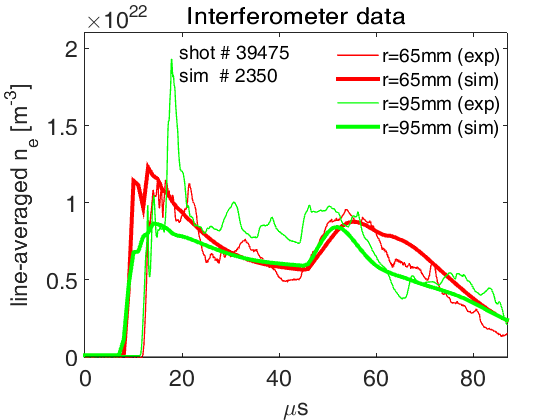}}\hfill{}\subfloat[Simulated and measured $T_{i}$]{\raggedleft{}\includegraphics[width=7cm,height=5cm]{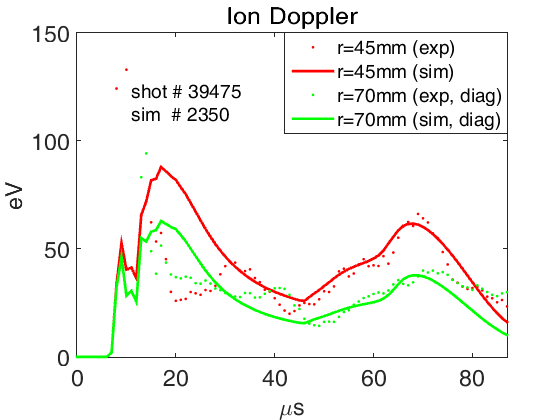}}

\caption{\label{fig:ne_IDcomp}Comparison of measured and simulated line-averaged
electron density (a) and line-averaged ion temperature (b) ($V_{comp}=$12kV,
11-coil configuration)}
\end{figure}

Note that in figures \ref{fig:Bpol_meas_cf_sim39475} to \ref{fig:ne_IDcomp},
the experimental diagnostics presented are from a simulation with
the same code input parameters as those specified above, with the
exception that $V_{comp}=12$kV instead of 18kV. 

A comparison of experimentally measured and simulated $B_{\theta}$
is shown in figure \ref{fig:Bpol_meas_cf_sim39475}. For this shot
(and simulation), $V_{comp}=12$kV and $t_{comp}=45\upmu$s. Note
that in the experiment, multiple magnetic probes were located at the
same radial coordinate but at different toroidal angles. For ease
of comparison, the magnetic field traces experimentally measured at
the same probe radii are averaged in figure \ref{fig:Bpol_meas_cf_sim39475}(a).
The locations of the magnetic probes (which record poloidal as well
as toroidal field) are indicated in figure \ref{fig:Grid-arrangement-with}(c).
Note that measured $B_{\theta}$ is the field component parallel to
the surface of the inner flux conserver in the poloidal plane. Simulated
poloidal field is calculated at the boundary nodes by calculating
$t_{r}$ and $t_{z}$, the $r$ and $z$ components of the unit tangents
to the computational domain at the boundary nodes, and using the equation
\[
\left(B_{\theta}\right)_{bn}=\left(\hat{\mathbf{t}}\cdot\mathbf{B}_{\theta}\right)_{bn}=\left(t_{r}B_{r}+t_{z}B_{z}\right)_{bn}
\]
Here, the subscript $bn$ denotes a particular boundary node. The
method using unit tangents allows for consistent evaluation of the
sign of the poloidal field, which varies across the separatrix between
the CT and levitation/compression field. Referring to figure \ref{fig:Bpol_meas_cf_sim39475}(a),
for shot  39475, $t_{lev}=-50\upmu$s, so with a current rise time
of $\sim40\upmu$s in the levitation coils, the poloidal levitation
field measured at the probes reaches its maximum negative value at
$t=\sim-10\upmu$s. Formation capacitors are fired at $t=0\upmu$s
and it takes $\sim15\upmu$s for the gun (stuffing) flux to be advected
up to the probe locations. Until this time the poloidal field measured
at the probes is purely the levitation field. When the magnetised
plasma bubbles into the CT confinement region, it displaces the levitation
field, and the polarity of the measured field reverses - the stuffing
field has opposite polarity to the levitation field, $i.e.,$ the
toroidal currents in the main solenoidal coil and the levitation coils
are in opposite directions. Over the next several tens of $\upmu$s,
during and after reconnection of poloidal field to form closed flux
surfaces, the CT resistively decays and undergoes Taylor relaxation
during which part of the poloidal flux is converted to toroidal flux.
Note that when levitation or compression field is being measured at
the probes, $|B_{\theta}|$ is larger at the outer probes, due to
the $1/(r_{coil}-r_{probe})$ scaling of levitation (compression)
field with levitation (compression) current in the external coils.
On the other hand, when CT field is being measured at the probes,
$B_{\theta}$ is larger at the inner probes, due to the $1/r_{probe}^{2}$
scaling of CT field with CT flux - poloidal field lines are bunched
together progressively more at smaller radii. The compression capacitors
are fired at $t=t_{comp}=45\upmu$s, and the total compression current
in the external coils rises over $\sim20\upmu$s to its peak, for
$V_{comp}=12$kV, of around 850kA, so that the total combined levitation
and compression current is around $1$MA at the time peak compression,
around $65\upmu$s. In this shot, the CT is compressed inwards beyond
the probes at $r=90\mbox{mm}$, so $B_{\theta}$ at the probes located
at $r\geq90\mbox{mm}$ ($i.e.,$ light orange, dark orange, and red
traces) is, starting at various times after $45\upmu$s, a measurement
of the compression field, while the CT poloidal field is measured
at the probes with radial coordinates $r<90\mbox{mm}$. The CT remains
stable during compression in this shot, and it expands approximately
to its pre-compression state (apart from resistive flux losses and
thermal losses) between $t\sim65\upmu$s and $t\sim87\upmu$s, when
the compression current falls to zero and changes direction. Referring
to figure \ref{fig:Bpol_meas_cf_sim39475}(b), axisymmetric MHD simulations
allow for only resistive loss of flux and do not capture the mechanisms,
discussed in \cite{thesis,exppaper}, that led to poloidal flux loss
in many compression shots. Shot $39475$ was a flux-conserving shot,
and a reasonable match is found between experimentally inferred and
MHD-simulated poloidal field. 
Experimentally measured $B_{\phi}$ for shot 39475 is shown in figure
\ref{fig:Bphi_exp_sim_comp39475_2350}(a). As discussed in \cite{thesis,exppaper},
$B_{\phi}$ at the magnetic probes rises at compression as crow-barred
shaft current flowing in conductors ($i.e.,$ machine components as
well as ambient plasma) surrounding the CT increases when it is able
to divert from the aluminum bars outside the insulating wall to a
lower inductance path through ambient plasma located outboard of the
CT. Due to an instability (thought to be an external kink) that was
routinely observed at compression, the current diversion, and hence
toroidal magnetic field measured at probes located at the same radii
and different toroidal angles, is generally toroidally asymmetric.
Measured toroidal field is generally observed to increase over compression
at probes located inside the CT boundary, and decrease at probes located
outside the boundary; the CT outer boundary is also toroidally asymmetric
at compression. The toroidal magnetic field traces experimentally
measured at the same radii at different toroidal angles have been
averaged in figure \ref{fig:Bphi_exp_sim_comp39475_2350}(a). The
asymmetric MHD code cannot capture the inherently 3D instability.
However, referring to figure \ref{fig:Bphi_exp_sim_comp39475_2350}(b),
toroidal field increases more at compression at the inner probes and
decreases at the outer probes at $r\geq90\mbox{mm}$, as shaft current
is diverted to a path inside the locations of the outer probes, which
is also qualitatively evident from the toroidal-averages of the measured
toroidal field. Simulated diversion of shaft current to a path through
ambient plasma located outboard of the CT is apparent from contours
of $f$ at $65\upmu$s, see figure \ref{fig: F_11coils}(d).
The comparison between experimentally measured and simulated electron
density and ion temperature for shot 39475 is indicated in figures
\ref{fig:ne_IDcomp}(a) and \ref{fig:ne_IDcomp}(b). The simulated
line-averaged electron density along the interferometer chord at $r=35$mm
(see figure \ref{fig: n_11coils}) hasn't been included in figure
\ref{fig:ne_IDcomp}(a) because the experimental data for that chord
is not available. The experimentally measured electron densities are
line-averaged quantities obtained with He-Ne laser interferometers
looking down the vertical chords at $r=65\mbox{\mbox{mm}}$ and $r=95\mbox{\mbox{mm}}$
that are indicated in figure \ref{fig: n_11coils}. The time it takes
for plasma to be advected up the gun and enter the CT containment
region is underestimated by $\sim5\upmu$s in this simulation. The
density magnitude during CT decay and magnetic compression is approximately
reproduced. Referring to figure \ref{fig:ne_IDcomp}(b), the experimentally
measured ion temperatures were evaluated from Doppler broadening of
line radiation from singly ionized Helium (He II line at 468.5nm)
along the chords indicated in figure \ref{fig: Ti_11coils}. Note
that the diagonal green coloured chord indicated in figure \ref{fig: Ti_11coils}
has its lower point at $r=70$mm. Simulated line-averaged ion temperature
is evaluated along the same chords, producing a reasonable match to
the experimental data. For this shot (and simulation), with $V_{comp}=12$kV,
an increase in ion temperature by a factor of over two, from $\sim25$eV
to $\sim60$eV, is indicated in the region of the ion Doppler chords.
As outlined in \cite{thesis,exppaper}, a maximum error in the temperature
measurement due to density broadening has been evaluated as $\sim$10eV
for density levels associated with shot 39475 at peak compression.
Careful analysis was undertaken to confirm that temperature broadening
rather than density broadening was the dominant broadening mechanism. 

\subsection{Simulated plasma-wall interaction\label{subsec:Simulated-plasma-wall-interactio}}

\begin{figure}[H]
\begin{raggedright}
\subfloat[6-coil configuration]{\raggedleft{}\includegraphics[width=7cm,height=5cm]{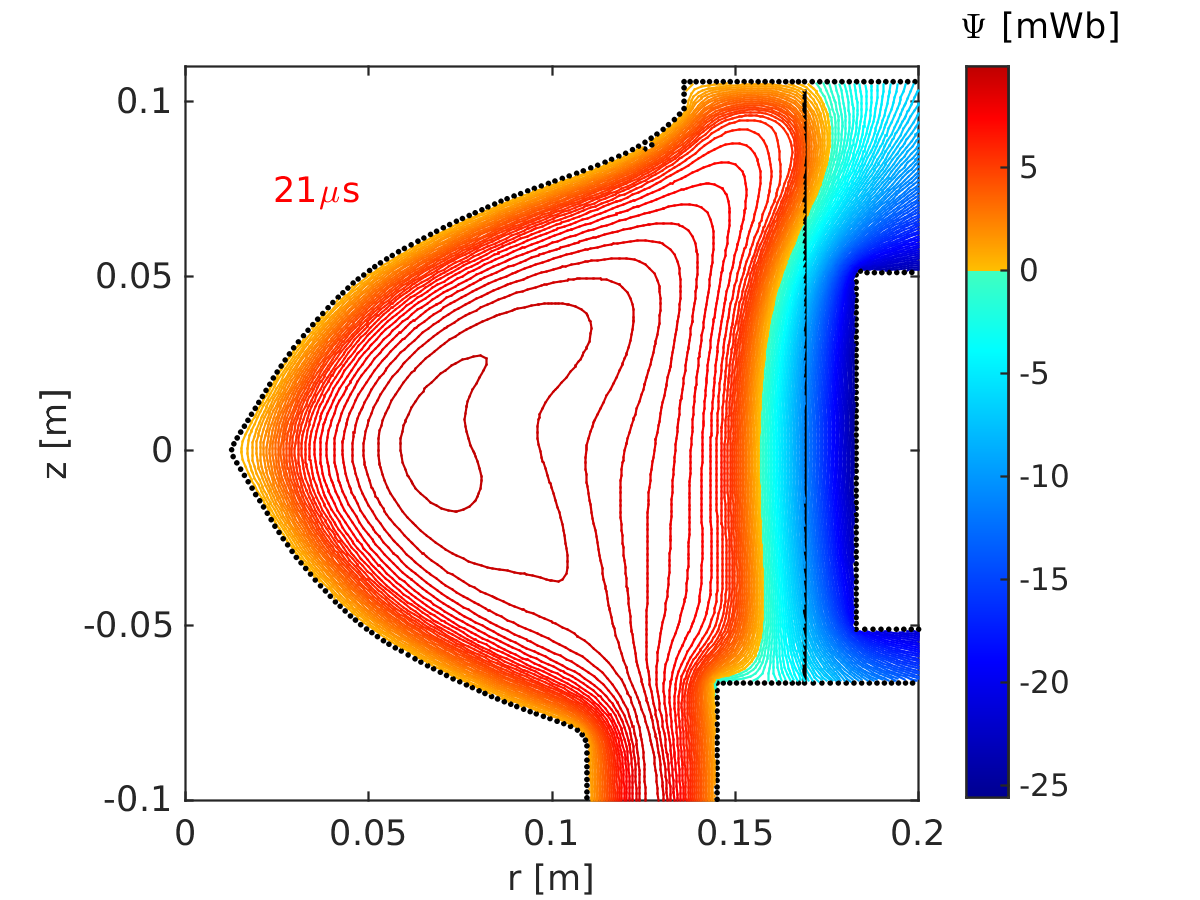}}\hfill{}\subfloat[11-coil configuration]{\raggedleft{}\includegraphics[width=7cm,height=5cm]{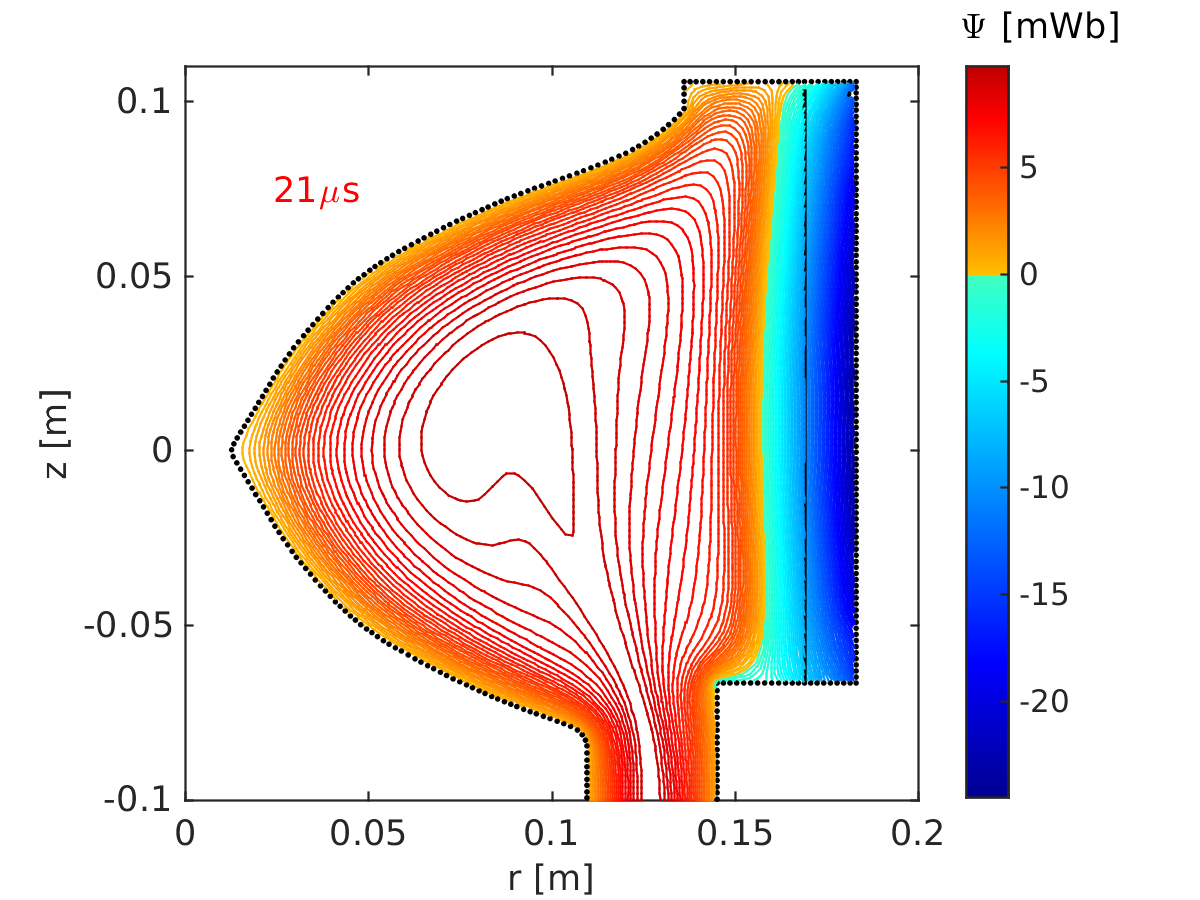}} 
\par\end{raggedright}
\centering{}\caption{\label{fig:plasmaWall}Simulated plasma-wall interaction. Note how
poloidal field penetrates the insulating wall during the bubble-in
process in the six coil configuration.}
\end{figure}
As described in \cite{thesis,exppaper}, levitated CT lifetime and
temperature was increased significantly when the original configuration
consisting of six levitation/compression coils was replaced with an
eleven coil configuration. It is believed that the improvement is
related to reduced levels of levitation field displacement and interaction
between the plasma and insulating wall during the CT formation process,
leading to consequent reductions in sputtering, plasma impurity concentration,
and radiative cooling. Spectrometer data indicated a four-fold reduction
in total spectral power with the 11-coil configuration, even at increased
formation power, confirming that the 11-coil configuration led to
reduced impurity concentration. The inclusion of a model for an insulating
region representing the wall and air located between the vacuum chamber,
in which MHD dynamics are evolved, and the stack of external levitation/compression
coils, enables simulations to provide further verification of the
proposed explanation for the observed improvement. MHD simulations
confirm the reduction of plasma-wall interaction with the eleven coil
configuration, as shown in figure \ref{fig:plasmaWall}. In figure
\ref{fig:plasmaWall}(a), the stack of six coils is partly located
in the blank rectangle on the right, centered around $z=0\mbox{cm}$,
and extends off further to the right (not shown). The region above,
below, and just to the left of the coil-stack represents the air around
the stack, while the area to the right of the vertical black line
at $r=17\mbox{cm}$ represents insulating material. The vertical black
line represents the inner radius of the insulating wall, and the outer
radius of the insulating wall at $r=17.7\mbox{cm}$ is not indicated.
In figure \ref{fig:plasmaWall}(b), the stack of eleven coils extends
all the way from the top to the bottom of the insulating wall, the
inner radius of the coil stack is the same as that for the six coil
stack. In both cases, only $\psi$, which determines the vacuum poloidal
field, is evaluated in the insulating region to the right of the inner
radius of the insulating wall. The solution for $\psi$ is coupled
to the full MHD solution in the remainder of the domain, as described
in section \ref{sec:Vacuum-field-in}. To maintain toroidal flux conservation,
boundary conditions for $f$, which has a finite constant value in
the insulating wall and is zero outside the current-carrying aluminum
bars depicted in figure \ref{fig:Machine-Schematic}(b), are evaluated
for, and applied to, the part of the boundary representing the inner
radius of the insulating wall, as described in section \ref{sec:PHIconservation-with}.
Both simulations have boundary conditions for $\psi_{main}$ and $\psi_{lev}$
from FEMM models, pertaining to $I_{main}=70$A, and with the total
levitation current such that $\psi_{lev}$ is approximately the same
for each configuration. Also, code input $V_{form}$ is set to 16kV
for both simulations. Figure \ref{fig:plasmaWall}(a) indicates how
poloidal field penetrates the insulating wall during the bubble-in
process in the six coil configuration. In practice, ions streaming
along and gyro-rotating around the field lines would then sputter
insulating material into the plasma, leading to impurity radiation
and radiative cooling, with consequent increased resistivity and reduced
CT magnetic lifetimes. 

\subsection{Compression field reversal\label{subsec:Compression-field-reversal}}

\begin{figure}[H]
\begin{centering}
\subfloat[$B_{\theta}$, shot 39735 ]{\raggedright{}\includegraphics[width=9cm,height=5cm]{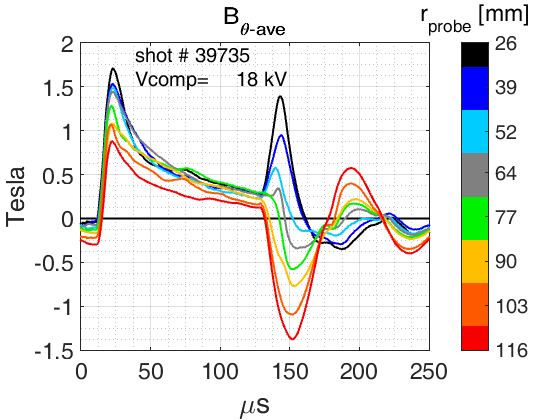}}\hfill{}\subfloat[$B_{\theta}$, simulation 2287]{\raggedleft{}\includegraphics[width=8.8cm,height=4.8cm]{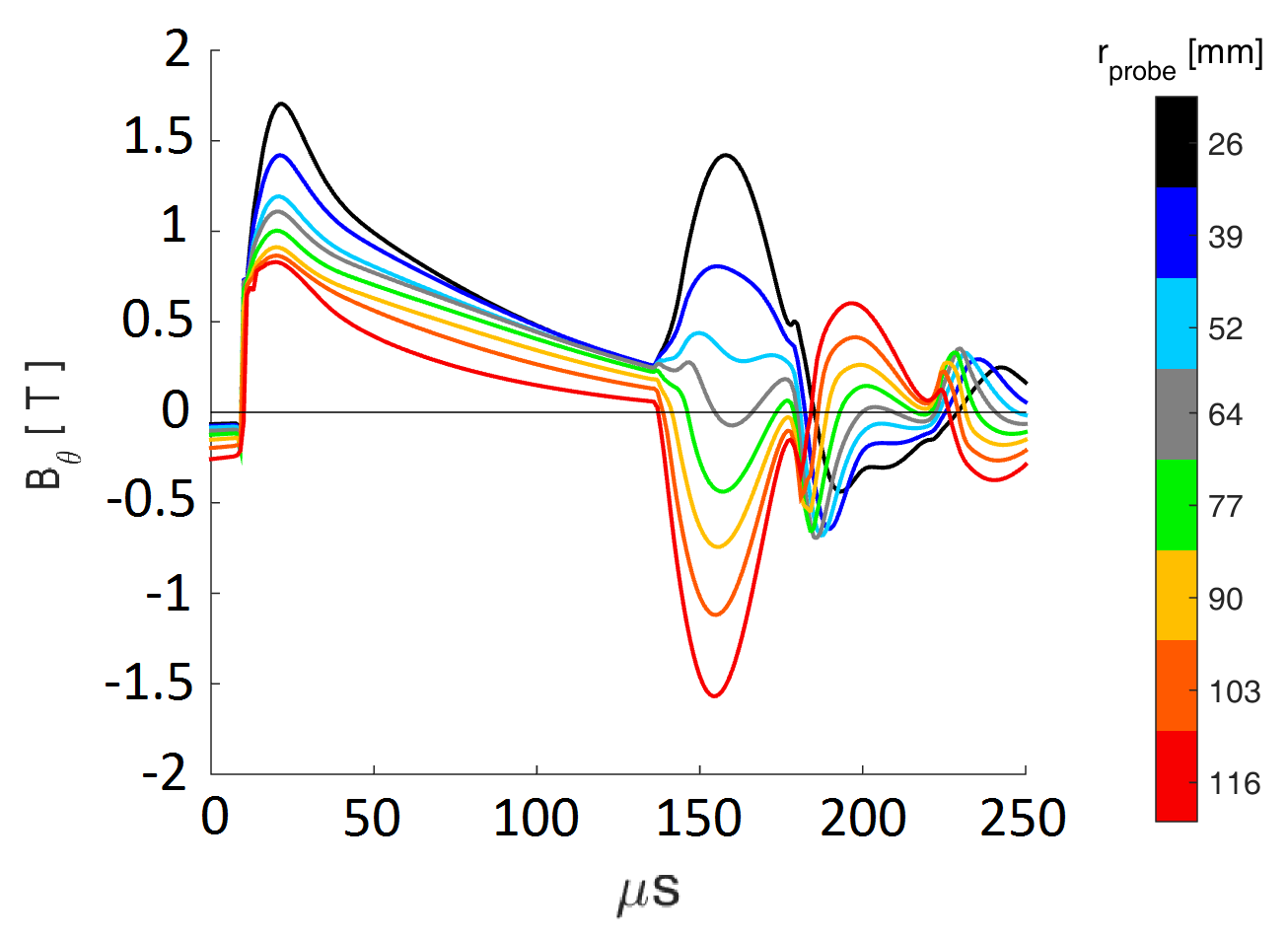}}
\par\end{centering}
\begin{centering}
\subfloat[$B_{\phi}$, shot 39735 ]{\raggedright{}\includegraphics[width=9cm,height=5cm]{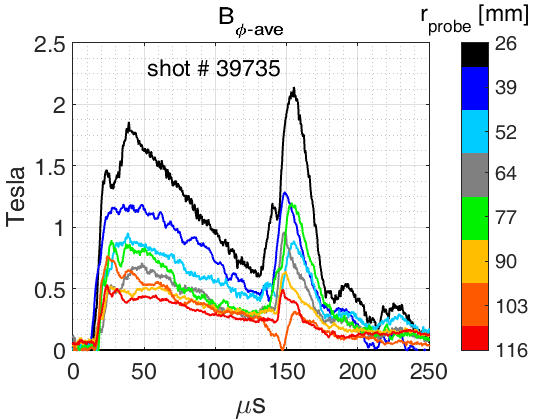}}\hfill{}\subfloat[$B_{\phi}$, simulation 2287]{\raggedleft{}\includegraphics[width=8.5cm,height=4.8cm]{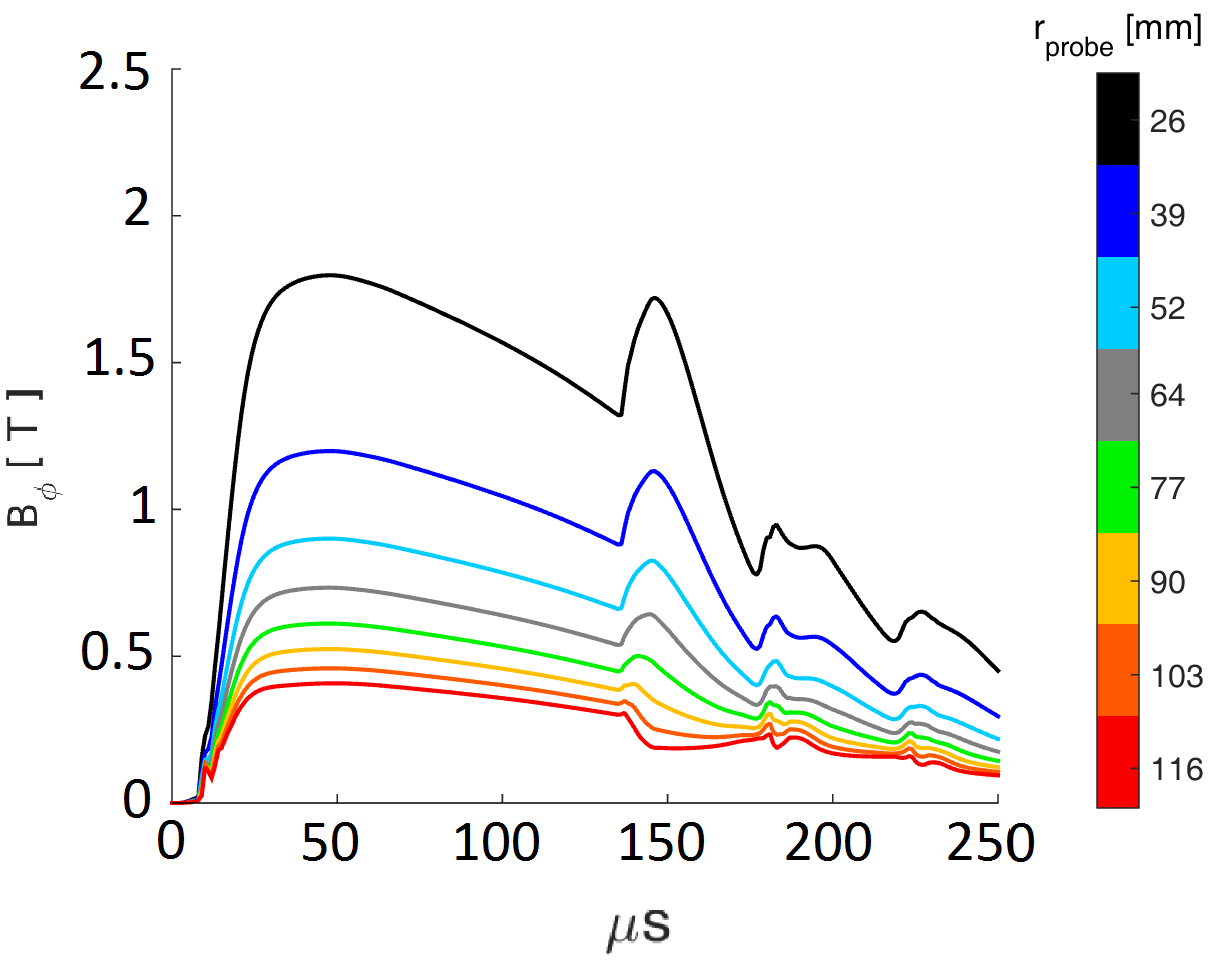}}
\par\end{centering}
\caption{\label{fig:Bpol_Btor_meas_cf_sim39735_2287}Comparison of measured
and simulated poloidal and toroidal magnetic field at magnetic probe
locations ($V_{comp}=18$kV, 11-coil configuration) }
\end{figure}

When the compression current in the coils changes direction (see figure
\ref{fig:Ilev_comp_sim}(b)), the CT poloidal field magnetically reconnects
with the compression field, and a new CT with polarity opposite to
that of the previous CT is induced in the containment region, compressed,
and then allowed to expand. The process repeats itself at each change
in polarity of the compression current until either the plasma loses
too much heat, or the compression current is sufficiently damped.
\begin{figure}[H]
\begin{centering}
\subfloat[]{\raggedright{}\includegraphics[width=5cm,height=4.5cm]{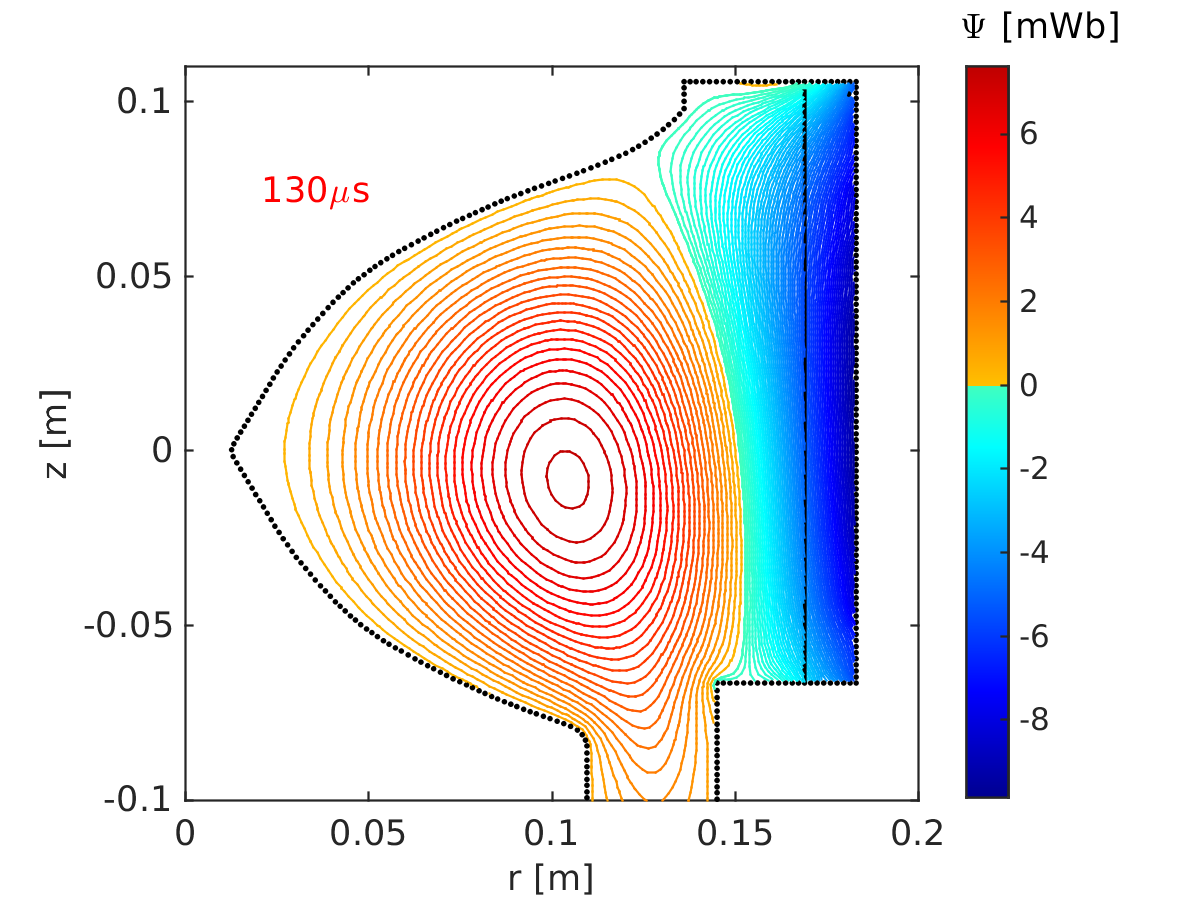}}\subfloat[]{\raggedleft{}\includegraphics[width=5cm,height=4.5cm]{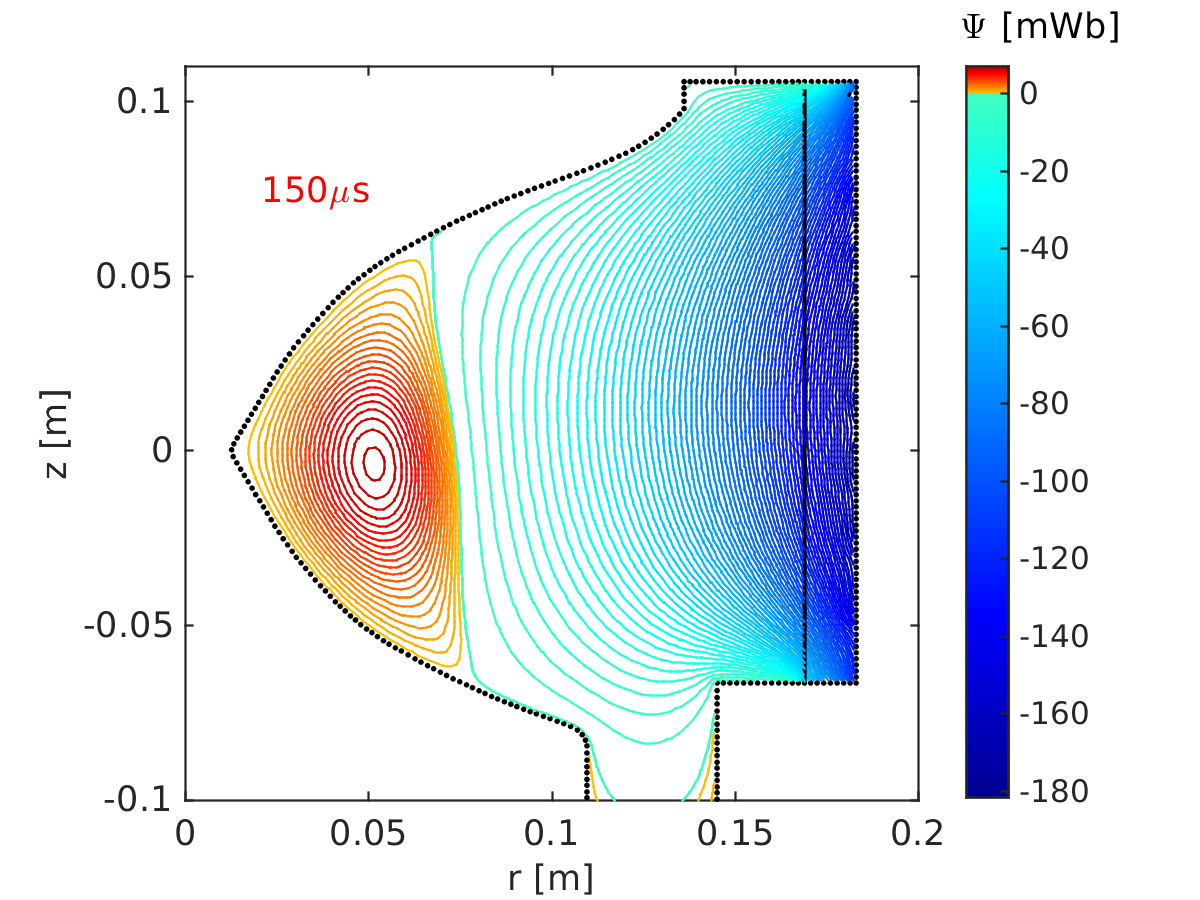}}\subfloat[]{\raggedright{}\includegraphics[width=5cm,height=4.5cm]{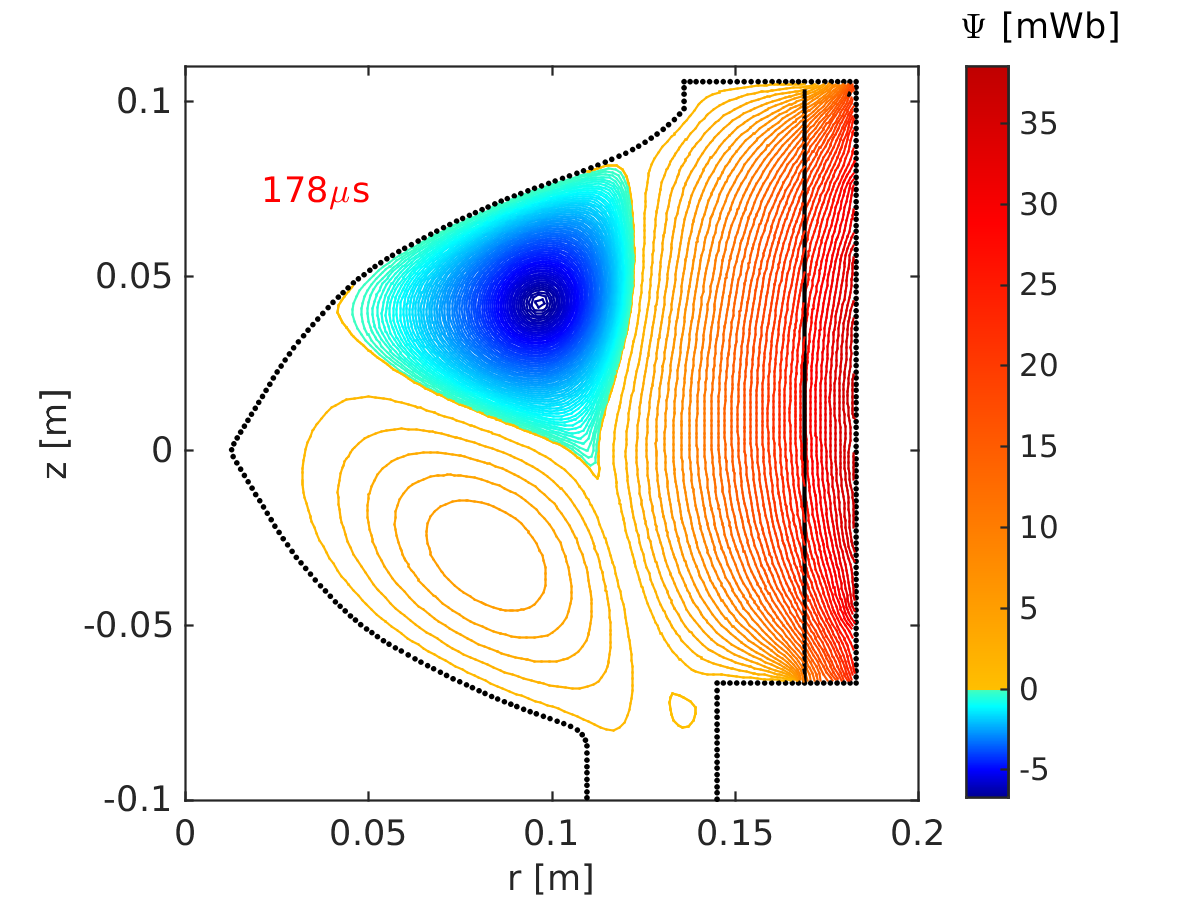}}\subfloat[]{\raggedright{}\includegraphics[width=5cm,height=4.5cm]{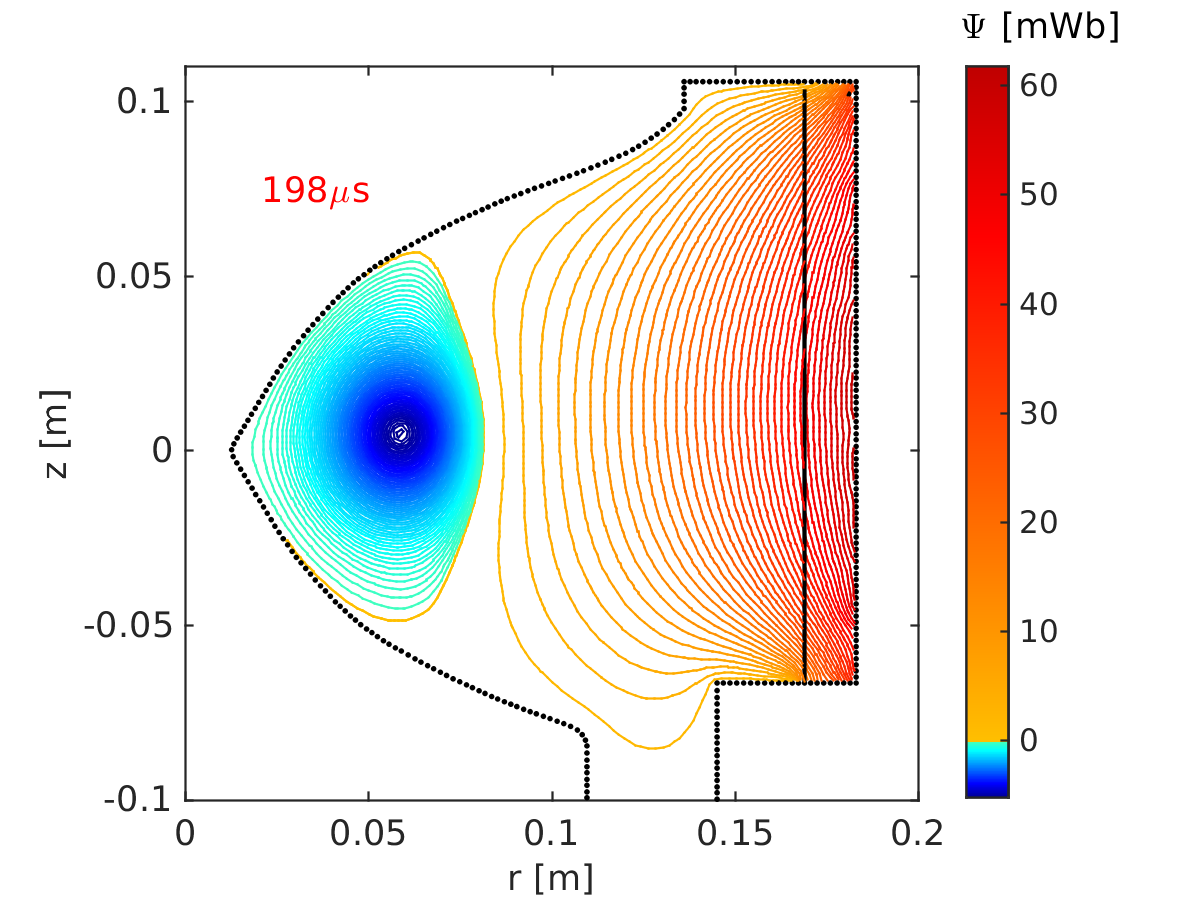}}
\par\end{centering}
\begin{centering}
\subfloat[]{\raggedright{}\includegraphics[width=5cm,height=4.5cm]{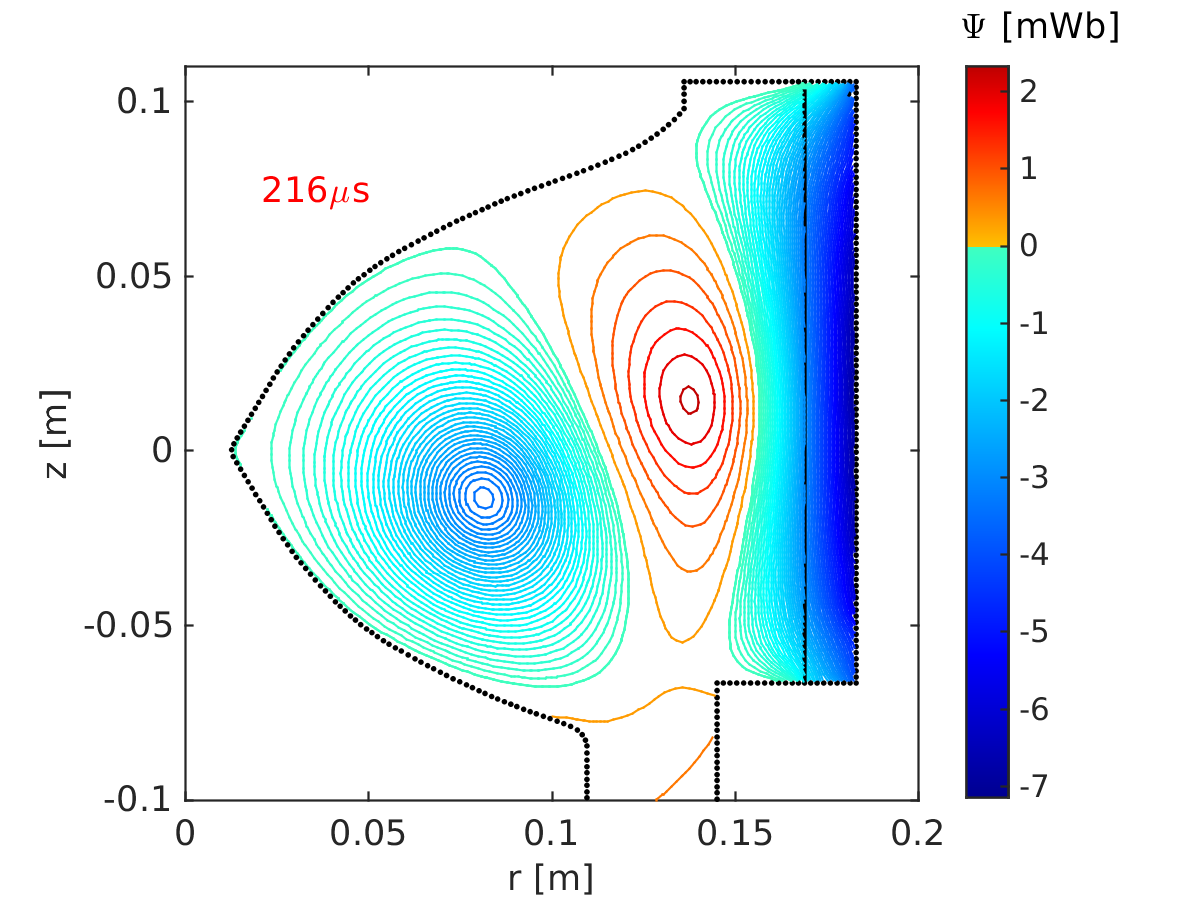}}\subfloat[]{\raggedleft{}\includegraphics[width=5cm,height=4.5cm]{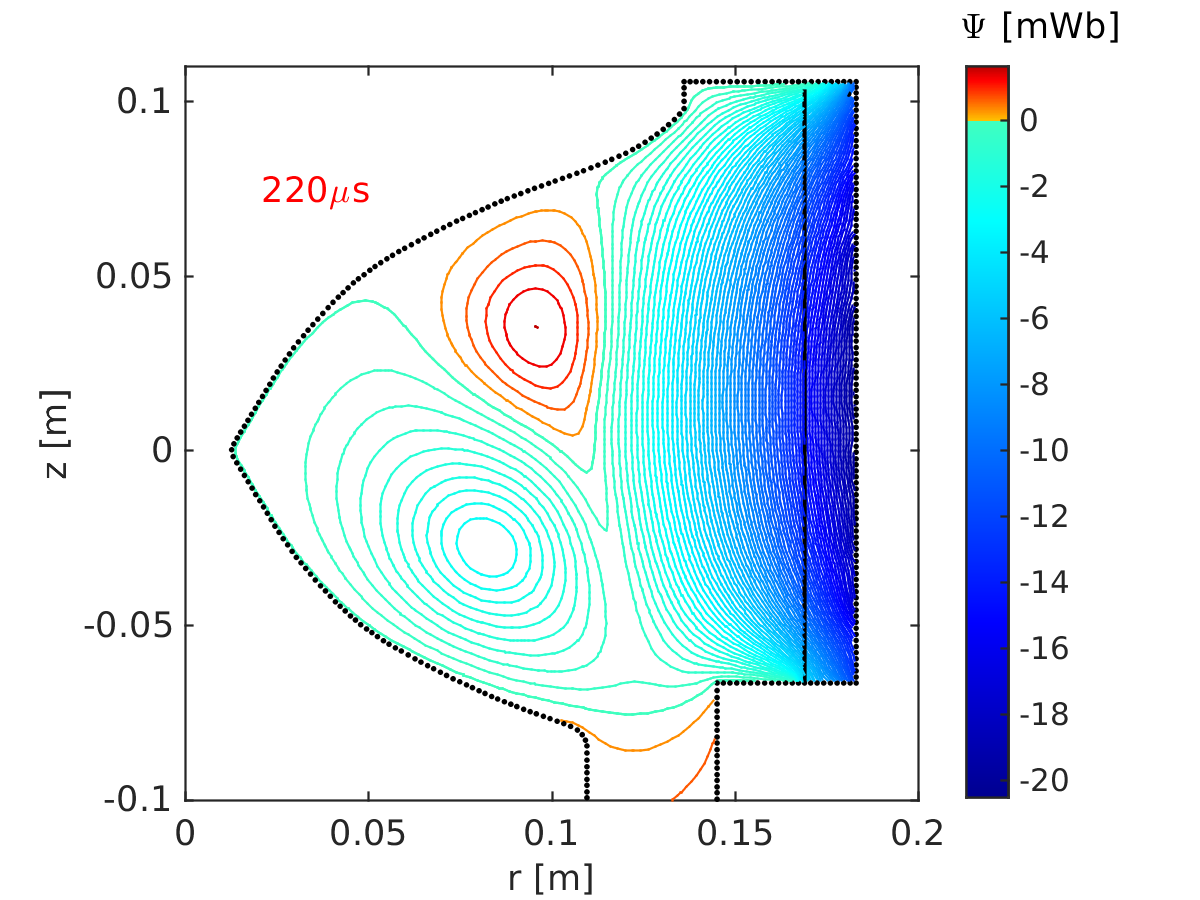}}\subfloat[]{\raggedright{}\includegraphics[width=5cm,height=4.5cm]{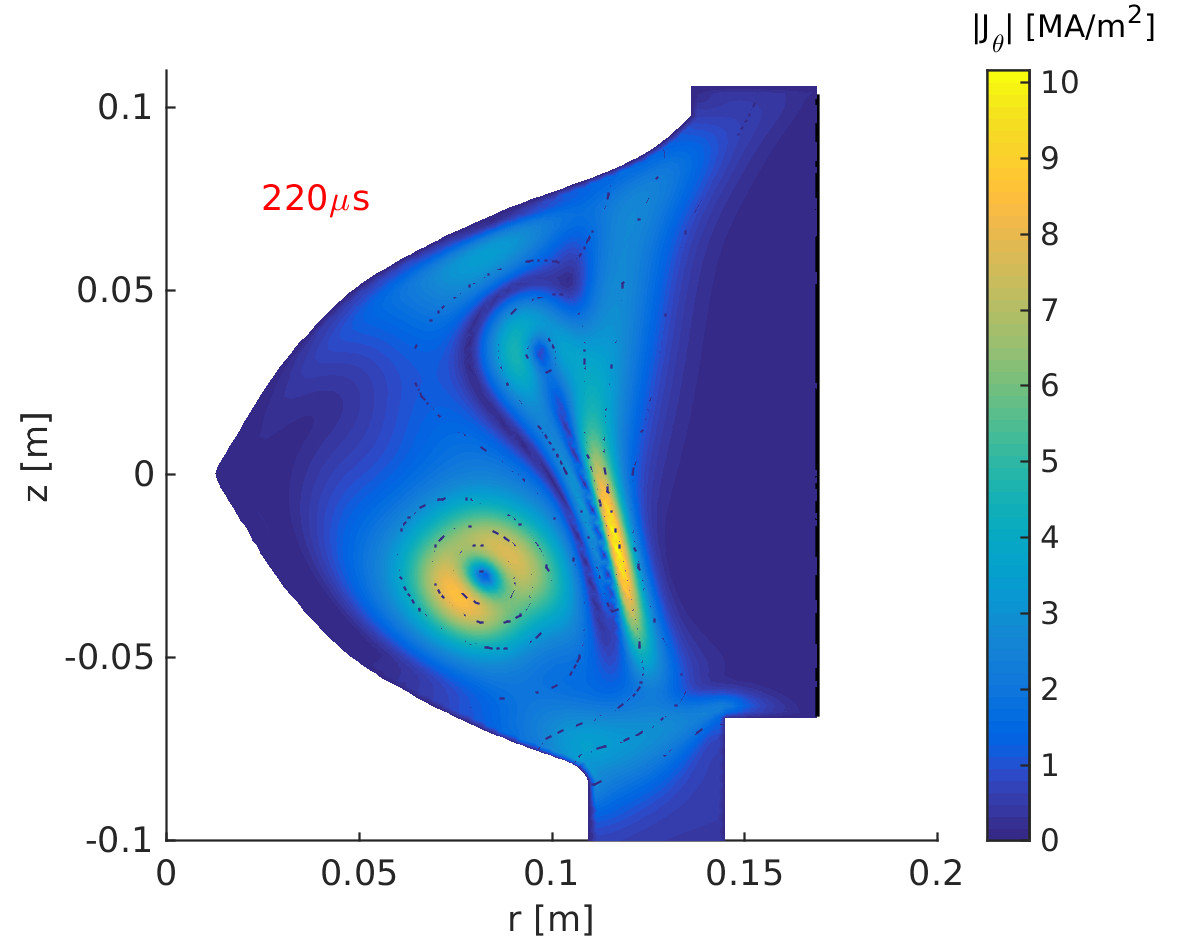}}\subfloat[]{\raggedright{}\includegraphics[width=5cm,height=4.5cm]{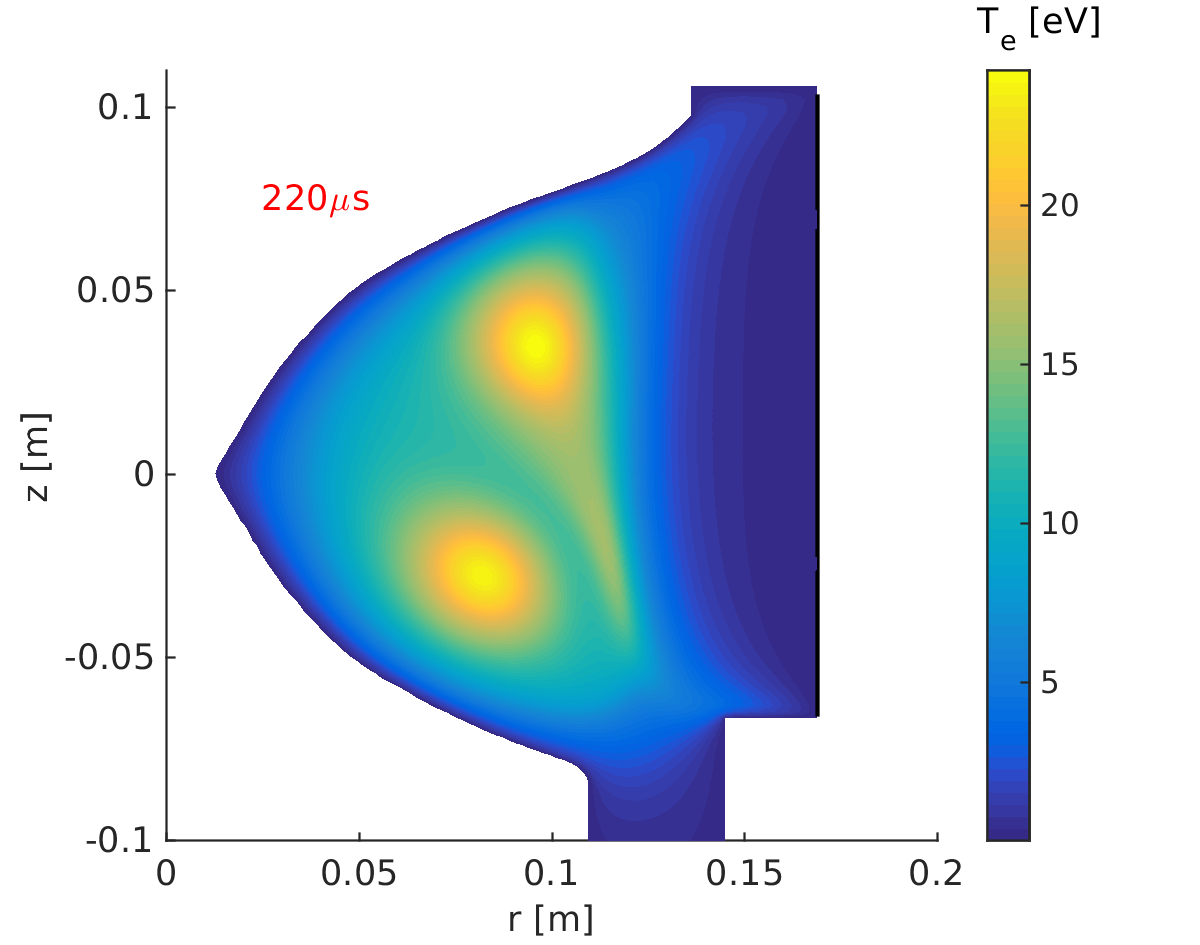}}
\par\end{centering}
\caption{\label{fig:psi_2287}Poloidal flux contours at $130\upmu$s (a), $150\upmu$s
(b), $178\upmu$s (c), $198\upmu$s (d), $216\upmu$s (e), $220\upmu$s
(f), and profiles of poloidal current density (g) and electron temperature
(h) at $220\upmu$s, for an MHD simulation with compression current
reversal}
\end{figure}
Figures \ref{fig:psi_2287}(a) to (f) show $\psi$ contours at various
times for a simulation similar to that presented in figures \ref{fig: psi_11coils}
to \ref{fig: V_11coils}, with the exception that $t_{comp}=130\upmu$s
instead of 45$\upmu$s. Peak compression is at $150\upmu$s (figure
\ref{fig:psi_2287}(b)). By $178\upmu$s, the external compression
field has changed polarity and starts to reconnect with the CT poloidal
field. Toroidal currents are induced to flow in the ambient plasma
initially located outboard of the original CT, enabling the formation
of a new CT (blue closed contours) with polarity opposite to that
of the original CT (figure \ref{fig:psi_2287}(c)). The new induced
CT is magnetically compressed inwards by the increasing reversed polarity
compression field, with peak compression at around $198\upmu$s (figure
\ref{fig:psi_2287}(d)). The compression field polarity rings back
to its original state by $216\upmu$s, when a third CT is induced,
with the same polarity as the original CT (figure \ref{fig:psi_2287}(e)).
By $220\upmu$s, the poloidal field of the second CT is reconnecting
with the compression field, and the third CT is being compressed inwards
(figure \ref{fig:psi_2287}(f)). Electron temperature and poloidal
current density at $220\upmu$s for the two co-existing CTs is presented
in \ref{fig:psi_2287}(g) and (h). Simulated magnetic field measurements
from this simulation are presented and compared with relevant experimental
measurements in figure \ref{fig:Bpol_Btor_meas_cf_sim39735_2287}.

The comparison between experimentally measured and simulated poloidal
and toroidal field, when the compression current reversals are included,
is shown in figure \ref{fig:Bpol_Btor_meas_cf_sim39735_2287}. The
experimental magnetic field measurements have been toroidally averaged
for clarity. Shot  39735 presented here had $V_{comp}=18$kV and $t_{comp}=130\upmu$s.
In shot  39735, the poloidal field measured at the inner probes collapses
at $\sim145\upmu$s (figure \ref{fig:Bpol_Btor_meas_cf_sim39735_2287}(a))
, while the compression coil current peaks at $\sim150\upmu$s, implying
(as presented in detail in \cite{thesis,exppaper}) that poloidal
flux was not well conserved during compression. Apart from resistive
losses, CT poloidal flux is conserved in the simulation, so the poloidal
field at the inner probes (figure \ref{fig:Bpol_Btor_meas_cf_sim39735_2287}(b))
continues to rise until the compression coil current peaks. 

As discussed in previously reference to figure \ref{fig:Bphi_exp_sim_comp39475_2350},
the compressional instability lead to toroidal field measurements
that are toroidally very asymmetric, and the axisymmetric MHD model
cannot reproduce this. Comparison of figures \ref{fig:Bpol_Btor_meas_cf_sim39735_2287}(c)
and (d) shows how the simulated $B_{\phi}$ does, in general, rise
at the magnetic probes as crowbarred shaft current increases when
it is diverted to a lower inductance path. There is a qualitative
agreement between the simulated $B_{\phi}$ and the toroidal-averages
of the measured $B_{\phi}.$

\section{Conclusion\label{sec:Conclusion}}

It has been shown how various differential operators with useful mimetic
properties have been developed, and used to ensure conservation of
total energy, particle count, toroidal flux, and (in some scenarios)
total angular momentum in an axisymmetric linear finite element numerical
scheme that implements the non-linear single fluid, two temperature
MHD equations. The principal mimetic qualities of the operators are
that they satisfy discrete forms of the differential product rule
and the divergence theorem. A novelty of the code is that all discrete
spatial differential operators are represented as matrices, and the
discretized forms of continuous differential equations may be obtained
by simply replacing the original continuous differentiations with
the corresponding matrix operators. The resultant DELiTE framework
may be applied to solve a wide range of systems of differential equations,
and may serve as an educational tool. The methods developed for simulating
CT formation, magnetic levitation and magnetic compression to study
a novel experiment have been outlined. Special care has been taken
to simulate the poloidal vacuum field in the insulating region between
the inner radius of the insulating tube and the levitation/compression
coils, and to couple this solution to the evaluation of $\psi$ in
the plasma domain, while maintaining toroidal flux conservation, enabling
a quantitative model of plasma/wall interaction in various coil configurations.
In general, a reasonable match is obtained between simulated and experimental
diagnostics.

Weak points of the code include the first order accuracy associated
with the linear scheme, and the explicit timestepping scheme, which
limits the range of allowable values for the various diffusion coefficients.
The diffusive nature of the code is associated with requirements for
numerical stability. Presently, angular momentum is conserved only
in particular simulation scenarios.

Further code development may include a more generally applicable model
for density diffusion with correction terms that maintain angular
momentum conservation even in simulation scenarios that include CT
formation and magnetic compression. The level to which artificial
density diffusion affects the simulated density diagnostic, and the
evolution of the other fields, should be investigated further. An
implicit timestepping scheme may be implemented as this would enable
faster simulation times while enabling the use of a greater range
of values for the various diffusion coefficients. 

The ability to model part of the domain as a plasma-free material
is useful because it expands the code's range of applicability as
a problem-solving tool. For example, the feature can easily be adapted
to model the penetration of field associated with magnetically confined
plasma into surrounding conductors. The effort invested in developing
a conservative scheme with many extended features has resulted in
a code framework with a solid foundation for further development and
improvement. For example, as described in \cite{Neut_paper}, a model
of plasma/neutral particle interaction was developed and included
in the DELiTE framework, enabling reproduction of the effect observed
on another experiment, the General Fusion SPECTOR machine, of neutral
gas diffusing up the Marshall gun, leading to markedly increased CT
electron density well after completion of the CT formation process.

\section{Acknowledgments}

Many thanks to so many people at Genera Fusion Inc, (LIST)XXX for
experimental results and general support and simulation validation
tests. Funding was provided in part by General Fusion Inc., Mitacs,
University of Saskatchewan, and NSERC. We would like to thank to Aaron
Froese, Pat Carle, Blake Rablah, James Wilkie, Alex Mossmann, Mike
Donaldson, Wade Zawalski, Kelly Epp, Michel Laberge, Stephen Howard,
and Meritt Reynolds for useful discussions. We acknowledge the University
of Saskatchewan ICT Research Computing Facility for computing support. 

\appendix

\section{Appendix A\label{sec:Appendix-A}}

\section*{Detailed derivation of the first order node-to-element differential
operator\label{sec:detailed-derivation-of}}

In the following, the element-specific node indexes $i,\,j,\,k$ will
be replaced with the indexes $1,\:2$ and $3$ for simplicity. From
equation \ref{eq:500.02}, the radial and axial first spatial derivatives
of $U^{e}(\mathbf{r})$ are constants across the element $e$, and
are given by $\frac{\partial U^{e}}{\partial r}=B^{e}$, and $\frac{\partial U^{e}}{\partial z}=C^{e}$.
Expressions for $B^{e}$ and $C^{e}$ are found using equation \ref{eq:503-1}:
\begin{align*}
U_{1}^{e}-U_{3}^{e} & =B^{e}(r_{1}-r_{3})+C^{e}(z_{1}-z_{3})\\
U_{2}^{e}-U_{3}^{e} & =B^{e}(r_{2}-r_{3})+C^{e}(z_{2}-z_{3})
\end{align*}

\begin{alignat}{1}
\Rightarrow B^{e} & =\frac{\partial U^{e}}{\partial r}=\frac{U_{1}^{e}(z_{2}-z_{3})+U_{2}^{e}(z_{3}-z_{1})+U_{3}^{e}(z_{1}-z_{2})}{2s^{e}}\nonumber \\
C^{e} & =\frac{\partial U^{e}}{\partial z}=\frac{U_{1}^{e}(r_{2}-r_{3})+U_{2}^{e}(r_{3}-r_{1})+U_{3}^{e}(r_{1}-r_{2})}{-2s^{e}}\label{eq:503.2}
\end{alignat}
The triangle area $s^{e}$ is introduced here, assuming that vertices
are numbered counterclockwise, noting that $2s^{e}=|\mathbb{R}^{e}|$,
where $|\mathbb{R}^{e}|$ is the$\text{}$ determinant of array $\mathbb{R}^{e}=\left(\begin{array}{ccc}
1 & r_{1} & z_{1}\\
1 & r_{2} & z_{2}\\
1 & r_{3} & z_{3}
\end{array}\right)$, so that 
\[
2s^{e}=(r_{1}-r_{3})(z_{2}-z_{3})-(r_{2}-r_{3})(z_{1}-z_{3})
\]
The element-specific spatial derivatives can also be expressed in
terms of the coefficients in the element-specific pyramid-side functions
$\psi_{1}^{e},\,\psi_{2}^{e}$, and $\psi_{3}^{e}$. As noted earlier,
each pyramid-side function is specific to a particular node, and to
a particular element, and is defined as $\psi_{n}^{e}(\textbf{\ensuremath{\mathbf{r}}})=a_{n}^{e}+b_{n}^{e}r+c_{n}^{e}z$,
and has the property $\psi_{n}^{e}(r_{j},\,z_{j})=\delta_{nj}$. This
yields, for each element $e$, the identity $\mathbb{R}^{e}*\mathbb{C}^{e}=\mathbb{I}$,
or 
\[
\left(\begin{array}{ccc}
1 & r_{1} & z_{1}\\
1 & r_{2} & z_{2}\\
1 & r_{3} & z_{3}
\end{array}\right)*\left(\begin{array}{ccc}
a_{1}^{e} & a_{2}^{e} & a_{3}^{e}\\
b_{1}^{e} & b_{2}^{e} & b_{3}^{e}\\
c_{1}^{e} & c_{2}^{e} & c_{3}^{e}
\end{array}\right)=\left(\begin{array}{ccc}
1 & 0 & 0\\
0 & 1 & 0\\
0 & 0 & 1
\end{array}\right)
\]
 
\begin{align}
\Rightarrow\mathbb{C}^{e} & =\left(\mathbb{R}^{e}\right)^{-1}=\frac{1}{2s^{e}}\left(\begin{array}{ccc}
\left(r_{2}z_{3}-r_{3}z_{2}\right) & \left(r_{3}z_{1}-r_{1}z_{3}\right) & \left(r_{1}z_{2}-r_{2}z_{1}\right)\\
\left(z_{2}-z_{3}\right) & \left(z_{3}-z_{1}\right) & \left(z_{1}-z_{2}\right)\\
\left(r_{3}-r_{2}\right) & \left(r_{1}-r_{3}\right) & \left(r_{2}-r_{1}\right)
\end{array}\right)\label{eq:504.1}
\end{align}
Comparing equations \ref{eq:503.2} with equation \ref{eq:504.1},
it is evident that 

\begin{align}
\frac{\partial U^{e}}{\partial r} & =\underset{\varepsilon}{\Sigma}U_{\varepsilon}^{e}b_{\varepsilon}^{e}\nonumber \\
\frac{\partial U^{e}}{\partial z} & =\underset{\varepsilon}{\Sigma}U_{\varepsilon}^{e}c_{\varepsilon}^{e}\label{eq:505}
\end{align}
Given the values of the approximation for the piecewise linear function
$U(\mathbf{r})$ on the triangle vertices, and the element-specific
array $\mathbb{C}^{e}$ evaluated using equation \ref{eq:504.1},
equations \ref{eq:505} can be used to determine the values of the
spatial derivatives of $U(\mathbf{r})$ at the interior of each element. 

\section{Appendix B\label{sec:Appendix-B}}

\section*{Node-to-element operator identities\label{sec:IDproofANGMOM}}

We demonstrate here a property of the node-to-element operators that
is used to demonstrate angular momentum conservation in section \ref{subsec:Angular-momentum-conservation}.
For any continuous scalar functions $f$ and $g$, the identity (Stoke's
theorem) holds that 
\begin{equation}
\int\left(\nabla f\times\nabla g\right)\cdot d\mathbf{S}=\int\nabla\times\left(f\,\nabla g\right)\cdot d\mathbf{S}=\int f\,\nabla g\cdot d\mathbf{l}\label{eq:505.001}
\end{equation}
With azimuthal symmetry, $\int\left(\nabla f\times\nabla g\right)\cdot d\mathbf{S}=\int\left(f_{z}'\,g_{r}'-f_{r}'\,g_{z}'\right)\widehat{\boldsymbol{\phi}}\cdot d\mathbf{S}=\int\left(f_{z}'\,g_{r}'-f_{r}'\,g_{z}'\right)ds_{\phi}$,
where $ds_{\phi}$ is an elemental area in the poloidal ($r-z$) plane
and $\widehat{\boldsymbol{\phi}}$ is the unit vector in the azimuthal
direction. Hence, the discrete form of $\int\left(\nabla f\times\nabla g\right)\cdot d\mathbf{S}$,
using the node-to-element differential operators, is 
\begin{alignat}{1}
\left(\int\left(\nabla f\times\nabla g\right)\cdot d\mathbf{S}\right)_{disc.} & =\widehat{s}^{T}*\left(\left(\overline{\widehat{Dz}}*\overline{f}\right)\,\,\left(\overline{\widehat{Dr}}*\overline{g}\right)-\left(\overline{\widehat{Dr}}*\overline{f}\right)\,\,\left(\overline{\widehat{Dz}}*\overline{g}\right)\right)\nonumber \\
 & =\left(\overline{\widehat{Dz}}*\overline{f}\right)^{T}*\widehat{\widehat{S}}*\left(\overline{\widehat{Dr}}*\overline{g}\right)-\left(\overline{\widehat{Dr}}*\overline{f}\right)^{T}*\widehat{\widehat{S}}*\left(\overline{\widehat{Dz}}*\overline{g}\right)\nonumber \\
 & =\overline{f}^{T}*\overline{\widehat{Dz}}^{T}*\widehat{\widehat{S}}*\overline{\widehat{Dr}}*\overline{g}-\overline{f}^{T}*\overline{\widehat{Dr}}^{T}*\widehat{\widehat{S}}*\overline{\widehat{Dz}}*\overline{g}\nonumber \\
 & =\overline{f}^{T}*\left(\overline{\widehat{Dz}}^{T}*\widehat{\widehat{S}}*\overline{\widehat{Dr}}-\overline{\widehat{Dr}}^{T}*\widehat{\widehat{S}}*\overline{\widehat{Dz}}\right)*\overline{g}\nonumber \\
 & =\overline{f}^{T}*\left(\overline{\overline{B}}^{T}-\overline{\overline{B}}\right)*\overline{g}\nonumber \\
 & =\overline{f}^{T}*\overline{\overline{A}}*\overline{g}\label{eq:505.0011}
\end{alignat}
where $\overline{\overline{B}}=\overline{\widehat{Dr}}^{T}*\widehat{\widehat{S}}*\overline{\widehat{Dz}}$.
Setting $\overline{f}=\overline{\phi}_{m}$ and $\overline{g}=\overline{\phi}_{n}$,
where $\overline{\phi}_{m}$ and $\overline{\phi}_{n}$ are the vectors
defining the nodal values of basis functions $\phi_{m}$ and $\phi_{n}$
respectively ($i.e.,$ $\overline{\phi}_{m}(k)=\delta_{mk}$ and $\overline{\phi}_{n}(k)=\delta_{nk}$),
equations \ref{eq:505.0011} and \ref{eq:505.001} imply that 
\begin{equation}
A_{mn}\equiv A(m,n)=\overline{\phi}_{m}^{T}*\overline{\overline{A}}*\overline{\phi}_{n}=\int\phi_{m}\,\nabla\phi_{n}\cdot d\mathbf{l}\label{eq:505.002}
\end{equation}
We will show that each element of $\overline{\overline{A}}$ is zero
at internal nodes, $i.e.,$ $\left(\overline{\overline{A}}\right)_{int}=0$.
There are four cases to consider: (1) $m\neq n$ and nodes $m$ and
$n$ are not adjacent, (2) $m=n$, (3) $m$ and $n$ are the indexes
of adjacent internal nodes, and (4) $m$ and $n$ are the indexes
of adjacent nodes located on the boundary of the computational domain.
In case (1), $A_{mn}$ is obviously zero. In case (2), $A_{mn}=\overline{\phi}_{m}^{T}*\left(\overline{\overline{B}}^{T}-\overline{\overline{B}}\right)*\overline{\phi}_{m}=\overline{\phi}_{m}^{T}*\overline{\overline{B}}^{T}*\overline{\phi}_{m}-\overline{\phi}_{m}^{T}*\overline{\overline{B}}*\overline{\phi}_{m}$
. The first term scalar term here can be transposed, so that $A_{mn}=0$
in case (2). 
\begin{figure}[H]
\centering{}\includegraphics[width=8cm,height=6cm]{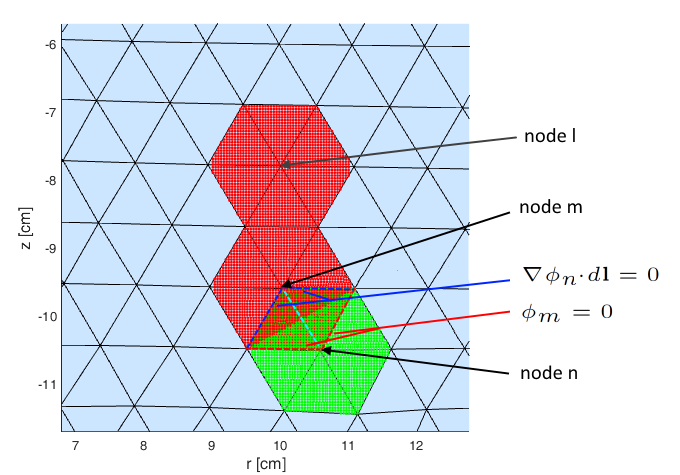}\caption{\label{fig:Linear-basis-function_topview}Linear basis function depiction
for triangular elements (top view). Note the dark blue dashed lines
represent the part of the boundary of basis function $\phi_{n}$ that
overlaps with $\phi_{m}$, and the red dashed lines represent the
part of the boundary of $\phi_{m}$ that overlaps with $\phi_{n}$.}
\end{figure}
Figure \ref{fig:Linear-basis-function_topview} is a top-view of figure
\ref{fig:Linear-basis-function}. Referring to equation \ref{eq:505.002}
and to figure \ref{fig:Linear-basis-function_topview}, the only contribution
to $A_{mn}$ in case (3) is from the contour integral $\int\phi_{m}\,\nabla\phi_{n}\cdot d\mathbf{l}$
along the boundary (depicted with dark blue and red dashed lines)
of the diamond-shape representing the overlapping area of the basis
functions $\phi_{m}$ and $\phi_{n}$. It can be seen that $\phi_{m}=0$
along part of the contour, and $\nabla\phi_{n}\cdot d\mathbf{l}=0$
along the remaining parts, so that $A_{mn}=0$ in case (3). In case
(4), when $m$ and $n$ are the indexes of adjacent nodes located
on the boundary of the computational domain, the integral $\int\phi_{m}\,\nabla\phi_{n}\cdot d\mathbf{l}$
would be finite along the light blue dashed line connecting nodes
$m$ and $n$ (in case (4), the light blue dashed line would represent
part of the computational domain boundary). Therefore, $A_{mn}$ is
finite only when $m$ and $n$ are the indexes of adjacent nodes located
on the boundary of the computational domain, leading to the identity
\begin{equation}
\left(\overline{\widehat{Dz}}^{T}*\widehat{\widehat{S}}*\overline{\widehat{Dr}}-\overline{\widehat{Dr}}^{T}*\widehat{\widehat{S}}*\overline{\widehat{Dz}}\right)_{int}=0\label{eq:505.003}
\end{equation}

\section{Appendix C\label{sec:Appendix-C}}

\section*{Node-to-node operator identities\label{sec:ID_proff_n2n}}

Now we introduce some special mimetic properties of the node-to-node
differential operators. $U_{r}'(\mathbf{r})$, the radial spatial
derivative of the approximate solution $U(\mathbf{r})$, can alternatively
be expressed, using the basis function expansion (equation \ref{eq:500.03}),
as 
\begin{equation}
U_{r}'(\mathbf{r})=\overset{N_{n}}{\underset{n=1}{\Sigma}}U{}_{n}\,\phi_{rn}'(\mathbf{r})\label{eq:506}
\end{equation}
Applying the Galerkin method, we obtain 
\begin{equation}
\int\phi{}_{j}\,U_{r}'(\mathbf{r})\,dr\,dz=\overset{N_{n}}{\underset{n=1}{\Sigma}}U{}_{n}\int\phi{}_{j}\,\phi_{rn}'\,dr\,dz\label{eq:507}
\end{equation}
Again, from equation \ref{eq:502.4}, the left side of this equation
is $U_{rj}'\,s_{j}/3$. On the right side, matrix $\overline{\overline{Lr}}$
of size $N_{n}\times N_{n}$ is introduced as 
\begin{equation}
\frac{1}{3}\overline{\overline{Lr}}(j,\,n)=\int\phi{}_{j}\,\phi_{rn}'\,dr\,dz\label{eq:508.1}
\end{equation}
Hence, equation \ref{eq:507} can be re-expressed as 

\begin{equation}
\overline{U}_{r}'(j)\,\overline{s}(j)=\overset{N_{n}}{\underset{n=1}{\Sigma}}\,\left(\overline{\overline{Lr}}(j,\,n)\,\overline{U}(n)\right)\label{eq:509}
\end{equation}
or in matrix form: $\overline{\overline{S}}*\overline{\overline{Dr}}*\overline{U}=\overline{\overline{Lr}}*\overline{U}$
\begin{equation}
\Rightarrow\overline{\overline{Dr}}=\overline{\overline{S}}^{-1}*\overline{\overline{Lr}}\label{eq:509.1}
\end{equation}
With integration by parts and Green's theorem, equation \ref{eq:508.1}
implies that 
\begin{align*}
\overline{\overline{Lr}}(j,\,n) & =3\int\phi{}_{j}\,\phi_{rn}'\,dr\,dz\\
 & =3\int\left(-\phi_{rj}'\phi{}_{n}\,+\left(\phi{}_{j}\,\phi{}_{n}\right)_{r}'\right)\,dr\,dz\\
 & =3\left(-\int\phi{}_{n}\,\phi_{rj}'\,dr\,dz+\varoint\phi{}_{j}\,\phi{}_{n}\,dz\right)
\end{align*}
The boundary term vanishes at internal nodes, so that $\overline{\overline{Lr}}(j,\,n)=-\overline{\overline{Lr}}(n,\,j)$
at the internal nodes. This implies that all elements of the matrix
$\overline{\overline{Lr}}+\overline{\overline{Lr}}^{T}$ are equal
to zero, except for the elements corresponding to boundary nodes.
Hence, using equation \ref{eq:509.1}, and noting that the transpose
of a diagonal matrix is the same matrix, we obtain 
\begin{equation}
\left(\overline{\overline{S}}*\overline{\overline{Dr}}+\overline{\overline{Dr}}^{T}*\overline{\overline{S}}\right)_{int}=0\label{eq:510.23}
\end{equation}
where the subscript $int$ denotes all elements of the matrix that
correspond to interior (non-boundary) nodes. There is an analogous
identity for the $\overline{\overline{Dz}}$ operator:
\begin{equation}
\left(\overline{\overline{S}}*\overline{\overline{Dz}}+\overline{\overline{Dz}}^{T}*\overline{\overline{S}}\right)_{int}=0\label{eq:510.23-1}
\end{equation}
As a consequence of these identities, the node-to-node differential
operators have some important properties that mimic properties of
the corresponding continuous operators. In particular, the discrete
expression for the volume integral in equation \ref{eq:500.0} is
\begin{align*}
 & \overline{dV}{}^{T}*\left[\overline{U}\,\,\left(\overline{\overline{\nabla}}\cdot\overline{\mathbf{P}}\right)+\overline{\mathbf{P}}\cdot\left(\overline{\overline{\nabla}}\,\,\overline{U}\right)\right]\\
 & =\overline{dV}{}^{T}*\left[\overline{U}\,\,\left(\overline{\overline{Dr}}*\left(\overline{r}\,\,\overline{P}_{r}\right)+\overline{\overline{Dz}}*\left(\overline{r}\,\,\overline{P}_{z}\right)\right)\,/\,\overline{r}+\overline{P}_{r}\,\,\left(\overline{\overline{Dr}}*\overline{U}\right)+\overline{P}_{z}\,\,\left(\overline{\overline{Dz}}*\overline{U}\right)\right]\\
 & =\frac{2\pi}{3}\left[\overline{U}^{T}*\overline{\overline{S}}*\left(\overline{\overline{Dr}}*\overline{\overline{R}}*\overline{P}_{r}+\overline{\overline{Dz}}*\overline{\overline{R}}*\overline{P}_{z}\right)+\overline{P}_{r}{}^{T}*\overline{\overline{S}}*\overline{\overline{R}}*\overline{\overline{Dr}}*\overline{U}+\overline{P}_{z}{}^{T}*\overline{\overline{S}}*\overline{\overline{R}}*\overline{\overline{Dz}}*\overline{U}\right]\\
 & =\frac{2\pi}{3}\,\overline{U}{}^{T}*\left[\left(\overline{\overline{S}}*\overline{\overline{Dr}}+\overline{\overline{Dr}}{}^{T}*\overline{\overline{S}}\right)*\overline{\overline{R}}*\overline{P}_{r}+\left(\overline{\overline{S}}*\overline{\overline{Dz}}+\overline{\overline{Dz}}{}^{T}*\overline{\overline{S}}\right)*\overline{\overline{R}}*\overline{P}_{z}\right]
\end{align*}
Note that in the last step, transpose operations on the scalar expressions
have been performed. Hence, due to equations \ref{eq:510.23} and
\ref{eq:510.23-1}
\begin{align}
\overline{dV}{}^{T}*\left[\left(\overline{U}\,\,\left(\overline{\overline{\nabla}}\cdot\overline{\mathbf{P}}\right)+\left(\overline{\overline{\nabla}}\,\,\overline{U}\right)\cdot\overline{\mathbf{P}}\right)\right] & =0 & \left(\mbox{if }\ensuremath{\ensuremath{\overline{U}\,|_{\Gamma}}=0}\mbox{ or }\ensuremath{\overline{\mathbf{P}}_{\perp}|_{\Gamma}=0}\right)\label{eq:511.04}
\end{align}
For the particular case of $\overline{U}=\overline{1}$ , this implies
that 
\begin{align}
\overline{dV}{}^{T}*\left[\overline{\overline{\nabla}}\cdot\overline{\mathbf{P}}\right] & =\overline{0} & \left(\mbox{if }\ensuremath{\overline{\mathbf{P}}_{\perp}|_{\Gamma}=0}\right)\label{eq:511.041}
\end{align}

\section{Appendix D\label{sec:Appendix-D}}

\section*{Maintenance of energy and momentum conservation with artificial density
diffusion\label{subsec:Maintenance-of-momentum}}

The components of force per volume vector $\overline{\mathbf{f}}_{\zeta}=(\overline{f}_{\zeta r},\,\overline{f}_{\zeta\phi},\,\overline{f}_{\zeta z})^{T}$
in the velocity equations, and in the heating term $\overline{Q}_{\zeta}$
in the ion pressure equation (equation \ref{eq:517.3}) are included
to cancel the effect of artificial density diffusion on the total
system momentum and energy, maintaining conservation. There are two
models for the correction terms. The first model is straightforward
and maintains conservation of energy and angular momentum when density
diffusion is included in the mass continuity equation. It is evaluated
simply by treating the density diffusion as a particle source and
assessing the local ($i.e.,$ nodal) correction terms as the local
effects of the source on momentum and energy. However this model can
lead locally to negative ion pressure if the density gradients are
extreme, so it is not suitable for inclusion in simulations of CT
formation and compression. The model works ok for other simulation
types, for example the resistive decay of a magnetised plasma described
by an equilibrium model. The second model is more complicated - it
maintains conservation of volume integrated energy and can not cause
negative ion pressure, but angular momentum conservation is not maintained.
Both models for the correction terms compensate for additional $r$
and $z$ directed momentum introduced due to density diffusion - for
the first model the compensation is local, while for the second model
the compensation is in the volume integrated sense. The simulations
presented in section \ref{sec:Simulation-results,-and} of this paper
use the second model. 

The correction terms for the first model are derived as follows. For
convenience, we use the notation $\overline{\zeta}_{n}=\widehat{\overline{\nabla}}\cdot\left(\widehat{\mathbf{\zeta}}\,\,\left(\overline{\widehat{\nabla}}\,\,\overline{n}\right)\right)$.
Referring to equation \ref{eq:517.3}, note that $\dot{\overline{n}}|_{\zeta}=\overline{\zeta}_{n}$,
$\dot{\overline{v}_{\beta}}|_{\zeta}=\overline{f}_{\zeta\beta}\,/\,\overline{\rho}$,
and $\dot{\overline{p}}_{i}|_{\zeta}=\overline{Q}_{\zeta}$ are the
parts $\dot{\overline{n}}$, $\dot{\overline{v}_{\beta}}$, and $\dot{\overline{p}}_{i}$
respectively that are associated with density diffusion. The local
rate of change of the $\beta$ component (where, here, $\beta=r,\,z$)
of momentum per unit volume due to the local particle source/sink
terms arising due to density diffusion is set to zero by design, resulting
in 
\begin{align*}
\left(\frac{\partial}{\partial t}\left(m_{i}\,\overline{n}\,\,\overline{v}_{\beta}\right)\right)_{\zeta} & =m_{i}\,\left(\overline{n}\,\,\overline{f}_{\zeta\beta}\,/\,\overline{\rho}+\overline{\zeta}_{n}\,\,\overline{v}_{\beta}\right)=0\\
 & \Rightarrow\overline{f}_{\zeta\beta}=-m_{i}\,\overline{v}_{\beta}\,\,\overline{\zeta}_{n}
\end{align*}
Naturally, for $\beta=\phi$, this expression for $\overline{f}_{\zeta\beta}$
also locally cancels additional angular momentum per unit volume introduced
by density diffusion. The rate of change of energy per unit volume
due to local particle source/sink terms associated with density diffusion
is also set to zero: 
\begin{align*}
\left(\frac{\partial}{\partial t}\left(\frac{1}{2}m_{i}\,\overline{n}\,\,\overline{v}^{2}+\frac{1}{\gamma-1}\,\dot{\overline{p}}_{i}\right)\right)_{\zeta} & =\frac{1}{2}m_{i}\,\,\overline{\zeta}_{n}\,\,\overline{v}^{2}+\underset{\beta}{\Sigma}\left(m_{i}\,\overline{n}\,\,\overline{v}_{\beta}\,\,\overline{f}_{\zeta\beta}\,/\,\overline{\rho}\right)+\overline{Q}_{\zeta}=0\\
 & \Rightarrow\overline{Q}_{\zeta}=\frac{1}{2}m_{i}\,\,\overline{v}^{2}\,\,\overline{\zeta}_{n}
\end{align*}
Hence, the contributions to $\dot{U}_{K}$ and $\dot{U}_{Th}$ arising
from the artificial density diffusion and relevant correction terms
cancel by design as:

{\small{}
\begin{align*}
\left(\dot{U}_{K}+\dot{U}_{Th}\right)_{\zeta} & =\overline{dV}^{T}*\left\{ m_{i}\left(\frac{1}{2}-1+\frac{1}{2}\right)\,\overline{v}^{2}\,\,\overline{\zeta}_{n}\right\} =0
\end{align*}
}However, note that if density gradients are extreme, $\overline{\zeta}_{n}$
and hence $\overline{Q}_{\zeta}$ will have some very large negative
values, so that $\overline{p}_{i}$ can become negative.\\

For the second model of correction terms, $\overline{Q}_{\zeta}=\overline{0}$,
and a suitable choice of $\overline{f}_{\zeta\beta}$ that ensures
maintenance, in a volume integrated sense, of conservation of total
energy, is 
\begin{equation}
\overline{f}_{\zeta\beta}=\frac{1}{2}m_{i}\,\zeta\left[\widehat{\overline{W}}*\left(\left(\overline{\widehat{\nabla}}\,\,\overline{n}\right)\cdot\left(\overline{\widehat{\nabla}}\,\,\overline{v}_{\beta}\right)\right)+\widehat{\overline{\nabla}}\cdot\left(\widehat{v}_{\beta}\,\,\left(\overline{\widehat{\nabla}}\,\,\overline{n}\right)\right)-\overline{v}_{\beta}\,\,\left(\overline{\overline{\Delta}}\,\,\overline{n}\right)\right]\label{eq:518-1}
\end{equation}
Note that while this term also cancels modifications to the $r$ and
$z$ directed momentum resulting from density diffusion, modification
to angular momentum is not compensated for. This method requires that
$\zeta$ is spatially constant so that, in this case $\overline{\zeta}_{n}=\zeta\,\overline{\overline{\Delta}}\,\,\overline{n}$.
Expressions for $\overline{f}_{\zeta\beta}$ that also conserve angular
momentum can be derived if $\zeta$ is made to be spatially dependent,
but this additional complication has been omitted for now.

Using expression in equation \ref{eq:518-1} for $\overline{f}_{\zeta\beta}$,
the total contribution to the volume integral of the $\beta$ component
of momentum per unit volume (where $\beta=r,\,z$) from terms associated
with density diffusion is {\scriptsize{}
\begin{align*}
\dot{P}_{\beta}|_{\zeta} & =\overline{dV}^{T}*\left\{ \left(\frac{\partial}{\partial t}\left(m_{i}\,\overline{n}\,\,\overline{v}_{\beta}\right)\right)_{\zeta}\right\} \\
 & =\overline{dV}^{T}*\left\{ m_{i}\,\overline{\zeta}_{n}\,\,\overline{v}_{\beta}+\overline{f}_{\zeta\beta}\right\} \\
 & =\overline{dV}^{T}*\left\{ m_{i}\,\zeta\,\left(\overline{\overline{\Delta}}\,\,\overline{n}\right)\,\,\overline{v}_{\beta}+\frac{1}{2}m_{i}\,\zeta\left[\widehat{\overline{W}}*\left(\left(\overline{\widehat{\nabla}}\,\,\overline{n}\right)\cdot\left(\overline{\widehat{\nabla}}\,\,\overline{v}_{\beta}\right)\right)+\cancel{\widehat{\overline{\nabla}}\cdot\left(\widehat{v}_{\beta}\,\,\left(\overline{\widehat{\nabla}}\,\,\overline{n}\right)\right)}-\overline{v}_{\beta}\,\,\left(\overline{\overline{\Delta}}\,\,\overline{n}\right)\right]\right\}  & \mbox{(use eqn. \ref{eq:515.04})}\\
 & =\frac{1}{2}m_{i}\,\zeta\,\overline{dV}^{T}*\left\{ \widehat{\overline{W}}*\left(\left(\overline{\widehat{\nabla}}\,\,\overline{n}\right)\cdot\left(\overline{\widehat{\nabla}}\,\,\overline{v}_{\beta}\right)\right)+\overline{v}_{\beta}\,\,\left(\overline{\overline{\Delta}}\,\,\overline{n}\right)\right\} \\
 & =\frac{1}{2}m_{i}\,\zeta\,\left[\widehat{dV}^{T}*\left\{ \left(\overline{\widehat{\nabla}}\,\,\overline{n}\right)\cdot\left(\overline{\widehat{\nabla}}\,\,\overline{v}_{\beta}\right)\right\} +\overline{dV}^{T}*\left\{ \overline{v}_{\beta}\,\,\left(\widehat{\overline{\nabla}}\cdot\left(\overline{\widehat{\nabla}}\,\,\overline{n}\right)\right)\right\} \right] & \mbox{(use eqn. \ref{eq:516.1})}\\
\Rightarrow\dot{P}_{\beta}|_{\zeta} & =0 & (\mbox{use eqn. }\ref{eq:515.031})
\end{align*}
}Note that in the last step, identity \ref{eq:515.031} has been used,
with $\widehat{\mathbf{P}}=\overline{\widehat{\nabla}}\,\,\overline{n}$
and $\overline{U}=\overline{v}_{\beta}$. 

The total contribution to energy associated with $v_{\beta}$ (where
$\beta=r,\,\phi,\,z$) from terms associated with density diffusion
is {\footnotesize{}
\begin{align*}
\dot{U}_{\beta}|_{\zeta} & =\overline{dV}^{T}*\left\{ \left(\frac{\partial}{\partial t}\left(\frac{1}{2}\,m_{i}\,\overline{n}\,\,\overline{v}_{\beta}^{2}\right)\right)_{\zeta}\right\} \\
 & =\overline{dV}^{T}*\left\{ \frac{1}{2}m_{i}\,\overline{\zeta}_{n}\,\,\overline{v}_{\beta}^{2}+\overline{v}_{\beta}\,\,\overline{f}_{\zeta\beta}\right\} \\
 & =\overline{dV}^{T}*\left\{ \cancel{\frac{1}{2}m_{i}\,\zeta\,\left(\overline{\overline{\Delta}}\,\,\overline{n}\right)\,\,\overline{v}_{\beta}^{2}}+\frac{1}{2}m_{i}\,\zeta\,\overline{v}_{\beta}\,\,\left[\widehat{\overline{W}}*\left(\left(\overline{\widehat{\nabla}}\,\,\overline{n}\right)\cdot\left(\overline{\widehat{\nabla}}\,\,\overline{v}_{\beta}\right)\right)+\widehat{\overline{\nabla}}\cdot\left(\widehat{v}_{\beta}\,\,\left(\overline{\widehat{\nabla}}\,\,\overline{n}\right)\right)\cancel{-\overline{v}_{\beta}\,\,\left(\overline{\overline{\Delta}}\,\,\overline{n}\right)}\right]\right\} \\
 & =\frac{1}{2}m_{i}\,\zeta\,\overline{dV}^{T}*\left\{ \overline{v}_{\beta}\,\,\left[\widehat{\overline{W}}*\left(\left(\overline{\widehat{\nabla}}\,\,\overline{n}\right)\cdot\left(\overline{\widehat{\nabla}}\,\,\overline{v}_{\beta}\right)\right)+\widehat{\overline{\nabla}}\cdot\left(\widehat{v}_{\beta}\,\,\left(\overline{\widehat{\nabla}}\,\,\overline{n}\right)\right)\right]\right\} \\
 & =\frac{1}{2}m_{i}\,\zeta\left[\widehat{dV}^{T}*\left\{ \widehat{v}_{\beta}\,\,\left(\left(\overline{\widehat{\nabla}}\,\,\overline{n}\right)\cdot\left(\overline{\widehat{\nabla}}\,\,\overline{v}_{\beta}\right)\right)\right\} +\overline{dV}^{T}*\left\{ \overline{v}_{\beta}\,\,\left(\widehat{\overline{\nabla}}\cdot\left(\widehat{v}_{\beta}\,\,\left(\overline{\widehat{\nabla}}\,\,\overline{n}\right)\right)\right)\right\} \right] & \mbox{(use eqn. \ref{eq:516})}\\
\Rightarrow\dot{U}_{\beta}|_{\zeta} & =0 & (\mbox{use eqn. }\ref{eq:515.031})
\end{align*}
}In the last step, identity \ref{eq:515.031} has been used again,
now with $\widehat{\mathbf{P}}=\widehat{v}_{\beta}\,\,\left(\overline{\widehat{\nabla}}\,\,\overline{n}\right)$
and $\overline{U}=\overline{v}_{\beta}$. 

Hence, the modification to total system energy due to density diffusion
is canceled with the inclusion of the correction terms using either
correction model, thereby ensuring maintenance of global energy conservation.
Conservation of angular momentum is maintained using the first model
only, but that model is suitable only for simulations that don't involve
steep density gradients. Both models compensate for the modifications
to radially and axially directed momentum associated with density
diffusion.

\end{document}